\documentclass[12pt]{amsart}

\usepackage{latexsym,amsfonts,amsthm,amsmath,amssymb,amscd}

\newtheorem{theorem}{Theorem}[section]

\newtheorem{corollary}[theorem]{Corollary}

\newtheorem{lemma}[theorem]{Lemma}

\newtheorem{proposition}[theorem]{Proposition}

{\theoremstyle{definition}

        \newtheorem{definition}[theorem]{Definition}
        \newtheorem{remark}[theorem]{Remark}
        \newtheorem{example}[theorem]{Example}
        \newtheorem{parrafo}[theorem]{{\!}}  }

\numberwithin{equation}{theorem}

\newcommand{\nat}{\mathbb N} \newcommand{\cali}{{\mathcal {I}}}
\newcommand{\calo}{{\mathcal {O}}}

\newcommand{\Fnu}{\boldsymbol{\nu}}  \newcommand{\Oh}{{\mathcal O}}
\newcommand{\calB}{{\mathcal B}} \newcommand{\call}{{\mathcal L}}

 \DeclareMathOperator{\codim}{codim} \DeclareMathOperator{\Coeff}{Coeff}
 \DeclareMathOperator{\Max}{\underline{Max}}
 \DeclareMathOperator{\ord}{ord} \DeclareMathOperator{\Proj}{Proj}
\DeclareMathOperator{\Reg}{Reg} \DeclareMathOperator{\Sing}{Sing} \DeclareMathOperator{\Spec}{Spec}
\DeclareMathOperator{\word}{w-ord}

\def\O{{\mathcal{O}}}

\title[]{A  Simplified Proof of Desingularization and  Applications}

\author{A. Bravo \and S. Encinas \and O. Villamayor U.}

\address{Department of Mathematics, University of Michigan, Ann Arbor,
MI 48109-1109, U.S.A.} \email{anabz@math.lsa.umich.edu}

\address{Current address: Dpto. Matem\'aticas, Facultad de Ciencias, Universidad
Aut\'onoma de Madrid, Canto Blanco 28049 Madrid, Spain.}
\email{ana.bravo@uam.es}

\address{Departamento de Matem\'atica Aplicada Fundamental Universidad
de Valladolid, Spain.} \email{sencinas@maf.uva.es}

\address{Dpto. Matem\'aticas, Facultad de Ciencias, Universidad
Aut\'onoma de Madrid, Canto Blanco 28049 Madrid, Spain.} \email{villamayor@uam.es}
\thanks{All the authors were partially supported by PB96-0065.
The first author was also partially supported by  a Rackham Fellowship during the summer in 2001.}


 \date{} \dedicatory{} \commby{}

\begin{document}
\maketitle
\begin{abstract}
This paper contains a short and simplified  proof of desingularization
    over fields  of characteristic zero, together with various
    applications to other problems in
algebraic geometry  (among others, the study of the behavior of
    desingularization of families of embedded schemes,  and  a formulation of
    desingularization
    which is stronger than Hironaka's).
    Our   proof avoids the use of the
    Hilbert-Samuel function and Hironaka's notion of normal flatness:
    First we define a procedure for principalization of ideals (i.~e. a
    procedure to make an ideal invertible), and then we show that
    desingularization of a closed subscheme \( X \) is achieved by
    using the procedure of principalization for the ideal \( {\mathcal
    I}(X) \) associated to the embedded scheme \( X \). The paper
    intends to be an
    introduction to the subject, focused on the motivation of ideas
    used in this new approach, and particularly on applications, some of which
    do not follow from Hironaka's proof.
\end{abstract}

\tableofcontents






\part{Introduction.}

\label{introduction}




\section{Introduction.}
\label{introduction1} We present  a unified proof, for the case of schemes over a field of characteristic
zero and of compact analytic spaces, of two central theorems in algebraic geometry: The {\em resolution
of singularities of schemes}, and the  {\em principalization of ideals} (see Section \ref{lang} for
precise statements of  our results). The existence of resolution of singularities is one of the most
important results in the area, due to its large number of applications,   while principalization of
ideals is related to the classical problem of elimination of base points of a linear system.
\medskip

\noindent{\bf Resolution of singularities of schemes.} {\em Given a reduced scheme $X$ of finite type
over a field of characteristic zero, find a proper and birational morphism $$\pi: X^{\prime}\to X,$$ such
that $X^{\prime}$ is non-singular and $\pi$ induces an isomorphism over the regular points of $X$.}
\medskip

\noindent{\bf Principalization.} {\em Given a non zero sheaf of ideals \(J\) in  non-singular variety
$W$, find  a proper and birational morphism, $$\pi: W^{\prime} \to W,$$ such that $W^{\prime}$ is
non-singular and $J{\mathcal O}_{W^{\prime}}$ is locally principal.}  \medskip

The purpose of this paper is threefold: We present  a simplified proof of both  theorems, already
indicated in \cite{EncVil97:Tirol}, see also \cite{EncVil99}, (this is done in Parts \ref{partbo} and
\ref{appendix} of this manuscript); we summarize some applications of our  approach, together with  new
results  (see  Part \ref{applications}); and finally, in Part \ref{nature}, we motivate the main ideas
behind the invariants that we use in this proof of desingularization (this part is purely expository).
Our very short proof of desingularization presented in \ref{existence}, makes use of the algorithm
developed essentially in Part \ref{appendix} (and recently implemented in MAPLE). This provides a proof
of desingularization that does not exceed thirty pages.

\

\begin{center}
{\em What makes this work different from other approaches to desingularization? }
\end{center}

\

In his  monumental work (\cite{Hir64}), Hironaka gave an existential proof of  resolution of
singularities, and,  in so doing, he also provided an existential proof for principalization of ideals.
Hironaka's  approach is based in two main ideas: Reduction to the hypersurface case, and induction on the
dimension. As  for  the reduction to the hypersurface case, he shows that given a variety $X$, it is
possible to describe it as an intersection of a finite number of hypersurfaces $\{{H_i}\}_{i\in I}$,  in
such a way that the intersection of the {\em worst singularities} of the hypersurfaces $\{H_i\}_{i\in I}$
is the set of the {\em worst singularities} of $X$. Then he proves, using an existential argument, that
there is a sequence of blowing-ups such that  the singularities of the hypersurfaces are better, and
hence  the singularities of $X$ are better too.
\medskip

To find the hypersurfaces $H_i$, Hironaka  considers a
stratification of $X$ by means of the Hilbert-Samuel function. The
worst singularities of $X$   correspond to the stratum of points
where this function  is maximum. He then  provides an existential
argument to show that by blowing-up at smooth centers contained in
the maximum stratum of the  Hilbert-Samuel  function, the maximum
drops (the containment condition over the centers is usually
referred as {\em normal flatness}). Once this maximum drops, new
hypersurfaces ${H_i}$ have to be chosen. All this requires the use
of   {\em standard bases}, algorithms of division and the notion
of {\em strict transform of ideals} (see \cite{Lipman},
\cite{Oda1987}, \cite{Gir}, for more information on Hironaka's
line of proof, and also Section \ref{weakvsstrict}, where the
notion of strict transform of an ideal is defined).  \medskip

During the eighties and nineties, the first constructive proofs of
resolution of singularities appeared: \cite{BM97}, \cite{Villa89}
(and \cite{Villa92}),  \cite{EncVil98}  (see also \cite{BM02}).
However, all of them were based on Hironaka's approach. These
proofs provided a general algorithm of desingularization
indicating where to blow-up in order to eliminate the
singularities in a step by step procedure. The idea is to define
invariants of singular points, and to show that these invariants
improve when blowing up the set of worst points. Embedded
desingularization is then achieved by repeatedly blowing up the
set of worst points. All algorithmic procedures mentioned above,
make use, as Hironaka did in his original work, of the Hilbert
Samuel function and the strict transform of ideals. From an
algorithmic point of view, the notion of strict transform of an
ideal is quite complicated, except in the case of a hypersurface.
\medskip

What makes this present proof different from the previous ones is
that we show that embedded desingularization can be achieved
avoiding the use of Hilbert Samuel functions and the control of
strict transforms of ideals. Instead of using a standard basis we
take any set of generators of an ideal, and we do not need to
change these generators in the procedure of desingularization.

\

Our strategy is based in the fact that embedded desingularization turns out to be a straightforward
corollary of a simpler result: The {\em Algorithmic Resolution of Basic Objects}. While related to
Hironaka's idealistic exponents (see \cite{Hironaka77}), the notion of basic object is different and
easier to work with.  Our approach follows from the study of local properties of the
    algorithmic resolution process for basic objects developed in
    \cite{Villa92} and  \cite{EncVil97:Tirol} and also included in
    this paper.
In fact the  {\it algorithm of desingularization} in
    \cite{Villa92} was also obtained  by a
    suitable application of this same algorithm of resolution of basic
    objects but in a different manner, and our proof of equivariance
(compatibility with group actions) is also similar to the proof in
that paper. Furthermore, our new proof coincides with that from
    \cite{Villa92} in the hypersurface case, so explicit examples for
    algorithmic desingularization
  of hypersurfaces  can be found in \cite[Section 8]{Villa92} and in
    \cite{Villa02} (Hilbert-Samuel functions were not required in
Hironaka's proof
    in the hypersurface case).

\

Among other applications, this algorithm of resolution of basic objects has been implemented as a
computer program by G. Bodn\'ar and J. Schicho. As a consequence these authors have developed a program
which  resolves singularities of hypersurfaces (see \cite{GabSch98} and \cite{GabSch99}). The proof of
resolution of singularities that we give in \ref{existence}  shows that that this implemented algorithm
of resolution of basic objects, as it stands, also applies for embedded desingularization of any scheme,
even if not a hypersurface (see Section \ref{Computer}).
\medskip

This new approach to  resolution of singularities also allows us
to go beyond Hironaka's desingularization: Based on the ideas
developed here, in \cite{BV1} we present  a stronger formulation
of desingularization (se Section \ref{weakvsstrict} for a precise
statement of this result, as well as additional applications). See
also Part \ref{applications} for more applications.
\medskip

Finally, there are other results in desingularization within this
new framework and  related to this one: For instance,  Matsuki's
lucid presentation of the ideas behind the algorithms of
desingularization in \cite{Matsuki}, which is based in  the work
developed in \cite{EncVil99}; or \cite{EH}, where an algorithm of
principalization and of desingularization is stated from  a
different point of view, using a new presentation of the local
invariant.

\

We also should mention the papers \cite{AbrJong97},
\cite{AbrWang97}, and \cite{BogPan96}, where  a completely
different proof of resolution of singularities it is given,
leading to the so called {\em weak desingularization of
varieties}.  As opposed to the algorithms mentioned before, in the
weak desingularization approach regular points of the variety may
be modified in the desingularization process. More generally, see
\cite{Hauser} for an interesting introduction to the problem of
desingularization.

\

As indicated above, this manuscript is a self contained presentation which intends to provide an
introduction to the subject, particularly focused in the motivation of the main ideas, and hopefully a
reference for future applications. To this end  we have
  organized the paper so that the reader can
easily get into the statement of the main results and some of their applications without having to go
through  the technicalities of the proofs.

\

More precisely, the paper is organized as follows: In the next section, the last of Part 1,  we state
the theorems of resolution of singularities and principalization of ideals that we are going to prove.
\medskip

Part \ref{partbo} is devoted to presenting the theory of basic objects: This notion is introduced in
Section \ref{secbo}, and those of equivariance, and of algorithmic resolution of basic objects, in
Sections \ref{equiv} and \ref{algorithm}. After presenting the properties  of  this implemented algorithm
of resolution of basic objects, in  \ref{existence} we give a short proof of desingularization  for
equidimensional schemes (in \cite{BV1} we provide  a proof which works in the general case). In Section
\ref{examplesfamilies} we introduce the two main equivariant functions involved in the construction of the
algorithms of resolution.  \medskip

In Part \ref{applications} we present applications of our results: In Section \ref{weakvsstrict} the
notions of weak and strict transforms of ideals are discussed, and a stronger version of Hironaka's
Theorem is stated; in Section \ref{completions} we discuss the extension of our results to a wider class
of schemes (which are not necessarily of finite type over  a field); Sections \ref{nonembeddeddes} and
\ref{families} deal respectively with resolution of locally embedded varieties  and resolution in
families of varieties;  finally,  in Section \ref{Computer} there are some indications on how the
algorithm of desingularization can be implemented as a computer program which works for arbitrary
varieties (and not only in the hypersurface case).
\medskip

Part \ref{nature} is devoted to giving a flavor of how the
algorithm of resolution of basic objects works. The algorithm is
described in Part \ref{appendix}.

\

{\bf Acknowledgments.}  Part of the material contained in  this manuscript has been  presented in a
series of seminars  at  the University of  Purdue, during the Spring in 2001. We profited from
discussions with Gabrielov, Matsuki, Moh and Nobile.




\section{Formulation of the Theorems.}
\label{lang} The main point in the statement of the problems of resolution of singularities and
principalization of ideals is the construction of a birational  morphism of regular schemes
\(\pi:W^{\prime}  \to W\). This morphism will be defined as a composition of monoidal transformations at
closed and regular centers. If \(e_1: W_1\longrightarrow W\) is a monoidal transformation with center
\(Z\), then the exceptional locus \(H_1\subset W_1\) is a regular hypersurface, and if
\begin{equation*}
W_2\stackrel{e_2}{\longrightarrow}W_1 \stackrel{e_1}{\longrightarrow} W
\end{equation*}
is a composition of monoidal transformations with centers \(Z\subset W\) and \(Z_1\subset W_1\), then the
exceptional locus of \(e_1\circ e_2: W_2\longrightarrow W\) is the union of the strict transform of
\(H_1\) with the exceptional locus \(H_2\) of \(e_2: W_2\longrightarrow W_1.\) Both Theorems
\ref{classical} and \ref{principalization} require  the exceptional locus of \(\pi:
W^{\prime}\longrightarrow W\) to have normal crossings (see Definition \ref{normalc} below), and to
achieve this condition we will have to impose some constraints on the regular centers \(Z\subset W\) and
\(Z_1\subset W_1\). With this purpose we introduce the notions of {\em pairs} and {\em transformation of
pairs} which are suitable for the formulation of both {\em Embedded Desingularization} and {\em
Principalization of Ideals}.

\begin{definition}
\label{normalc} Let $W$ be a regular scheme and let  $Y_1,\ldots,Y_k\subset W$ be a set of closed
subschemes.  We say that $Y_1\cup \ldots\cup Y_k$ have {\em normal crossings at a point  $\xi\in W$}  if
there exists a regular system of parameters  $\{x_1,\ldots,x_d\}\subset {\mathcal O}_{W,\xi}$, such that
for each $i\in \{1,\ldots,k\}$, either $\cali(Y_i)_{\xi}={\mathcal O}_{W,\xi}$,  or
$$\cali(Y_i)_{\xi}=\langle x_{i_1},\ldots,x_{i_{s_i}}\rangle$$ for some  $ x_{i_1},\ldots,x_{i_{s_i}}\in
\{x_1,\ldots,x_d\}$.  We say that $Y_1\cup\ldots\cup Y_k$ have {\em normal crossings in $W$}  if they
have normal crossings at any point of $W$.
\end{definition}

\begin{definition}
\label{pairs} Let \( W \)  be  a pure dimensional scheme, smooth over a field \( \mathbf{k} \) of
characteristic zero, and let \( E=\{H_{1},\ldots,H_{r}\} \) be a set of smooth hypersurfaces in \( W \)
with  normal crossings. The couple  \( (W,E) \) is said to be a {\em pair}.
\end{definition}

\begin{parrafo}
\label{transformationop} {\bf Transformation of pairs.} (\cite[2.10]{EncVil97:Tirol},
\cite[1.1]{Villa92})  A regular closed subscheme \( Y\subset W \) is said to be {\em permissible} for the
pair \( (W,E) \) if \( Y \) has normal crossings with  \( E \), (i.e. if $Y$ has normal crossings with
$E=\{H_{1},\ldots,H_{r}\}$). If \( Y \subset W\) is permissible  for a pair \((W,E)\), we define a {\em
transformation  of  pairs} in the following way: Consider the blowing-up with center \(Y\),
     $$W\stackrel{\Pi}{\longleftarrow} W_1,$$ and define \(
E_{1}=\{H'_{1},\ldots,H'_{r},H'_{r+1}\} \), where \(H'_{i} \) denotes the strict transform of \( H_{i}
\), and \( H'_{r+1}\) denotes \(\Pi^{-1}(Y) \),  the exceptional hypersurface in \( W_{1} \). Note that
\(W_{1}\) is smooth and that \(E_{1}\) has normal crossings. We say that \((W,E) \longleftarrow
(W_{1},E_{1})\) is a {\em transformation of the pair} \((W,E)\).
\end{parrafo}

Now we state the two theorems that we prove:

\begin{theorem}[{\bf Embedded Resolution of Singularities}]
\label{classical} Let \((W_0, E_0=\{\emptyset\})\) be a pair and let \(X_0\subset W_0\) be  a closed
equidimensional subscheme defined by
        \(\cali(X_0) \subset \calo_{W_0}\). Assume that the open set \(
        \Reg(X) \) of regular points of $X$  is dense in \(X\).  Then
        there exists a finite sequence of transformations of pairs
\begin{equation*}
(W_0,E_0=\emptyset)\longleftarrow \cdots \longleftarrow (W_r,E_r),
\end{equation*}
        inducing a proper birational morphism
         \(\Pi_r:W_{r}\longrightarrow W_{0}\), such that if
         \(E_r=\{H_1,\ldots,H_r\}\) is the exceptional locus of $\Pi_r$
         then:
\begin{enumerate}
\item[(i)] The strict transform of $X$ in $W_r$, \(X_r\),  is regular
in \(W_r\), and \(W_r\setminus\cup_{i=1}^rH_i\simeq W_0-\Sing (X)\). In particular
\(\Reg(X)\cong\Pi^{-1}_r(\Reg(X))\subset X_{r}\) via \( \Pi_r \).
\item[(ii)] The scheme \( X_{r} \) has normal crossings with \(
E_{r}=\cup_{i=1}^rH_i \).
\item[(iii)] (Equivariance) If a group acting on \(W\) is also
acting on \(X_0\subset W_0\), then the action can be lifted to one on \(X_r\subset W_r\).
\end{enumerate}
\end{theorem}

\begin{theorem}[{\bf Embedded Principalization of ideals}]
\label{principalization} Let \( (W_{0},E_{0}) \) be a pair and let \( I\subset\O_{W_{0}} \) be a non zero
sheaf of ideals. Then there exists an  embedded principalization of \(I\), i.e. there is a finite
sequence of transformations of pairs
\begin{equation}
\label{esta} (W_0,E_0=\emptyset)=(W,E)\longleftarrow \ldots\longleftarrow (W_r,E_r),
\end{equation}
at smooth centers $Y_i\subset W_i$,  such that:
\begin{enumerate}
\item[(i)] The morphism  \(W_r\to W\) defines an isomorphism over the
open subset \(W\setminus V(I)\).
\item[(ii)] The ideal $I\calo_{W_r}$ is invertible and
supported on a divisor with normal crossings, i.e.,
\begin{equation}
{\mathcal {L}}=I\calo_{W_r}=\cali(H_1)^{c_1} \cdot\ldots\cdot\cali(H_s)^{c_s},
\end{equation}
where \(E^{\prime}=\{ H_1,H_2,\dots, H_s\}\) are regular hypersurfaces with normal crossings,  \(c_i\geq
1\) for \(i=1,\ldots,s\), and $E^{\prime}=E_r$ if $V(I)$ has no components of codimension 1.
\item[(iii)] (Equivariance) If a group $G$ acts on
$W_0$ and on the ideal $I$, then it acts on each center $Y_i$ for $i=0,\ldots, r-1$, and the action can
be lifted at each step of the resolution.
\end{enumerate}
\end{theorem}

\begin{remark}
In particular, if $W={\mathbb A}_k^n$, or any toric variety, $J$ is a monomial ideal and
\begin{equation}
(W_0,E_0=\emptyset)=(W,E)\longleftarrow \ldots\longleftarrow (W_r,E_r),
\end{equation}
is an embedded principalization of $J$, obtained by  blowing up at smooth centers $Y_i\subset W_i$ as
indicated in Theorem \ref{principalization},  then all the ideals  ${\mathcal I}(Y_i)$  are  monomial
ideals, because the sequence of transformations is equivariant,    and the  torus acts on $W$ and $J$,
and hence on the centers  $Y_i$, for $i=0,1,\ldots,r-1$.
\end{remark}

The proofs of Theorems \ref{classical} and \ref{principalization}
will be given in \ref{existence} and \ref{proofprincipal},  where
it will be shown that both theorems are a direct consequence  of
the existence of an algorithm of resolution of basic objects.
Basic objects and algorithms of resolution of basic objects are
discussed in the upcoming sections.




\part{Basic objects.}
\label{partbo}
\section{Basic objects.}
\label{secbo} In this section we recall the definitions of {\em basic objects} and {\em resolution of
basic objects} (Definitions \ref{DefBasic} and \ref{defResol}, see also \cite{EncVil97:Tirol}).  Both
Embedded Desingularization and Strong Principalization of Ideals can be obtained from a resolution of
suitably defined basic objects.

\begin{definition}
\label{DefBasic} A \emph{basic object} is a triple that consists of a pair \( (W,E) \), an ideal \(
J\subset\calo_{W} \) such that \((J)_{\xi}\neq 0\) for any \( \xi\in W \), and a positive integer \( b
\).  It is denoted by \( (W,(J,b),E) \).  If the dimension of \(W\) is \(d\), then \( (W,(J,b),E) \) is
said to be a {\em \(d\)-dimensional basic object}.
\end{definition}

\begin{definition}\label{singular}
The {\em singular locus} of a  basic object is the closed subset of \(W\),
\begin{equation*}\Sing(J,b)=\{\xi\in W\mid \nu_{J}(\xi)\geq b\}
        \subset W, \end{equation*} where \(\nu_{J}(\xi)\) denotes the
order of the ideal \(J\) at the local regular ring \(\calo_{W,\xi}\).
  Sometimes we will use the notation
\begin{equation*}
\begin{array}{ccc}  & (W,(J,b),E) & \\
     & \cup & \\ & \Sing(J,b)&
\end{array}
\end{equation*}
to refer to the basic object $(W,(J,b),E)$ together with the closed subset $\Sing(J,b)$, meaning that
\(\Sing(J,b) \subset W\).
\end{definition}

\begin{example}\label{example1o}
If $X\subset W$ is a hypersurface and $J$ is its defining ideal, then $\Sing(J,b)$ is the set of points
of $X$ where the multiplicity is greater than or equal to $b$.
\end{example}

\begin{example}\label{example2o}
If $J\subset {\mathcal O}_{W}$ is an arbitrary non-zero sheaf of ideals, then $\Sing(J,b)$ is the set of
points of $W$ where the order of $J$ is greater than or equal to $b$.
\end{example}

\begin{definition}\label{permissible}
A regular closed subscheme \( Y\subset W \) is  \emph{permissible} for \( (W,(J,b),E) \) if \( Y \) is
permissible for the pair \( (W,E) \) and \( Y\subset \Sing(J,b)\).
\end{definition}

\begin{parrafo}
\label{transformation} {\bf Permissible transformations of basic objects.} Let
\((W,(J,b),E=\{H_1,\ldots,H_r\})\) be a basic object and let \(Y\subset \Sing(J,b)\) be a permissible
center. Consider \(W \longleftarrow W_1\) the monoidal transformation with center \(Y\). This induces a
transformation of  pairs,
\begin{equation*}
(W,E=\{H_1,\ldots,H_r\})\longleftarrow (W_1,E_1=\{H_1,\ldots,H_r,H_{r+1}\}),
\end{equation*}
where \(H_{r+1}\subset W_{1} \) is the exceptional divisor (see Definition \ref{pairs}).  Now, since
$Y\subset \Sing(J,b)$,
\begin{equation}
\label{formula111} J{\mathcal O}_{W_1}=\cali(H_{r+1})^bJ_1,
\end{equation}
for some $J_{1}\subset\O_{W_{1}}$. Then we define
\begin{equation*}
(W,(J,b),E)\longleftarrow (W_{1},(J_{1},b),E_{1})
\end{equation*}
as the \emph{permissible transformation of the basic object} \((W,(J,b),E)\).
\end{parrafo}

\begin{remark} \label{remarkper}
In general, given a sequence of transformations of basic objects
\begin{equation}
\label{resol}
\begin{array}{ccccccc}
(W_{0},(J_{0},b),E_{0}) & \longleftarrow & (W_{1},(J_{1},b),E_{1})&
\longleftarrow & \cdots & \longleftarrow & (W_{k},(J_{k},b),E_{k})\\
\cup & & \cup & & & & \cup \\ \Sing(J_0,b) & & \Sing(J_1,b) &   & \ldots & &  \Sing(J_k,b)
\end{array}
\end{equation}
at centers \(Y_i\subset \Sing(J_i,b)\), \(i=0,1,\ldots k-1\), we obtain expressions
\begin{equation}
\label{tresestrellas} J_{0}{\mathcal O}_{W_i}=\cali(H_{r+1})^{c_{r+1}}\cdot\ldots\cdot
\cali(H_{r+i})^{c_{r+i}} J_i.
\end{equation}
Note here that \(c_{r+1}=\ldots=c_{r+i}=b\) if none of the centers \(Y_i\) are included in any of the
exceptional divisors \(H_{j}\).
\end{remark}

\begin{definition}
\label{defResol} A finite sequence transformation of basic objects as (\ref{resol}) is a
\emph{resolution  of} \( (W_{0},(J_{0},b),E_{0}) \) if \( \Sing(J_{k},b)=\emptyset \).
\end{definition}

\begin{remark}
\label{so} Note that:
\begin{enumerate}
\item If sequence (\ref{resol}) is a resolution of the basic object
\((W_0,(J_0,b),E_0)\),  then $W_k\longrightarrow W_0$ defines an isomorphism over \(W_0\setminus
V(J_0)\), and \(J_0{\mathcal O}_{W_k}={\mathcal M}_kJ_k\), where \({\mathcal M}_k\) is an invertible
sheaf of ideals, and \(J_k\) has no points of order \( \geq b\) in \(W_k\).
\item The ideal  \(J_k\) is not the {\em strict transform} of \(J_0\),
an ideal which is far more complicated to define (see Section \ref{weakvsstrict} for  a discussion on
this matter).  However it is so in some particular cases. In fact, if \(X_0 \subset W_0\) is a closed
smooth subscheme, and \(J_0=I(X)\), then \(\Sing(J_0,1)=X\), and   given any sequence of transformations
of basic objects $$(W_0,(J_0,1),E_0)\longleftarrow (W_{1},(J_{1},1),E_{1}) \longleftarrow
\cdots\longleftarrow (W_{k},(J_{k},1),E_{k})$$ the ideal \(J_k\) is the ideal of a smooth subscheme \(X_k
\subset W_k\), which is the strict transform of \(X \subset W_0\).
\end{enumerate}
\end{remark}

\begin{theorem}\cite[Theorem 7.3.]{Villa92}
\cite[Theorem 7.13.]{EncVil97:Tirol} \label{resolucionob} Given a  basic object \((W,(J,b),E)\)
    where \(W\) is smooth over a field of characteristic zero, there is
a  resolution, i.e. there is  a finite sequence of monoidal transformations at permissible centers
\(Y_i\subset W_i\),
\begin{equation} \label{ressol}
(W,(J,b),E)=(W_0,(J_0,b),E_0)\longleftarrow (W_1,(J_1,b),E_1) \longleftarrow \ldots \longleftarrow
(W_k,(J_k,b),E_k),
\end{equation}
     such  that \(\Sing(J_k,b)=\emptyset.\)
\end{theorem}

\begin{remark}\label{remarkresol}
Note that:
\begin{enumerate}
\item[(i)] Theorem \ref{resolucionob} is existential. It claims that
given \((W,(J,b),E)\) there exists a resolution. However we shall give a {\em constructive} proof (the
same as in \cite[Theorem 7.3]{Villa92},  \cite[Theorem 7.13]{EncVil97:Tirol}), so that given
\((W,(J,b),E)\) we will define one particular resolution.

\item[(ii)] We will show in \ref{existence} and \ref{proofprincipal}
that a constructive proof of  Theorem \ref{resolucionob} will lead us to simple constructive proofs of
Theorems \ref{classical} and \ref{principalization}.

\end{enumerate}
\end{remark}

A central point in our constructive proof of Theorem \ref{resolucionob}, and hence of Theorems
\ref{classical} and \ref{principalization}, relies on the fact that all the invariants involved in the
algorithmic resolution of a basic object \((W,(J,b),E)\), will be defined in terms of the  closed sets of
the form $F=\Sing(J,b)$. The definition of these invariants will be addressed in
Sections~\ref{examplesfamilies} and \ref{SeccConsFunc}. As it turns out, the invariants that we are
aiming to use do  not behave well if only transformations as the ones of Definition \ref{transformation}
were allowed (see Section \ref{equiv} for equivariance, and Remark \ref{allow} where this fact is
illustrated). This forces us to enlarge the class of transformations of basic objects by including the
{\em projections}.

\begin{parrafo} \label{projections} {\bf Projections.} Let
\((W,(J,b),E)\) be a basic object. We will define a  new notion of {\em transformation of pairs and of
basic objects} as follows: Set
\begin{equation*}
     W_1=W\times {\mathbb A}_k^1
\end{equation*}
and consider the natural projection \(p: W_1\longrightarrow W.\) Then define
\begin{equation*}
J_1=J{\mathcal O}_{W_{1}} \mbox{ } \mbox{ } \mbox{ and } \mbox{ } \mbox{ } E_1=p^{-1}(E).
\end{equation*}
This gives us a transformation of pairs
\begin{equation*}
(W,E) \longleftarrow (W_1, E_1),
\end{equation*}
and a transformation of basic objects,
\begin{equation*}
(W,(J,b),E) \longleftarrow (W_1, (J_1,b), E_1).
\end{equation*}
In this case the closed subset in \(W_1\) is \(F_1=\Sing(J_1,b)=p^{-1}(\Sing(J,b)).\)
\end{parrafo}

In what follows, unless otherwise specified, whenever we mention a {\em transformation} or a {\em
sequence of transformations} of  basic objects, we will refer to both the permissible transformations
introduced in Definition \ref{transformation} and the projections that we have just defined.  In both
cases we will use the same notation,
\begin{equation}
\label{referida}
\begin{array}{ccccccc}
(W_{0},(J_{0},b),E_{0}) & \longleftarrow & (W_{1},(J_{1},b),E_{1})&
\longleftarrow & \cdots & \longleftarrow & (W_{k},(J_{k},b),E_{k})\\
\cup & & \cup & & & & \cup \\ \Sing(J_0,b) & & \Sing(J_1,b) &   & \ldots & &  \Sing(J_k,b),
\end{array}
\end{equation}
and we will assume that both kinds of transformations  can be used in the same sequence.




\section{Equivariance.}
\label{equiv}

An important outcome of constructive desingularization is the
lifting of any  group action on a subscheme \(X\subset W\), all
the way up to the desingularization. Since our desingularization
theorem will follow from a suitably defined resolution of basic
objects, in this section we present and discuss the notion of
equivariance in the  context of basic objects. We should point out
that we will work with isomorphisms \(\Theta: W\to W^{\prime}\)
which may be not necessarily defined over the base field
\(\mathbf{k}\) over which the smooth schemes \(W\) and
\(W^{\prime}\) are defined.

\begin{definition}
     \label{DefActPar}
An {\em isomorphism of pairs}
\begin{equation*}
\Theta:(W,E=\{H_1,\ldots,H_k\})\longrightarrow (W^{\prime},E^{\prime}=
\{H_1^{\prime},\ldots,H_k^{\prime}\})
\end{equation*}
is an isomorphism \(\Theta: W \longrightarrow W^{\prime}\) such that \(\Theta(H_i)=H_i^{\prime}\) for all
\(i=1,\ldots,k\).
\end{definition}

\begin{lemma}
\label{lifting} Let \(\Theta: (W,E)\longrightarrow (W^{\prime},E^{\prime}) \) be  an isomorphism of
pairs, let \( Y \) be a  permissible center for \((W,E)\) and set \(Y^{\prime}=\Theta(Y).\) Consider the
transformations of pairs with centers \(Y\),
\begin{equation*}
(W,E)\longleftarrow (W_1,E_1),
\end{equation*}
and \(Y^{\prime}\),
\begin{equation*}
(W^{\prime},E^{\prime})\longleftarrow (W^{\prime}_1,E^{\prime}_1).
\end{equation*}
Then there is a natural lifting of the isomorphism \(\Theta: (W,E)\longrightarrow
(W^{\prime},E^{\prime})\) to an isomorphism
\begin{equation} \label{levantiso}
\Theta_1: (W_1,E_1)\longrightarrow (W^{\prime}_1,E^{\prime}_1)
\end{equation} such that the diagram
\begin{equation*}
\begin{array}{ccc} (W,E) &
\longleftarrow & (W_1,E_1)\\ \Theta \downarrow & & \downarrow\Theta_1\\ (W^{\prime},E^{\prime}) &
\longleftarrow & (W_1^{\prime},E_1^{\prime})
\end{array}
\end{equation*}
commutes.
\end{lemma}

\noindent{\em Proof:} The map \(\Theta: (W,E)\longrightarrow (W^{\prime},E^{\prime}) \) induces an
isomorphism of the corresponding structural sheaves \({\mathcal O}_{W}\)  and \({\mathcal
O}_{W^{\prime}}\), and extends naturally to an  isomorphism
\begin{equation*}
     \hat{\Theta}: W\times\mathbb{A}_{\mathbf{k}^{1}}\longrightarrow
W^{\prime}\times\mathbb{A}_{\mathbf{k}^{1}},
\end{equation*}
i.e. an isomorphism between   \( \Oh_{W}[Z] \) and \(\Oh_{W^{\prime}}[Z^{\prime}]\), where \(Z\) and
\(Z^{\prime}\) are indeterminates and $\Theta(Z)=Z^{\prime}$.

Since \(\Theta(Y)=Y^{\prime}\), the sheaf of ideals
\(\cali(Y)\subset {\mathcal O}_{W}\) is mapped to
\(\cali(Y^{\prime})\subset{\mathcal O}_{W^{\prime}}\). As a
consequence, \(\Theta\) induces an isomorphism
\(\Oh_{W}[I(Y)Z]\simeq \Oh_{W^{\prime}}[I(Y^{\prime})
Z^{\prime}]\) which  preserves the grading. Therefore we get an
isomorphism
\begin{equation*}
W_{1}=\Proj\left(\Oh_{W}[I(Y)Z]\right)\simeq W_{1}^{\prime}= \Proj\left(\Oh_{W^{\prime}}[I(Y
^{\prime})Z^{\prime}]\right).
\end{equation*}
We finally check that this isomorphism maps the exceptional hypersurface of one to that of the other,
defining \(\Theta_1: (W_1,E_1)\longrightarrow (W_1^{\prime},E_1^{\prime}) \).  \qed

\begin{remark}\label{isotipo2}
Let \(\Theta:(W,E) \to (W^{\prime},E^{\prime})\) be an isomorphism of pairs. Set \( W_1=W\times {\mathbb
A}_k^1,\) \( W^{\prime}_1=W^{\prime}\times {\mathbb A}_k^1\), both as in \ref{projections}; and set
\(\Theta_1= \Theta \times \mbox{Id}_{{\mathbb A}_k^1}\). Note that \(\Theta_1:(W_1,E_1) \to
(W_1^{\prime},E_1^{\prime})\) is an isomorphism of pairs.
\end{remark}

\begin{remark} \label{lifting1}
Fix an isomorphism \(\Theta: (W,E)\longrightarrow
     (W^{\prime},E^{\prime})\) and a sequence
\begin{equation}
\label{iso}
\begin{array}{ccccccc}
(W,E) & \longleftarrow & (W_1, E_1)& \longleftarrow & \ldots & \longleftarrow & (W_k, E_k),
\end{array}
\end{equation}
where each transformations of pairs is either as in \ref{transformationop} or as in \ref{projections}.
Then Lemma \ref{lifting} and Remark \ref{isotipo2} assert that \(\Theta: W \longrightarrow W^{\prime}\)
defines a unique sequence of transformation of pairs,
\begin{equation}
\label{iso2}
\begin{array}{ccccccc}
(W^{\prime},E^{\prime}) & \longleftarrow & (W_1^{\prime}, E_1^{\prime})& \longleftarrow & \ldots &
\longleftarrow & (W_k^{\prime}, E_k^{\prime}).
\end{array}
\end{equation}
together with isomorphisms \(\Theta_i: (W_i,E_i)\longrightarrow
     (W_i^{\prime},E_i^{\prime})\), for $i=1,\ldots,k$.
\end{remark}

\begin{remark} \label{DefRestricBO}
Let \( \mathcal{B}=(W,(J,b),E) \) be a basic object.  If \( U\subset W \) is an open set, then we set the
\emph{restriction} of the basic object to be \( \mathcal{B}|_{U}=(U,(J|_{U},b),E_{U}) \), where \( J|_{U}
\) is the restriction of the sheaf of ideal to \( U \) and \( E_{U}=\{H\cap U\mid H\in E\} \).
\end{remark}

\begin{definition}\label{isoBO}
Let \((W,(J,b),E)\) and \((W^{\prime},(J^{\prime},b^{\prime}),E^{\prime})\) be two basic objects and  let
\begin{equation*}
\theta: (W, E)\to (W^{\prime},E^{\prime})
\end{equation*}
be an isomorphism of pairs. We will say that \(\theta\) induces an {\em isomorphism of basic objects},
\begin{equation*}
\Theta :(W,(J,b),E) \to (W^{\prime},(J^{\prime},b^{\prime}),E^{\prime}),
\end{equation*}
if the following conditions hold:
\begin{enumerate}
\item[(i)] The isomorphism \(\theta: W\to W^{\prime}\) induces an
isomorphism of the closed subsets defined by the basic objects
\begin{equation*}
\theta:F= \Sing(J,b)\cong F^{\prime}=\Sing(J^{\prime},b^{\prime}).
\end{equation*}
\item[(ii)] Given a sequence of transformations of basic objects as in
\ref{transformation} or \ref{projections},
\begin{equation}\label{equivar1}
\begin{array}{ccccccc} (W,(J,b),E) & \longleftarrow &
(W_1,(J_1,b), E_1)& \longleftarrow & \ldots & \longleftarrow & (W_k, (J_k,b), E_k),
\end{array}
\end{equation}
  together with the
corresponding sequence of transformation of pairs
\begin{equation*}
\begin{array}{ccccccc}
(W,E) & \longleftarrow & (W_1, E_1)& \longleftarrow & \ldots & \longleftarrow & (W_k, E_k)
\end{array}
\end{equation*}
then:
\begin{enumerate}
\item[(a)] The corresponding induced sequence, as in Remark
\ref{lifting1},
\begin{equation*}
\begin{array}{ccccccc}
(W^{\prime},E^{\prime}) & \longleftarrow & (W_1^{\prime}, E_1^{\prime})& \longleftarrow & \ldots &
\longleftarrow & (W_k^{\prime}, E_k^{\prime}),
\end{array}
\end{equation*}
defines a sequence of transformation of basic objects,
\begin{equation}\label{equivar2}
\begin{array}{ccccccc}
(W^{\prime},(J^{\prime},b^{\prime}),E^{\prime}) & \longleftarrow &
(W_1^{\prime},(J_1^{\prime},b^{\prime}), E_1^{\prime})& \longleftarrow & \ldots & \longleftarrow &
(W_k^{\prime}, (J_k^{\prime},b^{\prime}), E_k^{\prime})
\end{array}
\end{equation}

\item[(b)] The isomorphisms \(\theta_i: (W_i,E_i)\longrightarrow
(W_i^{\prime},E_i^{\prime})\) defined in  remark \ref{lifting1}, induce an isomorphism of the closed
subsets defined by the basic objects,
  $$\theta_i:\Sing(J_i,b)\cong
\Sing(J_i^{\prime},b^{\prime})$$ for all \(i=1,\ldots,k\).
\end{enumerate}
\item[(iii)] For any open set \( U \), set \( U'=\theta(U) \) and
consider the restrictions \( (U,(J|_{U},b),E_{U}) \) and \( (U',(J'|_{U'},b'),E'_{U}) \). We require that
properties (i) and (ii) hold for these restrictions.
\end{enumerate}
\end{definition}

\begin{example} \label{exiso}
Let \( \Theta: (W,E) \to (W^{\prime},E^{\prime})\) be an isomorphism of pairs, fix a basic object
\((W,(J,b),E)\) and let   \(J^{\prime}\) be the image of \(J\) in \({\mathcal O}_W^{\prime}\) under
\(\Theta\), i.e. \(J^{\prime}=\Theta(J) \subset \calo_{W^{\prime}}.\) Then we claim that
\begin{equation*}
     \Theta: (W,(J,b),E) \to (W^{\prime},(J^{\prime},b),E^{\prime})
\end{equation*}
     is an isomorphism of basic objects: Clearly   \(\Theta: W \to
W^{\prime}\) maps \(\Sing(J,b)\) isomorphically into \(\Sing(J^{\prime},b)\). To check that condition
(ii) of Definition \ref{isoBO} holds, note that if
\begin{equation}\label{cinsiet}
(W,(J,b),E) \longleftarrow (W_1,(J_1,b),E_1)
\end{equation}
is a transformation of basic objects as the ones defined in \ref{projections}, then the corresponding
transformation of basic objects  over \((W^{\prime},(J^{\prime},b),E^{\prime})\),
\begin{equation*}
(W^{\prime},(J^{\prime},b),E^{\prime}) \longleftarrow
(W^{\prime}_1,(J^{\prime}_1,b),E^{\prime}_1)\end{equation*} is such that the there is a natural way to
define an  isomorphism
\begin{equation*}
\Theta_1: W_1 \to W_1^{\prime},
\end{equation*}
which maps \(\Sing(J_1,b)\) isomorphically into \(\Sing(J_1^{\prime},b)\). In fact, \( \Theta_1\) maps
\(J_1= J \calo_{W_{1}}\) into \(J^{\prime}_1= J^{\prime} \calo_{W^{\prime}_{1}}\).  \medskip

Now if
\begin{equation}
\label{cinsiet1} (W,(J,b),E) \longleftarrow (W_1,(J_1,b),E_1)
\end{equation}
is a transformation of basic objects with center \(Y\) as in Definition  \ref{transformation}, and if
\begin{equation*}
(W^{\prime},(J^{\prime},b),E^{\prime}) \longleftarrow (W^{\prime}_1,(J^{\prime}_1,b),E^{\prime}_1)
\end{equation*}
     is a transformation of basic objects with center
\(Y^{\prime}=\Theta(Y)\), then by Lemma \ref{lifting} we obtain an isomorphism of pairs
\begin{equation*}
     \Theta_1: (W_1,E_1) \to (W^{\prime}_1,E^{\prime}_1).
\end{equation*}
Since \(\Theta_1: W_1 \cong W^{\prime}_1\) is compatible with \(\Theta: W \cong W^{\prime}\) it follows
that
\begin{equation*}
     \Theta_1(J \calo_{W_{1}})=J^{\prime} \calo_{W^{\prime}_{1}}.
\end{equation*}
In particular, by (\ref{formula111}), \(\Theta_1(J_1)=J^{\prime}_1\), and hence
\(\Theta_1(\Sing(J_1,b))=\Sing(J^{\prime}_1,b)\).
\end{example}

\begin{example}\label{exampleequiv1}
The identity map induces an isomorphism of the basic objects,
\begin{equation*}
       \text{Id}:(W,(J,b),E) \to (W,(J^2,2.b),E),
\end{equation*}
(see also Example \ref{exampleequiv}).
\end{example}

\begin{remark}
Note that if two basic objects are isomorphic, then there might be may different isomorphisms between
them. We are interested in resolutions which are compatible with {\bf any} of these  isomorphisms. This
is the philosophy behind the next definition.
\end{remark}

\begin{definition}\label{defiequivar}
Let $(W,(J,b),E)$ be a basic object. We will say that a sequence of transformations of basic objects,
\begin{equation*}
       (W,(J,b),E) \longleftarrow (W_1,(J_1,b),E_1) \longleftarrow
       \ldots \longleftarrow (W_k,(J_k,b),E_1)
\end{equation*}
    is {\em equivariant}, if for any isomorphic basic object
\((W^{\prime},(J^{\prime},b),E^{\prime})\), the induced sequence of transformations defined as in
(\ref{equivar2}),
\begin{equation} \label{EqDefEquivar}
       (W^{\prime},(J^{\prime},b^{\prime}),E^{\prime}) \longleftarrow
       (W_1^{\prime},(J_1^{\prime},b^{\prime}),E_1^{\prime})
       \longleftarrow \ldots \longleftarrow
       (W_k^{\prime},(J_k^{\prime},b^{\prime}),E_1^{\prime}),
\end{equation}
is independent of the isomorphism \(\Theta:(W,(J,b),E) \to
(W^{\prime},(J^{\prime},b^{\prime}),E^{\prime})\) that we take.
\end{definition}

\begin{parrafo} \label{eqbaso}
{\bf What does equivariance mean?} In Section \ref{algorithm} we will introduce the notion of {\em
algorithm of resolution of basic objects}. The algorithm will provide, for each basic object, a
resolution which is equivariant, i.e., let
\begin{equation}
\label{resol3}
\begin{array}{ccccccc}
(W_{0},(J_{0},b),E_{0}) & \longleftarrow & (W_{1},(J_{1},b),E_{1})&
\longleftarrow & \cdots & \longleftarrow & (W_{k},(J_{k},b),E_{k})\\
\cup & & \cup & & & & \cup \\ \Sing(J_0,b) & & \Sing(J_1,b) &   & \ldots & & \Sing(J_k,b)=\emptyset
\end{array}
\end{equation}
be the resolution of \({\mathcal B}=(W_{0},(J_{0},b),E_{0})\) provided by the algorithm, and let
\begin{equation}
\label{resol4}
\begin{array}{ccccccc}
(W^{\prime}_{0},(J^{\prime}_{0},b^{\prime}),E^{\prime}_{0}) & \longleftarrow &
(W^{\prime}_{1},(J^{\prime}_{1},b),E^{\prime}_{1})& \longleftarrow & \cdots & \longleftarrow &
(W^{\prime}_{k^{\prime}},(J^{\prime}_{k^{\prime}},b^{\prime}),E_{k^{\prime}})\\
\cup & & \cup & & & & \cup \\ \Sing(J^{\prime}_0,b^{\prime}) & & \Sing(J^{\prime}_1,b^{\prime}) &   &
\ldots & & \Sing(J^{\prime}_{k^{\prime}},b^{\prime})=\emptyset
\end{array}
\end{equation}
be that of \({\mathcal B}^{\prime}=(W^{\prime}_{0},(J^{\prime}_{0},b^{\prime}),E^{\prime})\). We require
that if \( \mathcal{B} \) and \( \mathcal{B}' \) are isomorphic, i.e. if $$\Theta: (W_0,(J_0,b),E_0) \to
(W_0^{\prime},(J_0^{\prime},b),E_0^{\prime})$$ is an isomorphism, then sequence (\ref{EqDefEquivar})
induced by (\ref{resol3}) is precisely the sequence (\ref{resol4}),  and hence the isomorphism can be
lifted along the resolutions  (\ref{resol3}) and (\ref{resol4}).

\

Since the resolutions are equivariant, $\Theta$ can be lifted,  and furthermore, any other isomorphism
$\Gamma: {\mathcal B} \to {\mathcal B}^{\prime}$ can also be lifted.  An algorithm with this property is
said to be an {\em equivariant algorithm}. This property, and others, will arise, in a very simple
manner, by defining the resolution in terms of upper semi-continuous functions with suitable properties.
\medskip
\end{parrafo}




\section{The algorithmic proof of
Theorems \ref{principalization} and  \ref{classical}.} \label{algorithm}

\label{last1} In this section we explain  what we mean by a {\em constructive or algorithmic proof} of
Theorem~\ref{resolucionob}. This constructive proof will be done in terms of an {\em algorithm of
resolution of basic objects} (see Definition \ref{AlgResol}).  In order to define such an algorithm we
first  introduce a totally ordered set, and then we  attach to each basic object an {\em upper
semi-continuous function} with values on this   ordered set. These upper semi-continuous functions will
define a resolution for  each basic object \((W,(J,b),E)\), whenever \(W\) is defined over a field of
characteristic zero. We will refer to this resolution as {\em the resolution defined by the algorithm.}
At the end of the section we will indicate why an algorithm of resolution of basic objects
  as the one described in Definition  \ref{AlgResol} already provides
a proof of Theorems~\ref{classical} and \ref{principalization} (see \ref{existence} and
\ref{proofprincipal}).

\begin{definition}
\label{semicon} Let \(X\) be a topological space, let \((T, \geq )\) be a totally ordered set, and let
\(g: X \longrightarrow T\) be an
     upper semi-continuous function. Assume that \(g\) takes only
     finitely many values. Then the  largest value achieved by \(g\) will be
     denoted by
\begin{equation*}
\max g.
\end{equation*}
Clearly the set
\begin{equation*}
     \Max g=\{ x\in X : g(x)= \max g \}
\end{equation*}
     is a closed subset of \(X\).
\end{definition}

\begin{definition}
\label{familyequiv} Let  \(T\) be a totally ordered set. Assume that for any basic   object
\({\calB}=(W,(J,b),E)\) there is an upper semi-continuous function
\begin{equation*}
f_{\calB}: \Sing(J,b)\to T
\end{equation*}
     associated to it. This is what we will call a {\em family of
functions with values on \(T\)}. Given a sequence of transformations of basic objects,
$$(W_k,(J_k,b),E_k)\ldots\longleftarrow (W_1,(J_1,b),E_1)\longleftarrow (W_0,(J_0,b),E_0),$$ we will
denote by $f_i$ the corresponding function $$f_{{\calB}_i}: \Sing(J_i,b)\to T$$ associated to
${\calB}_i=(W_i,(J_i,b),E_i)$.
\end{definition}

\begin{definition} \label{familyequiv2}
A family of functions is said to be {\em equivariant} if for any isomorphism of basic objects
\begin{equation*}
\Theta: {\calB}=(W,(J,b),E)\to {\calB}^{\prime}=(W^{\prime},(J^{\prime},b),E^{\prime}),
\end{equation*}
we have that
\begin{equation} \label{EqFamilyEquiv2}
f_{{\calB}^{\prime}}\circ\Theta=f_{\calB}.
\end{equation}
\end{definition}

\begin{remark}\label{remarkfamily}
A family of functions with values on \(T\) attaches to each basic object a unique function. An
isomorphism \( \Theta \) as in Definition \ref{familyequiv2} defines an isomorphism \(W \cong W^{\prime}\)
mapping the domain of \(f_{{\calB}}\) into the domain of \(f_{{\calB}^{\prime}}\), and in particular,
mapping $\Sing(J,b) \subset W$ isomorphically into \(\Sing(J^{\prime},b) \subset W^{\prime}\) so  formula
(\ref{EqFamilyEquiv2}) makes sense.
\end{remark}

\begin{definition}
\label{AlgResol} An {\em algorithm of resolution of
     $d$-dimensional  basic objects}  consists of:
\begin{enumerate}
\item[(A)] A totally ordered set \( (I_{d},\leq) \).
\item[(B)] A family of  functions with values on $(I_{d},\leq)$,
\begin{equation*}
f_{\calB}: \Sing(J,b)\to I_{d},
\end{equation*}
which satisfies the following properties:
\begin{enumerate}
\item The closed subset
\begin{equation*}
\Max f_{\calB}
\end{equation*}
     is a smooth permissible center for \((W,(J,b),E)\), and it
therefore defines a transformation of basic objects
\begin{equation*}
(W,(J,b),E)\longleftarrow (W_1,(J_1,b),E_1).
\end{equation*}
     Hence, given a basic object \((W_0,(J_0,b),E_0)\), by this
property we always get a finite sequence of transformations of basic objects
\begin{equation}
\label{resolu}
\begin{array}{ccccccc}
(W_{0},(J_{0},b),E_{0}) & \longleftarrow & (W_{1},(J_{1},b),E_{1})&
\longleftarrow & \cdots & \longleftarrow & (W_{k},(J_{k},b),E_{k})\\
\cup & & \cup & & & & \cup \\ \Sing(J_0,b) & & \Sing(J_1,b) &   & \ldots & &  \Sing(J_k,b),
\end{array}
\end{equation}
with centers \(\Max  f_i\subset \Sing(J_i,b)\).
\item For any sequence   as (\ref{resolu}) we have that
\begin{equation*}
\max f_0 > \max f_1 >\ldots > \max f_k.
\end{equation*}

\item For any sequence as    (\ref{resolu})  there is an index \(N\),
depending on \((W_0,(J_0,b),E_0)\), such  that
\begin{equation*}
\Sing(J_N,b) =\emptyset.
\end{equation*}

\item  If  \( X_{0} \) is a regular pure dimensional subscheme  of
        dimension \(r\) then there is a value \(a(r)\in I_{d}\) so that
        if \( J_{0}=\cali(X_0)\), \( b=1 \) and \( E_{0}=\emptyset \),
        then the function \( f_{0} \) is constant and equal to $a(r)$.
\item  If \( \xi\in \Sing(J_{i},b) \) and  \( \xi\not\in \Max f_{i} \)
for \( i=0,\ldots,k-1 \), then \( f_{i}(\xi)=f_{i+1}(\xi') \) via the natural identification of the point
\( \xi \) with a point \( \xi' \) of \( \Sing(J_{i+1},b) \).
\item The family of functions $f_{\calB}$ is {\em equivariant}:  If
\begin{equation*}
\Theta: (W,(J,b),E)\longrightarrow (W^{\prime},(J^{\prime},b^{\prime}) E^{\prime})
\end{equation*}
     is an isomorphism of basic objects,  then the functions given by
the algorithm
\begin{equation*}
f: \Sing(J,b)\to T \mbox{ } \mbox{ } \mbox{ and }\mbox{ } \mbox{ } f^{\prime}: \Sing(J^{\prime},b)\to T
\end{equation*}
are such that \(f'\circ\Theta=f\).
    \end{enumerate}
\end{enumerate}
\end{definition}

\begin{parrafo} \label{propalg}
{\bf Some comments on the properties of the algorithm.}
\begin{enumerate}
\item[(1)] Property  (c) says that for $i=0,1,\ldots,k$, the functions
\begin{equation*}
f_i:\Sing(J_i,b) \to I_{d}
\end{equation*}
define a resolution  of the basic object \((W_0,(J_0,b),E_0)\). We will refer to it as  {\em the
resolution defined by the algorithm}.
\item[(2)] Compatibility with open restrictions: Let \(U\subset W\) be
an open subset and let
\begin{equation*}
(W,(J,b),E)_U=(U,(J|_{U},b),E_{U})
\end{equation*}
be the restriction to \(U\) of \((W,(J,b),E)\). Property (e)  asserts that the restriction to \(U\) of
the resolution of \((W,(J,b),E)\) defined by the algorithm, coincides with the resolution of the basic
object \((W,(J,b),E)_U\) defined by the algorithm.

\item[(3)] Let
\begin{equation}
\label{formula1} \Theta_0: (W_{0},(J_{0},b),E_{0})\to
(W_{0}^{\prime},(J_{0}^{\prime},b^{\prime}),E_{0}^{\prime})
\end{equation}
     be an isomorphism of basic objects. Now consider on the one hand
the resolution of \((W_{0},(J_{0},b),E_{0})\),
\begin{equation}
\label{resol1}
\begin{array}{ccccccc}
(W_{0},(J_{0},b),E_{0}) & \longleftarrow & (W_{1},(J_{1},b),E_{1})&
\longleftarrow & \cdots & \longleftarrow & (W_{k},(J_{k},b),E_{k})\\
\cup & & \cup & & & & \cup \\ \Sing(J_0,b) & & \Sing(J_1,b) &   & \ldots & &  \Sing(J_k,b)=\emptyset,
\end{array}
\end{equation}
defined by the  upper semi-continuous functions
\begin{equation*}
     f_i: \Sing(J_i,b)\to I_{d},
\end{equation*}
and on the other hand the resolution of \((W_{0}^{\prime},(J_{0}^{\prime},b^{\prime}),E_{0}^{\prime})\),
\begin{equation} \label{resol1bis}
\begin{array}{ccccccc}
(W_{0}^{\prime},(J_{0}^{\prime},b^{\prime}),E_{0}^{\prime}) & \longleftarrow &
(W_{1}^{\prime},(J_{1}^{\prime},b^{\prime}),E_{1}^{\prime})& \longleftarrow & \cdots & \longleftarrow &
(W_{k}^{\prime},(J_{k}^{\prime},b),E_{k}^{\prime})\\ \cup & & \cup & & & & \cup \\
\Sing(J_0^{\prime},b^{\prime}) & & \Sing(J_1^{\prime},b^{\prime}) &   & \ldots & &
\Sing(J_k^{\prime},b^{\prime})=\emptyset,
\end{array}
\end{equation}
defined by the  upper semi-continuous functions
\begin{equation*}
f_i^{\prime}: \Sing(J_i,b)\to I_{d}.
\end{equation*}
Then condition (f) asserts that
\begin{equation} \label{equinivel}
f_i(x)=f_i^{\prime}(\Theta_i(x))
\end{equation}
for all \(x\in \Sing(J_i,b)\)  and in particular
\begin{equation*}
\Theta_i: \Max{f_i}\to \Max{f_i}
\end{equation*}
     is an isomorphism. Note that this guarantees that
\begin{equation*}
\Theta_i: (W_{i},(J_{i},b),E_{i})\to (W_{i}^{\prime},(J_{i}^{\prime},b^{\prime}),E_{i}^{\prime})
\end{equation*}
lifts to an isomorphism
\begin{equation*}
{\Theta}_{i+1}: (W_{i+1},(J_{i+1},b),E_{i+1})\longrightarrow
(W^{\prime}_{i+1},(J^{\prime}_{i+1},b^{\prime}) E^{\prime}_{i+1}).
\end{equation*}
Therefore sequences (\ref{resol1}) and (\ref{resol1bis}) are linked as (\ref{equivar2}) is to
(\ref{equivar1}). Furthermore, the same holds for any isomorphism $\Gamma: {\mathcal B} \to {\mathcal
B}^{\prime}$.

\noindent Note that the previous isomorphism \(\Theta_i\) in (\ref{equinivel}) links functions on two
basic objects, both with the same index \(i\).

\end{enumerate}
\end{parrafo}

\begin{remark}
The construction of the \(i\)-th function
\begin{equation*}
f_i:\Sing (J_i,b)\to I_{d}
\end{equation*}
defined by the algorithm depends on the previous functions $ f_j:\Sing (J_j,b)\to I_{d}, $ for \(j\in
\{0,1,\ldots,i-1\}\). This is not explicitly said in the definition of algorithm, but this is only a
technical fact. Note that in any case the previous isomorphism \(\Theta_i\) in (\ref{equinivel}) links
functions on two basic objects, both with the same index \(i\).
\end{remark}

\begin{parrafo}
\label{existence} \textbf{Proof of Theorem~\ref{classical}.}  With the same notation as in Theorem
\ref{classical}, consider the basic object
\begin{equation*}
        (W_0,(J_0,1),E_0),
\end{equation*}
where $W_0=W$, $J_0=\mathcal{I}(X)$ and $E_0= \emptyset$. Clearly $X= \Sing(J_0,1)$.

\

By  Definition \ref{AlgResol} (d) and \ref{propalg} (2),
  the function $$f_0: \Sing(J_0,1)\to (I_d,\leq)$$ is constant
on the restriction of $(W_0,(J_0,1),E_0)$ to $U=W\setminus\Sing(X)$. Let $a(d)$ denote this constant value
along the points in $U \cap X$.

\

By Definition \ref{AlgResol} (c), we know that the algorithm provides a resolution of the basic object
$(W_0,(J_0,1),E_0)$ by means of a finite sequence of blow-ups
\begin{equation}
\label{pruebaclasica} (W_0,(J_0,1),E_0)\longleftarrow(W_1,(J_1,1),E_1)\longleftarrow\ldots\longleftarrow
(W_N,(J_N,1),E_N),
\end{equation}
at permissible centers $Y_i\subset \Sing(J_i,b)$ for $i=0,1,\ldots,N-1$. Therefore, there must be an index
$k\in\{0,1,\ldots,N\}$ such that $\max f_k= a(d)$, and by Definition \ref{AlgResol} (b), there is a
unique index $k$ with this condition.

\

Now  $U$ can be identified with an open set, say $U$ again, of
$W_k$ (note that the centers of the transformations in sequence
(\ref{pruebaclasica}) are defined by $\Max f_i$ and $\max f_i >
a(d)$ for $i<k$).  If $X_k$ denotes the strict transform of $X$ in
$W_k$, then
\begin{equation*}
        X_k \cap U=X \cap U= \Max f_k \cap U.
\end{equation*}
Since $X \cap U=\mbox{Reg}(X)$ is dense in $X$, it follows that $X_k$ is the union of some of the
components of $\Max f_k$,  and hence it is regular and has normal crossings with the exceptional
components by  Definition \ref{AlgResol} (a).  This proves (i) and (ii) of Theorem \ref{classical}.

\

Now it only remains to show that the resolution of singularities
of $X$ that we have achieved is equivariant. An argument similar
to the one   in the proof of Lemma \ref{lifting} shows that if a
group $G$ acts on $W_0$, and $\Theta (Y_0)=Y_0$ for all $\Theta
\in G$, then the  group $G$ acts on the ideal $\mathcal{I}(Y)
\subset \calo_W$, and hence on the blow-up $W_1$. So in order to
lift the action of $G$  to $W_k$, it suffices to lift the action
of all the elements of the group $\Theta \in G$, step by step, to
an isomorphism on $W_k$.  By assumption $\Theta(X_0)=X_0$ so we
can argue as in Example \ref{exiso} to show that each element
$\Theta \in G$ defines an isomorphism of the basic object
$(W_0,(J_0,1),E_0)$ into itself. Hence part (iii) of Theorem
\ref{classical} follows from \ref{propalg} (3). \qed
\end{parrafo}

\begin{parrafo}
\label{proofprincipal} \textbf{Proof of
Theorem~\ref{principalization}.}  With the  notation of Theorem
\ref{principalization}, it is enough to consider the resolution of
the basic object $(W,(I,1),E)$ provided by Theorem
\ref{resolucionob} (see also Remark \ref{so} (1)).

To show that the embedded principalization is equivariant  we argue as in \ref{existence}, where now each
element $\Theta \in G$ defines an isomorphism of the basic object $(W_0,(J_0,1),E_0)=(W,(I,1),E)$ into
itself. \qed
\end{parrafo}




\section{The two main families of equivariant functions.}
\label{examplesfamilies} In Section \ref{algorithm} we have introduced the notion of resolution of basic
objects, and we have shown how such a resolution  can be applied to give a constructive proof of Theorems
\ref{classical} and \ref{principalization}.  A resolution of basic objects is built by defining a
suitable family of equivariant functions $f_B$. In this section we will introduce the two main invariants
used to construct these functions.

\begin{parrafo} \label{ordn}
{\bf The function $\ord$.} Let $B=(W, (J,b),E)$ be a $d-$dimensional basic object.  Then define the
function
\begin{equation*}
\begin{array}{rccl}
\ord^{d}_B: & \Sing(J,b) & \to & {\mathbb Q}\\ & x & \longrightarrow &
     \frac{\nu_J(x)}{b},
\end{array}
\end{equation*}
where \(\nu_J(x)\) denotes the order  of the ideal \(J{\mathcal
O}_{W,x}\) at \({\mathcal O}_{W,x}\). Note that $\ord^{d}_B$ is an
upper semi-continuous function with values in ${\mathbb Q}$.  We
will use the notation \(\ord\)  when the basic object is
understood.
\end{parrafo}

\begin{parrafo} \label{nnn}
{\bf The function $n$.}  Let $B=(W, (J,b),E)$ be a $d-$dimensional basic object, and let $E=\{H_1,\ldots,
H_l\}$. Define
\begin{equation*}
\begin{array}{rccl} n^{d}_B: & \Sing(J,b) & \to & {\mathbb N}\\
     & x & \longrightarrow & n(x)=\sharp\{H_i\subset E: x\in H_i\}.
\end{array}
\end{equation*}
Note that $n^{d}_B$ is an upper semi-continuous function with
values in ${\mathbb N}$.  We will use the notation \(n\)  when the
basic object is understood.
\end{parrafo}

Now we want to study some of the properties of these functions,
for instance, equivariance. To do so, we first need to introduce
some auxiliary definitions:

\begin{parrafo} \label{ExtendTrans}
{\bf The class of extendable transformations.}  Given a basic
object \((W,(J,b),E)\) and an open set \(U\subset W\), consider
the restriction of \((W,(J,b),E)\) to \(U\), \((U,(A,b),E_U) \).
\medskip

Each  sequence of transformations of basic objects over \((W,(J,b),E)\),
\begin{equation}\label{recuerdo}
\begin{array}{ccccccc}
(W,(J,b),E) & \longleftarrow & (W_1,(J_1,b), E_1)& \longleftarrow & \ldots & \longleftarrow & (W_k,
(J_k,b), E_k),
\end{array}
\end{equation}
induces  a sequence  of transformations of basic objects over \((U,(A,b),E_U) \),
\begin{equation}
\begin{array}{ccccccc}\label{recuerdoU}
(U,(A,b),E_U) & \longleftarrow & (U_1,(A_1,b), E_{U_1})& \longleftarrow & \ldots & \longleftarrow & (U_k,
(A_k,b),E_{U_k} ),
\end{array}
\end{equation}
where each \(U_i\) is an open subset of \(W_i\), and \((U_i,(A_i,b), E_{U_i})\) is the restriction of
\((W_i,(J_i,b), E_i)\) to \(U_i\). In fact, sequence (\ref{recuerdoU}) can be defined by setting \(U_i\)
as the pull-back of \(U\subset W\) via \(W \longleftarrow W_i\)  for each index \(i\) in sequence
(\ref{recuerdo}), and by taking the corresponding restriction to \(U_i\). It is natural to ask if the
converse holds:
\medskip

{\em Given a sequence as (\ref{recuerdoU}),  is there a sequence as (\ref{recuerdo}) which is a natural
extension of the first?.}  \medskip

The answer is no. In fact if we consider  a transformation on \((U,(A,b),E_U)\) as in Definition
\ref{transformation} with a smooth center \(Y \subset U\), it might occur that the closure of \(Y\) in
\(W\), \(\overline{Y} \subset W\), is no longer smooth.  However, if we consider transformations over
\((U,(A,b),E_U)\) as the ones introduced in  \ref{projections}, then the extension can be obtained in a
very natural way.
\end{parrafo}

\begin{definition}
\label{exextendable} Let  \(B=(W,(J,b),E)\) be a basic object and
let \(x_0 \in \Sing(J,b)\) be a point. We say that sequence
(\ref{recuerdo}) is {\em \(x_0\)-extendable}  if and only if the
following condition holds:

\

``Whenever \((W_i,(J_i,b), E_i) \longleftarrow (W_{i+1},(J_{i+1},b),E_{i+1}) \) is a permissible
transformation with center \(Y_i \subset W_i\) as in \ref{transformation}, the center \(Y_i\) is mapped
to \(x_0 \in W\) via \(W_i \to W\).''

\

We will denote by \(C_{x_0}(B)\) the class of all
\(x_0\)-extendable sequences over \((W,(J,b),E)\).
\end{definition}

\begin{remark}\label{remarkextendable}
The main property of \(x_0\)-extendable sequences can be stated as
follows: Let $(U,(A,b),E_U)$ be the restriction of \((W,(J,b),E)\)
to $U$. Assume now that \(x_0 \in \Sing(A,b)\), and that the
sequence (\ref{recuerdoU}) is \(x_0\)-extendable. Then condition
of Definition \ref{exextendable} guarantees that we can always
extend sequence (\ref{recuerdoU}) to a sequence (\ref{recuerdo})
so that the former is a restriction of the later.

In particular note that if sequence (\ref{recuerdo}) is
\(x_0\)-extendable, setting \(U=W\setminus\{x_0\}\), the
restricted sequence (\ref{recuerdoU}) induces by composition the
transformation,
\begin{equation*}
(U, (A,b),E_U) \longleftarrow (U_k, (A_k,b),(E_k)_U)
\end{equation*}
which is either the identity map, or a composition of projections  as in \ref{projections}.
\end{remark}

\begin{remark}\label{isoCx}
Let \(\Theta :B=(W,(J,b),E) \to B^{\prime}=(W^{\prime},(J^{\prime},b^{\prime}),E^{\prime})\) be an
isomorphism of basic objects. Fix \(x_0\in \Sing(J,b)\) and set \(x^{\prime}=\theta_0(x_0)\). It is easy
to check that \(\Theta \) defines a natural bijection $$ {\alpha}_{\Theta} :C_{x_0} (B) \to
C_{x^{\prime}_0}(B^{\prime}). $$
\end{remark}

\begin{lemma}
\label{propertiesnord} The upper semi-continuous functions $\ord^{d}_B$ and $n^{d}_B$ verify properties
(e) and (f) of Definition \ref{AlgResol}.
\end{lemma}

\noindent {\em Proof:}  Clearly   if \((U,(A,b),E_U)\) is the restriction of \({\calB}=(W,(J,b),E)\) to
an open set \(U\), and  \(x_0 \in \Sing(A,b)\), the values of both functions at \(x_0\) are independent
of \(U\). This shows that both $\ord^{d}_{\calB}$ and $n^{d}_{\calB}$ satisfy property (e) of Definition
\ref{AlgResol}.  \medskip

The property of equivariance of the family of functions \(n_B^d\) follows from the way isomorphisms of
pairs and basic objects are defined (see Definitions \ref{DefActPar} and \ref{isoBO}).  \medskip

The property of equivariance of the functions \(\ord\) is a key point in our development, and the proof
is based on a very elementary and enlightening idea which we present below (see Appendix \ref{Hironaka}
for a short and formal proof of this fact): \medskip

Let \(B=(W_0,(J_0,b),E_0)\)  be a basic object and consider a sequence
    of transformations of basic objects,
\begin{equation}
\label{resol2}
\begin{array}{ccccccc}
(W_{0},E_{0}) & \longleftarrow & (W_{1},E_{1})& \longleftarrow & \cdots & \longleftarrow &
(W_{k},E_{k})\\ \cup & & \cup & & & & \cup
\\ F_0=\Sing(J_0,b) & & F_1=\Sing(J_1,b) &   & \ldots & &
F_k=\Sing(J_k,b).
\end{array}
\end{equation}
Now fix  a point \(x_0\in \Sing(J_0,b)\subset W_0\). Note that the morphisms \(W_k \to W_0\) maps
\(\Sing(J_i,b) \to \Sing(J_0,b),\) and hence defines  a fiber \(\Sing(J_i,b)_{x_0} \subset \Sing(J_i,b)\)
over \(x_0\).  The point is that the rational number \(\ord^d(x_0)\) can be defined in terms of:
\begin{enumerate}
\item The dimension of $W$, say $d$.
\item The dimension of the fibers \(\Sing(J_i,b)_{x_0}\), whenever
sequence (\ref{resol2}) runs over the \(x_0\)-extendable sequences
(sequences in \(C_{x_0}(B)\)).
\end{enumerate}

This already proves that the family of functions \(\ord^{d}\) is equivariant:  Fix an isomorphism
\(\Theta: B=(W,(J,b),E) \to B'=(W^{\prime},(J^{\prime},b^{\prime}),E^{\prime}),\) and  let
\(x^{\prime}_{0}=\Theta_0(x_0)\in \Sing(J^{\prime},b^{\prime})\). By Definition \ref{isoBO} (ii) a
sequence of transformations over \((W,(J,b),E)\) defines, via \(\Theta\), a sequence over
\((W^{\prime},(J^{\prime},b^{\prime}),E^{\prime})\).  By Remark \ref{isoCx}, this correspondence maps
\(x_0\)-extendable sequences into \(x^{\prime}_0\)-extendable sequences of transformations defining the
bijection \(\alpha_{\Theta}:C_{x_0}(B)\to C_{x^{\prime}_0}(B^{\prime})\) (see Remark \ref{isoCx}). Each
isomorphism \(\Theta_i:\Sing(J_i,b) \cong \Sing(J_i^{\prime},b^{\prime})\) maps fibers to fibers,
\begin{equation}\label{equifibra}
     \Sing(J_i,b)_{x_0}  \cong
     \Sing(J_i^{\prime},b^{\prime})_{x^{\prime}_0};
\end{equation}
     so both fibers have the same dimensions, and hence the rational
     numbers \(\ord^{d}_B(x_0)\) will coincide with
     \(\ord^{d}_{B^{\prime}}(x^{\prime}_0)\). \qed

\begin{remark} \label{remarkexample1}
Since the functions \(f_i(x)\) which appear in Definition \ref{AlgResol} will be defined in terms of the
functions $\ord_B^d$ and $n_B^d$, they will ultimately inherit the properties of equivariance required in
Definition \ref{AlgResol} (e) and compatibility with open restrictions in \ref{propalg} (2).
\end{remark}

\begin{remark}\label{remarkexample2}
\label{allow} The notion of {\em transformation} of basic object defined in \ref{projections}  is
introduced precisely in order to ensure that \(\ord_B^d\) is an equivariant family of functions. As an
example to illustrate this fact, let
\begin{equation*}
W={\mathbb A}^2_k=\mbox{Spec}(k[x,y]),
\end{equation*}
     and set \(J=<x,y>\). If transformations as the ones described in
\ref{projections} were not introduced, then the identity would induce an isomorphism of basic objects
\begin{equation*}
Id:B=(W,(J^5,4),\emptyset)\cong B^{\prime}=(W,(J^5,5), \emptyset ).
\end{equation*}
     Note that for both basic objects the singular locus is the origin
\((0,0)\), and that both are resolvable by one quadratic transformation. Note also that those would be
the only transformations of basic objects if only transformations as the ones introduced in
\ref{transformation} were allowed. Now $$\ord_B^2((0,0))=5/4 \qquad  \text{ and } \qquad
\ord_{B^{\prime}}^2((0,0))= 1,$$ so the family of functions $\ord_B^d$ would not be compatible with the
identity. The law of transformation introduced in \ref{projections} enlarges the class of all possible
transformations, and therefore makes the notion of isomorphism of basic objects stronger. In the case of
our example, it turns out that the identity map is not an isomorphism of basic objects when adding
\ref{projections} as possible law of transformation.
\end{remark}




\part{Applications.}
\label{applications}



\section{Weak  and strict transforms of ideals: Strong Factorizing
   Desingularization.}
\label{weakvsstrict} In this section we discuss the difference between the notions of {\em weak} and {\em
strict transforms of ideals}. As a related matter we establish Theorem \ref{MainThStrict} which says that
given a variety we can    resolve its singularities, and in addition  require a very natural algebraic
condition on the lifting of the ideal defining the variety (see Theorem \ref{MainThStrict} (iii)).  We
finish this section by presenting some  applications of this result.

\begin{definition}
\label{weaksequence} Let $J\subset {\mathcal O}_W$ be a non-zero sheaf of ideals, and let $$W\leftarrow
W_{1}  \leftarrow \ldots  \leftarrow W_k$$ be a sequence of blowing-ups at regular centers.  Then for
each $i\in \{1,\ldots,k\}$ there is an intrinsic factorization $$J{\mathcal O}_{W_i}={\mathcal
L}_i\overline{J}_i,$$ where:
\begin{enumerate}
\item[(a)] The sheaf of ideals ${\mathcal L}_i$ is locally principal
and supported on the exceptional locus.
\item[(b)] The sheaf of ideals $\overline{J}_i$ has no height one
exceptional components: All codimension one components of $V(\overline{J}_i)$  are strict transforms of
the codimension one components of $V(J)$.
\end{enumerate}
With this notation, for $i\in \{1,\ldots, k\}$, $\overline{J}_i$ is
    the  {\em weak transform} of $J$ in $W_i$, $J_i=J{\mathcal O}_{W_i}$ is
    the {\em total transform}  of $J$ in $W_i$, and the {\em strict
    transform} of $J$ in $W_i$, $\widetilde{J}_i$,  is defined by
$$\widetilde{J}_i= \bigcup_{j=1}^{\infty}(\widetilde{J}_{i-1}{\mathcal O}_{W_j}:_{{\mathcal
O}_{W_i}}\cali(H_i)^j),$$ where $H_i$ denotes the exceptional divisor in $W_{i-1}\leftarrow W_{i}$.
Therefore we have the inclusions,
\begin{equation}
\label{inclusionsweak} J_i\subset \overline{J}_i\subset \widetilde{J}_i,
\end{equation}
for $i=1,\ldots,r$.

In general these are strict inclusions, i.e.,  the notions of weak and strict transforms of ideals are
different: For instance, let $X\subset W$ be  a reduced subscheme and let $\pi_r: W_r\longrightarrow W$
be an embedded desingularization of $X$.  Then
\begin{equation}
\label{factormil} \cali(X){\mathcal O}_{W_r}={\mathcal L}_r\cdot\overline{J}_r,
\end{equation}
and unless $X$ is a hypersurface, usually $\overline{J}_r$ does not coincide with
$\cali(X_r)=\widetilde{J}_r$, the ideal of the strict transform of $X$ in $W_r$.  The reason is that,
during the process of desingularization, new primary components of $\overline{J}_r$ arise. This can
clearly be seen by examining the following example:
\end{definition}

\begin{example}
\label{nuevoejemplo} Let $W$ be the real affine space, and let $G$ be the curve parameterized by
$$(t^6,t^9,t^{13}).$$ This curve is defined by the ideal $$\cali(G)=\langle x^3-y^2, x^2y^3-z^3\rangle.$$

Note that $G$ is minimally embedded in ${\mathbb A}_{\mathbb Q}^3$ at the origin $({\bf0,0,0})$, and that
this  is the only singular point of the curve.  Let $$T=V(\langle x^3-y^2\rangle).$$

Since the origin is the only singular point of $G$, to obtain a desingularization by monoidal
transformations  we have to blow-up this point: $\pi_1: W_1\longrightarrow W.$ \medskip

Note that
\begin{equation}
\label{embedded} \cali(G){\mathcal O}_{W_1}=\cali(H_1)^2\overline{J}_1,
\end{equation}
where $H_1\subset W_1$ is the exceptional divisor and $\overline{J}_1$ the weak transform of  $\cali(G)$
in ${\mathcal O}_{W_1}$. Let $T_1 \subset W_1$ be the strict transform of $T$ in $W_1$. Note that
$\cali(T_1) \subset \overline{J}_1$.

It is easy to see that  $\overline{J}_1$   has an embedded component corresponding to the line
\begin{equation}
\label{onedimensional} L=H_1\cap T_1 \subset V(\overline{J}_1).
\end{equation}

We will never eliminate  this primary component $\cali (L)$ of $\overline{J}_1$ by blowing-up at points
of the strict transform of $G$.
\end{example}

\begin{center}
{\bf A Strong Factorizing  Desingularization Theorem.}
\end{center}

In \cite{BV1}  we show   that, given a reduced subscheme,
$X\subset W$,  it is possible to construct a finite sequence of
monoidal transformations, $$W \leftarrow W_1 \leftarrow \ldots
\longrightarrow W_r$$ so that $X_r$ is non-singular, and
    the weak transform   of $J=\cali(X)$ in $W_r$   describes the
    strict transform of $X$ in $W_r$, i.e.  that $\overline{J}_r$ coincides
    with $\cali(X_r)$.  That is,  we present  an  algorithm of
    desingularization which provides in a simple manner  the equations
    describing the desingularization of $X$. An announcement of this
    result has appeared in \cite{BV2}.  \medskip

In fact, our statement   is valid under  milder hypothesis; for instance, we do not need to assume that
$X$ is reduced. The theorem can be stated more precisely as follows:

\begin{parrafo} {\bf  A Strong Factorizing  Desingularization
Theorem.} \label{MainThStrict} \cite[Theorem 1.2]{BV1}
    {\em Let $(W_0, E_0=\{\emptyset\})$ be a pair and let $X=X_0\subset
W_0$ be  a closed subscheme defined by $\cali(X_0) \subset \calo_{W_0}$. Assume that the open set of
regular points  \( \Reg(X) \)  is dense in $X$.  Then there exists a finite sequence of transformations
of pairs,
\begin{equation*}
(W_0,E_0)\longleftarrow \cdots \longleftarrow (W_r,E_r),
\end{equation*}
inducing a proper birational morphism $\pi_r:W_{r}\longrightarrow W_{0}$, so that setting
$E_r=\{H_1,\ldots,H_r\}$,  and letting \( X_{r}\subset W_{r} \) be the strict transform of \(X_0\), we
have that:

\begin{enumerate}
\item[(i)] $X_r$ is regular in $W_r$, and
$W_r\setminus\cup_{i=1}^rH_i\simeq W_0\setminus\mbox{Sing}(X)$. In particular
$\Reg(X)\cong\pi^{-1}_r(\Reg(X))\subset X_{r}$ (via \( \pi_r \) restricted to  $X_{r}$).

\item[(ii)] \( X_{r} \) has normal crossings with \(
E_{r}=\cup_{i=1}^rH_i \) (the exceptional locus of $\pi_r$).

\item[(iii)] The total transform of the ideal $\cali(X_0)\subset
\calo_{W_0}$ factors as a product of ideals in $\calo_{W_r}$: $$ \cali(X)\calo_{W_r}={\mathcal
{L}}\cdot\cali(X_r), $$ where now $\cali(X_r)\subset \calo_{W_r}$ denotes the sheaf of ideals defining
$X_r$, and $$ {\mathcal {L}}=\cali(H_1)^{a_1} \cdot\ldots\cdot\cali(H_r)^{a_r} $$ is an invertible sheaf
of ideals supported on the exceptional locus of $\pi_r$.
\end{enumerate}
}
\end{parrafo}

Parts (i) and (ii) are the classical statement of the theorem of resolution of singularities.  However,
part (iii) is new (see \ref{different} below). Note also  that (i) ensures that, in our algorithm,  the
only points of $W_0$ that will be modified by the morphism $\pi_r$ are the ones in $\Sing(X_0)$.

\begin{parrafo}
\label{different} {\rm {\bf Why is part (iii) new?} \cite[Section 2]{BV1} In Hironaka's line of proof,
the centers of monoidal transformations, chosen in accordance with  the so called {\em standard basis},
are always {\em included in the strict transform} of the scheme. For instance, in Example
\ref{nuevoejemplo}, the first monoidal transformation must be the  the blowing-up at the origin, and any
other center will also have dimension zero. The one dimensional component which appears after this
blowing-up, $H_1\cap T_1$ (see (\ref{onedimensional})), is a primary component of $\overline{J}_1$, and
will never be eliminated by blowing-up at centers supported on the singular locus of $G$.  Hence (iii)
will never hold for  desingularizations of this curve that follow from  Hironaka's proof.  \medskip

In order to achieve (iii) one must blow-up $H_1\cap T_1$ (or some strict transform of it). The new
algorithm that we propose, first considers the quadratic transformation $\pi: W_1\longrightarrow W_0$,
and, some steps later, the  blowing-up at the strict transform of the one dimensional scheme $H_1\cap
T_1$. Since $H_1\cap T_1$ is mapped to the singular locus,  the first isomorphism in Theorem
\ref{MainThStrict} (i) is preserved after such monoidal transformation.  \medskip

We think of a subscheme $X$ of a smooth scheme $W$, at least locally, as a finite number of {\em
equations} defining the ideal $\cali(X)$. An  algorithm of desingularization should  provide us with:

\begin{enumerate}
\item[(1)] A sequence of monoidal transformations over the smooth
scheme $W$, $$W_0=W \leftarrow  W_{1}  \leftarrow \ldots  \leftarrow W_n$$ so that conditions (i) and
(ii)  Theorem \ref{MainThStrict} hold for the strict transform of $X$ at $W_n$.

\item[(2)] A pattern of manipulation of equations defining $X$, so as
to obtain, at least locally at an open covering of $W_n$, equations defining the strict transform  $X_n$
of $X$  at $W_n$.
\end{enumerate}

So (2) indicates how the original equations defining $X$ have to be treated, at an affine open subset of
$W_n$, in order to obtain local equations defining $X_n$. While this is very complicated in Hironaka's
line of proof, here it is a  direct consequence of Theorem \ref{MainThStrict} (iii). In fact, for
algorithms that follow Hironaka's proof, to get   both  (1) and (2) one  must consider the {\em strict
transform} of the ideal of the subscheme at each monoidal transformation. In that setting one has to
choose  a {\em standard basis} of the ideal, which is a system of generators of the ideal  of the
subscheme, but such choice of generators must be changed if the maximum Hilbert Samuel invariant drops
in the sequence  of monoidal transformations. All of  these complications are  avoided in our proof,
which simplifies both (1) and (2).  }
\end{parrafo}

\begin{center}
{\bf The idea of the proof of Theorem \ref{MainThStrict}: \\ The notion of ``excess of embedding
codimension'' of an ideal.}
\end{center}

The idea behind   the proof of Theorem \ref{MainThStrict} is to
     construct  a sequence of monoidal transformations
$$W \leftarrow W_1 \leftarrow \ldots \leftarrow  W_r, $$ so that the exceptional components of the total
transform of $J$ in $W_r$, $J{\mathcal O}_{W_r}$,   become locally principal.

\

Thus, the strategy  is to show first that for a suitable sequence of monoidal transformations
    $$ W_0 \leftarrow W_{1}\leftarrow \ldots\leftarrow  W_l,
$$ the subscheme $V((\overline{J}_l)_y)$ is locally included in a smooth hypersurface, and this condition
holds at every $y\in V(\overline{J}_l)$. Then we show that by applying more monoidal transformations $$
W_l\leftarrow    W_{l+1}\leftarrow \ldots\leftarrow W_m,$$ the closed subscheme  $V((\overline{J}_m)_y)$
is  locally included in a smooth subscheme of codimension two, for every $y\in V(\overline{J}_m)$. To
obtain our result we repeat this process as many times as needed.

\

This strategy  leads us to the notion of {\em excess of embedding codimension of an ideal}:

\begin{definition}\cite[Definition 5.2]{BV1}
\label{omni} {\rm Let $(W,E)$ be a pair, and let $J\subset {\mathcal O}_{W}$ be a non-zero sheaf of
ideals.

\begin{enumerate}
\item[(a)] We  say that $J$ has {\em excess of embedding codimension
    $\geq a$ in (W,E)  at a point} $y\in W$, if either $J_y={{\mathcal
    O}}_{W,y}$, or there is a regular system of parameters $\{x_1,x_2,
    \ldots , x_d \}\subset {{\mathcal O}}_{W,y}$ such  that

\begin{enumerate}
\item[(i)] $<x_1,x_2, \ldots , x_a>\subset J_y \subset {{\mathcal
O}}_{W,y}$, and

\item[(ii)] every  hypersurface $H_i\in E$ containing the point $y$
has a local equation $$\cali(H_i)=<x_{i_j}> \subset {{\mathcal O}}_{W,y}$$ with $i_j > a$.
\end{enumerate}

\item[(b)] We  say that $J$ has {\em excess of embedding codimension }
$\geq a$ in $(W,E)$, if the  conditions stated in (a)  hold at every point $y \in W$. We will abbreviate
this by saying that $J$ has {\em excess of codimension } $\geq a$ in $(W,E)$, or {\em $(W,E)$-codimension
$\geq a$}. Reference to the pair $(W,E)$ will be omitted if it is clear from the context.
\end{enumerate}
}
\end{definition}

\begin{remark}\label{rem2} Note
that any non-zero ideal $I \subset {{\mathcal O}}_{W}$ has excess
codimension $\geq 0$ in $(W,E)$, since such condition is vacuously
satisfied in that case.
\end{remark}

\begin{remark}\label{rem3}
If $X\subset W$ is a regular subscheme of pure codimension $e$, then the sheaf of ideals $\cali(X)\subset
{\mathcal O}_{W}$ has excess of codimension $\geq e$ in $(W,E=\emptyset)$.
\end{remark}

The main result over which the proof of Theorem \ref{MainThStrict}
is based is  \cite[Lemma 4.7]{BV1}. What we state below is a
slightly stronger version of this lemma:

\begin{lemma}
\label{plemma} Let $(W,E)$  be a pair, and let
    $J$ be a non-zero sheaf of ideals. Assume that there is a non-empty
open set $U^a\subset W\setminus E$ such that $V(J)\cap U^a\neq \emptyset$ and that $J$ has $(W,E)-$excess
of codimension  $\geq a$ at any point of $U^a$. Then there is a finite sequence of transformation of
pairs:
\begin{equation}
\label{2eq200000} (W_k,E_k)\longrightarrow  \ldots\longrightarrow (W_1,E_1)\longrightarrow
(W_0,E_0)=(W,E),
\end{equation}
such that if
    $$J{\mathcal O}_{W_k}=\call_{k} \overline{J}_{k},$$ then:
\begin{enumerate}
\item[(i)] The weak transform of $J$ in $W_k$,  $\overline{J}_{k}$,
has $(W_k,E_k)-$excess of  codimension $\geq a$.
\item[(ii)] The birational morphism $W_k \longrightarrow W$ induces an
isomorphism over the open set $U^a\subset W$.
\item[(iii)] If in addition we can express  $E=\{H_1,\ldots,H_k\}$
as a disjoint union, $E=E^a\sqcup E^{a^{\prime}}$, so that $J$ has
excess of embedding codimension $\geq a$ at any point of $V(J)\cap
E^a$, then  the sequence (\ref{2eq200000}) can be constructed so
that (ii) also  holds if we only assume that $U^a\subset
W\setminus E^{a^{\prime}}$.
\end{enumerate}
\end{lemma}

\begin{proof} Parts (i) and (ii) have been  proven in \cite{BV1}.
For part (iii) it is enough to observe that if we set $E_0^+=E^a$ and $E_0^-=E^{a^{\prime}}$, then the
morphism $W_k \longrightarrow W$ induces an isomorphism on $E\setminus E^{a^{\prime}}$ (see Definition
\ref{nbo}, where $E_0^+$ and $E_0^-$ are introduced, and see also \cite[Definition 5.9]{BV1}).
\end{proof}

The following is an easy consequence of Lemma \ref{plemma}:

\begin{corollary}
\label{corolario} Fix a pair $(W,E)$. Let $X\subset W$ be a closed subscheme, let $J=\cali(X)$, and let
$$X=X(1)\cup X(2)\cup\ldots\cup X(e),$$ where for $i=1,\ldots,e$, each  $X(i)$  is a closed subscheme of
pure codimension $i\leq n$. Assume that $\Reg(X(e))\neq \emptyset$ and let $Z$ be the closure of
$\Reg(X(e))$ in $W$.  Suppose that $E=E^e\sqcup E^{e^{\prime}}$ and that $Z$ has normal crossings with
$E^e$ at any point of $Z\cap E^e$ ($E^e$ might be empty). Let ${\mathcal Q}\subset W$ be the open set
$W\setminus \left[X(1)\cup\ldots\cup X(e-1)\cup \Sing(X(e))\cup E^{e^{\prime}}\right]$.  Then there is a
finite sequence of transformations of pairs,
\begin{equation}
\label{2eq222} (W_k,E_k)\longrightarrow  \ldots\longrightarrow (W_1,E_1)\longrightarrow (W_0,E_0)=(W,E)
\end{equation}
such that:
\begin{enumerate}
\item[(i)] The strict transform of $Z$ in $W_k$, $Z_k$,  is regular,
and has normal crossings with $E_k$.
\item[(ii)] The morphism $W_k\to W$ induces an isomorphism on ${\mathcal Q}$.
\item[(iii)] The weak transform of $J$ in $W_k$,  $\overline{J}_{k}$,
coincides with $\cali(Z_k)$, where $Z_k$ is the  strict transform of $Z$ in $W_k$, i.e. $$J{\mathcal
O}_{W_k}=\call_k\cali(Z_k),$$ where $\call_k$ is a locally principal sheaf of ideals supported on the
exceptional locus of $W_k\to W$.
\end{enumerate}

\end{corollary}

\noindent{\em Proof:} By hypothesis, $J$ has $(W,E)-$excess of embedding  codimension $\geq e$ in
${\mathcal Q}$. Then by Lemma \ref{plemma}, there is a finite sequence of  transformations of pairs,
\begin{equation}
\label{2eq211} (W_k,E_k)\longrightarrow  \ldots\longrightarrow (W_1,E_1)\longrightarrow (W_0,E_0)=(W,E)
\end{equation}
such that:
\begin{enumerate}
\item[(a)] The weak transform of $J$ in $W_k$, $\overline{J}_k$, has
excess of codimension $e$ in $(W_k,E_k)$.
\item[(b)] The morphism $W_k\to W$ induces an isomorphism on ${\mathcal Q}$.
\end{enumerate}

Now write $$V(\overline{J}_k)=Y_1\cup\ldots \cup Y_m$$ where each $Y_i$ is an irreducible component of
codimension at least $e$.  Since $Z_k\subset V(\overline{J}_k)$, without lost of generality, we may
assume that $$Z_k=Y_1\cup\ldots\cup Y_s,$$ where $1\leq s\leq m$. Note that since $\overline{J}_k$ has
excess of embedding codimension $\geq e$ in $(W_k,E_k)$, $Z_k=Y_1\cup\ldots\cup Y_s$ is regular.  Now, if
$s=m$ we are done. On the other hand, if $s<m$, then  $V(\overline{J}_k)$ is a disjoint union of
$Y_1\cup\ldots\cup Y_s$ and $Y_{s+1}\cup\ldots\cup Y_m$, i.e.,
\begin{equation}
\label{disjoint} V(\overline{J}_k)=\left[Y_1\cup\ldots\cup Y_s\right]\sqcup \left[Y_{s+1}\cup\ldots\cup
Y_m\right]
\end{equation}
(here we use the fact that $Z_k=Y_1\cup\ldots\cup Y_s$ is regular and has codimension $e$: If
(\ref{disjoint})  were not a disjoint union, at any point of the intersection the codimension would be
greater than $e$).  As a consequence $$ \overline{J}_k=\cali(Y_1\cup\ldots\cup
Y_s)\cdot\cali(Y_{s+1}\cup\ldots\cup Y_m). $$ By Theorem \ref{principalization}, there is an enlargement
of sequence \ref{2eq211},
\begin{equation}
\label{2eq2111} (W_k,E_k)\leftarrow  \ldots\leftarrow (W_l,E_l)
\end{equation}
which defines a strong principalization of $\cali(Y_{s+1}\cup\ldots\cup Y_m)\subset {\mathcal O}_{W_k}$.
Therefore $$J{\mathcal O}_{W_l}=\call_l \overline{J}_l$$ and $$\overline{J}_l=\cali(Z_l).$$ Note that
since (\ref{disjoint}) is a disjoint union, $W_l\to W$ defines an isomorphism on ${\mathcal Q}$. \qed

\

\begin{center}
{\bf An application of Corollary \ref{corolario}:   Log resolutions of divisors.}
\end{center}

\begin{definition}  Let $D\subset W$ be a divisor. A {\em log
resolution} of the pair $(W,D)$ is a  proper and birational morphism, $$\mu: W^{'}\to W$$ such that
$W^{'}$ is non-singular and $\mu^*(D)+\mbox{Exc}(\mu)$  is a normal crossing divisor.
\end{definition}

As an application of Corollary \ref{corolario} we prove the following statement:

\begin{theorem}
\label{ese} Let $D=\sum_{i=1}^{r}a_iD_i\subset W$ be a divisor, and let $V\subset W$ be the open set of
points where $D$ has normal crossings.  Then there is a log resolution of $D$ $$\mu :W^{\prime}\to W$$
which induces an isomorphism on $V$.
\end{theorem}

\begin{parrafo}
\label{motivation} {\rm {\bf The motivation.}  The interest of
this result,  lies in its role in proving statements which involve
the application of vanishing theorems in the compactification of a
given variety. More explicitly, assume that  $W$ is a smooth
quasi-projective variety and that $D\subset W$ is a normal
crossing divisor. Some vanishing theorems are only stated for
projective varieties.  Hence, if $W$ is not projective, there is a
compactification of $W$, $\overline{W}$, where the vanishing
theorems hold. However, it may happen that the compactification of
$D$ in $\overline{W}$, $\overline{D}$,  is not  a normal crossing
divisor any more.  By applying a log resolution to
$(\overline{W},\overline{D})$,
$$\overline{\mu}: \overline{W}^{\prime}\to \overline{W}$$ we may assume that
$$\overline{\mu}^{*}(\overline{D})$$ is a  divisor with normal crossings. For  technical reasons, it may be
useful to assume that the induced morphism $$\mu:{W}^{\prime}\to {W}$$ is an isomorphism over $D\subset
W$. To this end, it would be enough if $$\overline{\mu}: \overline{W}^{\prime}\to \overline{W}$$ is an
isomorphism over the points where $\overline{D}$ is already a normal crossing divisor. For a concrete
example see \cite[Remark 9.4.2]{laz}.}
\end{parrafo}

Before proving  Theorem \ref{ese}  we will introduce some notation
and prove some auxiliary results (Proposition \ref{erre}).  The
proof of Theorem \ref{ese}, which is given  in \ref{proofofthe},
will follow from Proposition \ref{erre} and an inductive argument.

\begin{remark}
\label{laformula} In what follows, we will denote by  $D=\sum_{i=1}^ra_iD_i$  a divisor in $W$, where
each $D_i$ corresponds to an irreducible and reduced hypersurface in $W$. Let
  $J={\mathcal O}_W(-D_1)+\ldots+ {\mathcal O}_W(-D_r)$ be the ideal
corresponding to the intersection of all the components of $D$. If $$W_k\to W_{k-1}\to \ldots\to
W_1\longrightarrow W$$ is a finite sequence of monoidal transformations, and if $\pi_k:W_k\to W$ denotes
the composition, then we will write $$ \pi_k^*(D)=\sum_{i=1}^ra_iD_i^{[k]}+F_k=D^{[k]}+F_k, $$ where
$D_i^{[k]}\subset W_k$ denotes the strict transform of $D_i$ in $W_k$ for $i=1,\ldots,r$, and $F_k$ is a
divisor supported in the exceptional locus.  Let  $J^{[k]}$ be the sheaf of ideals ${\mathcal
O}_{W_k}(-D_1^{[k]})+\ldots+ {\mathcal O}_{W_k}(-D_r^{[k]})$ which describes the intersection
$D_1^{[k]}\cap \ldots\cap D_r^{[k]}$ in $W_k$.  Note that if $\overline{J}_k$ is the weak transform of
$J={\mathcal O}_W(-D_1)+\ldots+ {\mathcal O}_W(-D_r)$ in ${\mathcal O}_{W_k}$, then
\begin{equation}
\label{inequality1} \overline{J}_k\subset J^{[k]}.
\end{equation}
This containment   can be proven by an inductive argument: Let
$C\subset W$ be a center such that $\nu_{C}({\calo}_W(-D_i))=a_i$,
for $i=1,\ldots,r$ (here $\nu_C({\calo}_W(-D_i))$ denotes the
order of  the ideal ${\calo}(-D_i)$ at $C$). Then
$\nu_{C}(J)=b=\mbox{Min} \{a_1,\ldots,a_r\}$.  Let $W\leftarrow
W_1$ be the monoidal transformation with center $C$. Then, if
$H_1$ denotes the exceptional divisor,
$$\overline{J}_1=(J{\calo}_{W_1}:\cali(H_1)^b) \qquad \mbox{ and }
\qquad J^{[1]}={\mathcal O}_{W_1}(-D_1^{[1]})+\ldots+ {\mathcal
O}_{W_1}(-D_r^{[1]}),$$ therefore,
$$\overline{J}_1\subset J^{[1]}  \qquad \mbox{ and } \qquad
    \overline{J}_1={\mathcal M}_{1,1}{\mathcal O}_{W_1}(-D_1^{[1]})+\ldots+
    {\mathcal M}_{r,1}{\mathcal O}_{W_1}(-D_r^{[1]}),$$ for some ideals
    ${\mathcal M}_{i,1}\subset {\calo}_{W_1}$, $i=1,\ldots,r$.

\

Now, assume by inductive hypothesis that after a finite sequence of blowing-ups, $$W \leftarrow W_1
\leftarrow \ldots  \leftarrow W_{k}$$ at smooth centers $C_i\subset \overline{J}_i$, $i=1,\ldots,k-1$, we
have that $$\overline{J}_k\subset J^{[k]} \qquad \mbox{ and } \qquad \overline{J}_k={\mathcal
M}_{1,k}{\mathcal O}_{W_k}(-D_1^{[k]})+\ldots+ {\mathcal M}_{r,k}{\mathcal O}_{W_k}(-D_r^{[k]}),$$ for
some ideals ${\mathcal M}_{i,k}\subset {\calo}_{W_k}$.

\

Let $C_k\subset V(\overline{J}_k)\subset W_k$ be a smooth center, assume  that
$\nu_{C_k}({\calo}_{W_k}(-D_i^{[k]}))=a_{i,k}$, for $i=1,\ldots,r$ and that $\nu_{C}
(\overline{J}_{k})=b_k$. Note that $b_k=\mbox{Min } \{a_{1,k}+\nu_{C_k}({\mathcal M}_{1,k}),\ldots,
a_{r,k}+\nu_{C_k}({\mathcal M}_{r,k})\}$.  Then if $W_k\leftarrow W_{k+1}$ is the  blowing-up at $C_k$
and if $H_{k+1}$ denotes the exceptional divisor, $$\overline{J}_{k+1}=(\overline{J}_k{\calo}_{W_{k+1}} :
\cali(H_{k+1})^{b_k}) \qquad \mbox{ and } J^{[k+1]}={\mathcal O}_{W_{k+1}}(-D_1^{[k+1]})+\ldots+ {\mathcal
O}_{W_{k+1}}(-D_r^{[k+1]}).$$ Now note that $$\overline{J}_{k+1}=(\overline{J}_k{\calo}_{W_{k+1}} :
\cali(H_{k+1})^{b_k})=$$ $$=\left(\left[{\mathcal M}_{1,k}{\mathcal O}_{W_k}(-D_1^{[k]})+\ldots+
{\mathcal M}_{r,k}{\mathcal O}_{W_k}(-D_r^{[k]})\right]{\calo}_{W_{k+1}} : \cali(H_{k+1})^{b_k}\right).$$
Hence, $$\overline{J}_{k+1}\subset J^{[k+1]}, \qquad \mbox{ and } \qquad \overline{J}_{k+1}={\mathcal
M}_{1,k+1}{\mathcal O}_{W_{k+1}}(-D_1^{[k+1]})+\ldots+ {\mathcal M}_{r,k+1}{\mathcal
O}_{W_{k+1}}(-D_r^{[k+1]}),$$ for some ideals ${\mathcal M}_{i,k+1}\subset {\calo}_{W_{k+1}}$,
$i=1,\ldots,r$. \qed

\end{remark}

\begin{definition}
\label{snc} Let $(W,E)$ be a pair, and let $D=\sum_{i=1}^ra_iD_i\subset W$ be a divisor, where each $D_i
$ is  irreducible and reduced. Let $s\in \{1,\ldots,r\}$.
\begin{itemize}
\item We will say that  $D$ has {\em $s-$normal crossings in
$(W,E)$}  at a point $x\in W$ if:
\end{itemize}
\begin{enumerate}
\item[(a)] There are fewer than $s$ components of
$D$ containing $x$; or else,
\item[(b)] The divisors  $\{\cup_{i=1}^r D\cup E\}$
have normal crossings at $x$.
\end{enumerate}
\begin{itemize}
\item If $D$ does not
have $s-$normal crossings in $(W,E)$ at $x$  for any
$s\in\{1,\ldots,r\}$, then we  will say that $D$ has {\em
$r+1-$normal crossings} at $x$ in $(W,E)$, since this condition is
always vacuously satisfied.
\item If we  say that $D$ has $s-$normal crossings in $(W,E)$
  we will understand that it
  has $t-$normal crossings in $(W,E)$ at any  point of $W$,
with $t\geq s$.
\end{itemize}
\end{definition}

\begin{remark}
\label{normalcrossings} Note that:
\begin{enumerate}
\item[(i)] If a divisor $D$ has normal crossings  then it has normal
crossings at any point of $W$ and hence it has $1-$normal crossings in $(W,\emptyset)$ at any point $x\in
W$.
\item[(ii)] If $D$ has $1-$normal crossings in $(W,E)$ at any point
   $x\in W$, then $D$ has
normal crossings, and $\{D\cup E\}$ have normal crossings.
\end{enumerate}
\end{remark}

\begin{proposition}
\label{erre} Let $(W,E)$ be a pair, let $D=\sum_{i=1}^{s}a_iD_i$ be a divisor, and assume that
$E=E^s\sqcup E^{s^{\prime}}$ is a disjoint union such that $D$ has $s-$normal crossings at any point of
$E^s$ ($E^s$ might be empty). Let
  $U$ be the subset of $W\setminus E^{s^{\prime}}$  where $D$ has normal
crossings, and let $U^{s}$ be the subset of $W\setminus E^{s^{\prime}}$
  where $D$ has $s-$normal crossings.
Then there is a finite sequence of transformations of pairs, $$(W_0,E_0)=(W,E) \leftarrow (W_1,E_1)
\leftarrow \ldots\leftarrow (W_l,E_l),$$ such that
\begin{enumerate}
\item[(i)] The divisor ${D}^{[l]}=\sum_{i=1}^sD_i^{[l]}$ has
$s-$normal crossings in $(W_l,E_l)$.
\item[(ii)] The morphism $\pi_l:W_l\to W$ induces an isomorphism on
$U^{s}\subset W$ (and hence on $U$).
\end{enumerate}
\end{proposition}

\

\noindent{\em Proof:} First note that by Definition \ref{snc} (b), $U\subset U^{s}$. We will distinguish
three  cases:

\

\noindent {\em Case 1:} If $D_1\cap \ldots\cap D_s=\emptyset$ then there is nothing to prove.

\

\noindent {\em Case 2:} If there  are no points in $D_1\cap \ldots\cap D_s$ where $\{D_1,\ldots, D_s,E\}$
have normal crossings, or if $D_1\cap \ldots\cap D_s\subset E^{s^{\prime}}$, then applying  Theorem
\ref{principalization} to $J={\mathcal O}_{W}(-D_1)+\ldots+{\mathcal O}_{W}(-D_s)$, there is a finite
sequence of blowing-ups at regular centers, $$(W_0,E_0)=(W,E)\leftarrow (W_1,E_1) \leftarrow \ldots
\leftarrow (W_l,E_l),$$ such that the weak transform of $J$ in ${\mathcal O}_{W_l}$, $\overline{J}_l$,
is trivial. Since, $${\mathcal O}_{W_l}=\overline{J}_l\subset J^{[l]},$$ (see Remark \ref{laformula}), we
have that $$D_1^{[l]}\cap \ldots \cap D_s^{[l]}=\emptyset.$$ Note that in this case $\pi_l:W_l\to W$
induces an isomorphism on $W\setminus V(J)$, so in particular it induces an isomorphism  on $U^s$ and
hence on $U$.

\

\noindent {\em Case 3:} As in the previous case, set $J={\mathcal O}_W(-D_1)+\ldots+{\mathcal
O}_W(-D_s)$.  If there is a point $x\in D_1\cap \ldots\cap D_s$ where $\{D_1,\ldots, D_s,E\}$ have normal
crossings and $x\notin E^{s^{\prime}}$,
  then  $J$ has excess of codimension
$\geq s$ at $x$,  and in fact $J$ has excess of codimension $\geq s$ at any point of $V(J)\cap U^{s}$.
Let $Z$ be the closure of $U^{s}\cap V(J)$ in $W$.  By  Corollary
  \ref{corolario} applied to $X=V(J)$ there is a finite sequence of
transformation of pairs $$(W_0,E_0)=(W,E)\leftarrow (W_1,E_1) \leftarrow \ldots \leftarrow (W_l,E_l),$$
such that:
\begin{enumerate}
\item[(a)] The strict transform of $Z$ in $W_l$, $Z_l$, is smooth and
has normal crossings with $E_l$.
\item[(b)] The morphism $W_l\to W$ induces an isomorphism on
   $U^{s}$, and hence on $U$.
\item[(c)] The ideal $\overline{J}_l$ coincides with $\cali(Z_l)$.
\end{enumerate}
Therefore $\{D^{[l]},E_l\}$ have normal crossings at any point of $Z_l=D_1^{[l]}\cap \ldots\cap
D_s^{[l]}$. \qed

\begin{parrafo}
\label{proofofthe}{\rm  {\em Proof of Theorem \ref{ese}:} Assume that $D$ has $t-$normal crossings in
$(W,E)$, with $1<t\leq r+1$, but that it does not have $(t-1)-$normal crossings in $(W,E)$. Let $\Sigma$
be the set of all possible $(t-1)-$tuples of elements in $\{1,\ldots,r\}$. If $\sigma\in \Sigma$, then we
will denote by $D_{\sigma}$ the divisor $\sum_{i\in \sigma}D_i$, by $J_{\sigma}$ the ideal $\sum_{i\in
\sigma} {\mathcal O}_{W}(-D_i)$, and by $G_{\sigma}$ the closed set of points in $V(J_{\sigma})$ where
$J_{\sigma}$ does not have excess of embedding codimension $\geq t-1$.

\

Fix an order in $\Sigma$  such that if $$\Sigma=\{\sigma_1,\ldots, \sigma_{m}\},$$ with
$$\sigma_1<\ldots<\sigma_{m-1}<\sigma_{m},$$ then every  divisor in $${\mathcal
G}=\{D_{\sigma_1},\ldots,D_{\sigma_{l-1}}\}$$ has  $(t-1)-$normal crossings in $(W,E)$, and each of the
divisors $$D_{\sigma_l},\ldots,D_{\sigma_m}$$ has only $t-$normal crossings in $(W,E)$ (note that
${\mathcal G}$ might be empty). Since by assumption $D$ has $t-$normal crossings in $(W,E)$, note that:
\begin{enumerate}
\item[(a)] For every $j< l$, $G_{\sigma_j}=\emptyset$, and
for every $j\geq l$, locally, at any closed point $x\in G_{\sigma_j}$,  the divisor $D=\sum_{i=1}^rD_i$
can be identified with $D_{\sigma_j}$.
\item[(b)] If  $E_{\sigma_l}=E\cup \{\cup D_i\}$,
with $i\in \{1,\ldots,r\}\setminus \sigma_l$, then there is an
open neighborhood ${\mathcal P}$ of $V(J_{\sigma_l})$ such that
$(W,E_{\sigma_l})$ is a pair.
\end{enumerate}

Now we apply Proposition \ref{erre}  to $D_{\sigma_l}$,
$({\mathcal P},E_{\sigma_l}=E\cup \{\cup D_i\})$, setting $s=t-1$,
$E^{t-1}=\{\cup D_i\}$ with $i\in \{1,\ldots,r\}\setminus
\sigma_l$,
  $E^{(t-1)\prime}=E$, and restricting the open sets,
$V\cap {\mathcal P}$, and  $V^{t-1}\cap {\mathcal P}$ (where $V^{t-1}$ is the subset of
  $W\setminus
E$ where $D_{\sigma_l}$ has $(t-1)-$normal crossings). Thus, there exists a finite sequence of
transformations of pairs,
\begin{equation}
\label{sigma1} ({\mathcal P}_0,(E_{\sigma_l})_0)=({\mathcal P},E_{\sigma_l}) \leftarrow ({\mathcal
P}_1,(E_{\sigma_l})_1) \leftarrow \ldots \leftarrow ({\mathcal P}_{N},(E_{\sigma_l})_{N}),
\end{equation}
such that:
\begin{enumerate}
\item The divisor
$D^{[N]}_{\sigma_l}= \sum_{i\in \sigma_l}D_i^{[N]}$ has $(t-1)-$normal crossings in $({\mathcal
P}_{N},(E_{\sigma_l})_{N})$ (here we use the fact that the centers where we blow-up in sequence
(\ref{sigma1}) have normal crossings with $\{\cup D_i\}_{i\notin \sigma_l}$, so for $i=1,\ldots,l-1$,
each $D_{\sigma_i}^{[N]}$ has still $t-1$ normal crossings in $(W_N,E_N)$).
\item The morphism $\pi_{N}:{\mathcal P}_{N}\to {\mathcal P}$
induces an isomorphism on
   $V^{t-1}\cap U$, and hence on $V\cap U$.
\end{enumerate}

Now observe that:
\begin{enumerate}
\item[(A)] Sequence (\ref{sigma1}) induces a sequence of transformations
   of pairs,
$$(W_0,E_0)=(W,E) \leftarrow (W_1,E_1) \leftarrow \ldots \leftarrow (W_{N},E_{N}).$$
\item[(B)] The divisors
   $\{D^{[N]}_{\sigma_1},\ldots,D^{[N]}_{\sigma_{l}}\}$
have $(t-1)-$normal crossings in $(W_N,E_N)$.
\item[(C)] The morphism $W_N\to W$ induces an isomorphism on
$V$.
\end{enumerate}

Iterating this argument, there is a finite sequence of transformations of pairs inducing an isomorphism
on $V$, $$(W_0,E_0)=(W,E) \leftarrow (W_1,E_1) \leftarrow \ldots \leftarrow (W_{M},E_{M}),$$ such that
for each  $\sigma_i\in \Sigma$, all divisors $D_{\sigma_i}^{[M]}$ have $(t-1)-$normal crossings in
$(W_{M},E_{M})$, and hence $D^{[M]}$ has $(t-1)-$normal crossings in $(W_{M},E_{M})$. Therefore,
applying an inductive argument, there is a finite sequence of transformations of pairs, $$(W_0,E_0)=(W,E)
\leftarrow (W_1,E_1) \leftarrow \ldots \leftarrow (W_{L},E_{L}),$$ such that $\pi_L: W_L\to W$  induces
an isomorphism on $V$, and  $D^{[L]}$ has $1-$normal crossings in $(W_{L},E_{L})$. By Remark
\ref{normalcrossings} (ii), this means that $\pi_L^{*}(D)=D^{[L]}+F_L$ has normal crossings, since
$F_L\subset E_L$. \qed }

\end{parrafo}



\section{On a class of regular schemes and on real and complex analytic
spaces.} \label{completions} Up to now we have worked with schemes
of finite type over fields. However, the constructions and proofs
of Theorems \ref{classical} and \ref{principalization} are valid
for a much wider class of schemes (not necessarily of finite type
over a field $\mathbf{k}$). Here we will focus on a class to which
our techniques naturally extend.

\begin{definition}\label{laclasescero}
We define ${\mathcal S}_0$ as the class of regular,
equidimensional schemes $ W $ containing a field, $\mathbf{k}$  of
characteristic 0 (which may vary), satisfying the following two
conditions:
\begin{enumerate}
\item[(i)] If $ W $ is an $n-$dimensional $\mathbf{k}$-scheme in ${\mathcal
     S}_0$, then there is  a finite affine open covering $\{{\mathcal
     U}_i\}_{i\in I}$, of $W$ such that for each $i\in I$,  ${\mathcal
     U}_i \approx \mbox{Spec}(R_i)$, for  some Noetherian, regular
     $\mathbf{k}$-algebra, $R_i$ with the additional property that
     $\mbox{Der}_{\mathbf{k}}(R_i)$ is a finite projective
$R_i$-module, locally
     of rank $n$.
\item[(ii)]  If $\mathfrak{m}$ is a maximal ideal in $ R_i$ then dim
$(  R_i)_{\mathfrak{m}}=$n and $ R_i/\mathfrak{m}$ is algebraic over $ \mathbf{k} $.
\end{enumerate}
\end{definition}

Note that under  condition (ii),   $\mathbf{k}$ is a quasi-coefficient field at the localization at any
closed point in the sense of \cite[p. 274]{Matsumura1980}.

\begin{remark} Note that:
\begin{enumerate}
\item[(i)]  Any smooth scheme over a field of characteristic zero, as
well as the spectrum of the completion, or henselization,  of a
local ring which is in  ${\mathcal S}_0$ is in ${\mathcal S}_0$
(see \cite[Theorem 90]{Matsumura1980}).
\item[(ii)] The class of schemes in ${\mathcal S}_0$ is closed
under monoidal transformations. This fact follows from
\cite[Appendix 40, Theorem 99, (3), (4)]{Matsumura1980}.
\end{enumerate}
\end{remark}

\begin{theorem}
The algorithms of desingularization and of principalization stated in Theorems \ref{classical} and
\ref{principalization} extend to the class of schemes in ${\mathcal S}_0$.
\end{theorem}

The reason why all the constructions that we have made for schemes
of finite type over a field of characteristic zero extend to the
class of schemes in ${\mathcal S}_0$,  is that this class of
schemes satisfies Property D stated below in Lemma
\ref{propertyd}: The algorithm of resolution  is based in an
inductive argument which in turn is based  on the nice properties
of the differential operator $\Delta$ defined on  the class of
smooth schemes over a field (see Sections \ref{SeccDelta} and
\ref{proposition}).
  The class ${\mathcal S}_0$ parallels the class of
compact (real or complex) analytic spaces with a well defined
tangent bundle. It is easy to check that our results and
developments extend to these class of spaces, too.

\begin{definition}
Let $R$ be a Noetherian, regular $\mathbf{k}$-algebra, such that $\mbox{Der}_k(R)$ is a finite projective
$R$-module, locally of rank $n$, and assume that for each maximal ideal $\mathfrak{m}$ of $R$, $
R_j/\mathfrak{m}$ is algebraic over $\mathbf{k} $. Then for each ideal $J$ of $R$ we define $\Delta(J)$
to be the ideal generated by $$\{\delta(f): \delta\in \mbox{Der}_{\mathbf{k}}(R), f\in J\},$$ (see
\cite[Appendix 40, Theorems 99, 102]{Matsumura1980}).
\end{definition}

Property D asserts that the differential operator $\Delta$ has the same nice properties when defined
over the class of schemes in ${\mathcal S}_0$. This is enough to guarantee that the algorithm of
resolution of basic objects, and hence the algorithms of principalization and resolution of
singularities, hold within the class of schemes in ${\mathcal S}_0$.  We make this idea more precise in
the following paragraph.

\begin{parrafo}
{\bf The general strategy.}  We have approached the problems of
desingularization and   strong principalization of ideals in a
unified way: by means of the notions of basic objects and
resolution of basic objects. Given a basic object $B=(W,(J,b),E)$,
where  $W$ is smooth over a field ${\mathbf{k}}$ of characteristic
zero, we have defined a function $\ord^{d}_B: \Sing(J,b)  \to
{\mathbb Q}$ (\ref{ordn}), basically in terms of the function
$$\ord: W \to {\mathbb Z},$$ given by the order of the ideal $J$
at
   points of $W$. In  \ref{PropOrdDelta} we show that if $W$ is smooth
   over a field $k$, this function is
   upper-semi-continuous. Furthermore, because  $W$ is smooth, we can
   define the   \(\Delta\) operator which acts on ideals of $W$.  Note
   that the order of $J$  at the local regular ring $\calo_{W,x}$ is,
   say, $d$, if and only if \( x \in V(\Delta^{d-1}(J))\setminus
   V(\Delta^{d}(J))\).

\

Resolution of basic objects (and hence desingularization and
strong principalization) follow by induction on the dimension of
the ambient space $W$. To be precise, and here is where a
constraint on the characteristic of the underlying field
${\mathbf{k}}$ is imposed, the order of $J$ at $\calo_{W,x}$ is
$d$ if, and only if, the order of $\Delta^{d-1}(J)$ at
$\calo_{W,x}$ is 1 (see {\bf (P2)} in \ref{PropOrdDelta} and
Example \ref{exchar}). An element of order one in
$\Delta^{d-1}(J)\subset \calo_{W,x}$ defines a smooth
hypersurface, say $H$ locally at $x \in W$ (so $I(H) \subset
\Delta^{d-1}(J)$). After  restriction, we may assume that
$\Sing(J,b) \subset H$. Then we show  that $H$ is a good candidate
for induction on the ambient space, by showing that the inclusion
$I(H) \subset \Delta^{d-1}(J)$) (and hence the inclusion
$\Sing(J,d) \subset H $) is stable by transformations of basic
objects. In doing so, we will only use the fact that
$\Delta^{d-1}(J)$ is defined in terms of partial derivations (see
Section \ref{SeccDelta}).
\end{parrafo}

\begin{lemma}
\label{propertyd} (PROPERTY D) If $W=\mbox{Spec} (A)$ is a scheme in the class ${\mathcal S}_0$,
$\mathfrak{p}\ A$  is a prime ideal  and $J$ is any ideal of $A$, then the order of $ J A_\mathfrak{p} $
in the local regular ring $A_\mathfrak{p}$ is greater than or equal to $b$ if and only if
$\Delta^{b-1}(J) \subset \mathfrak{p}$.
\end{lemma}

\begin{proof}
Let  $ R $ the localization of $A$ at a maximal ideal of $A$ containing $\mathfrak{p}$,  and let
$\hat{R}$  the completion of $R$. Then, the residue field of $\hat{R}$ is a finite extension $K$ of
${\mathbf{k}}$ and $\mbox{Der}_{\mathbf{k}}(R) $ induces $\mbox{Der}_{{\mathbf{k}}}(\hat{R}) $ over
$\hat{R}$ (see \cite[Theorem 99, (4)]{Matsumura1980}). By \cite[Theorem 102]{Matsumura1980}, $R$ is
excellent, so $\mbox{Spec}( \hat{R}) \to \mbox{Spec}(R)$ is a faithfully flat morphism with regular
fibers.

\

Note that  the order of $J$ defines an upper-semi-continuous function on $\mbox{Spec}(R)$. Since $J$ is
finitely generated it suffices to check this when $J$ is principal.  In that case, set $ \bar{R}=R/J$;
the order of $ J $ at a prime $\mathfrak{p}$ is the multiplicity of the ring $ \bar{R}$ at the prime
ideal, say $\bar{\mathfrak{p}}$, induced by $\mathfrak{p}$, and multiplicity is upper-semi-continuous.

\

Assume first that $\mathfrak{p}$ is a regular prime ideal in $R $. Then it expands to a regular prime
ideal in $\hat{R}$ and it is easy to see that the Property D holds for regular primes by checking that it
holds at the completion, which is a ring of formal power series.

\

An arbitrary prime ideal $\mathfrak{p}\subset R $ is the intersection of all prime ideals of height $
n-1$ , where $n $ is the dimension of $R$. In fact, if $ \mbox{dim} R/\mathfrak{p}= n-h $ and $ f \in R $
is not in $\mathfrak{p}$, one can find a prime ideal of height $ n-1 $ containing $\mathfrak{p}$ but not
$ f $ (set $ \bar{f}_1 \in R/\mathfrak{p}$ as the class of $ f $, extend to $ \bar{f}_1
,\bar{f}_2,\ldots\bar{f}_{n-h} \in R/\mathfrak{p}$ a system of parameters, and now take
$\mathfrak{q}\subset R $ by lifting a minimal prime ideal containing $ <\bar{f}_2 ,\ldots\bar{f}_{n-h}>
\subset R/\mathfrak{p}$).

\

Hence it suffices to check Property D for $\mathfrak{p}$ a prime ideal defining a curve. Now we use  a
trick by Hironaka (see \cite{Hir64}), which consists in reducing to the regular case: Let $$\pi: W' \to W
$$ be an  an embedded desingularization of the curve defined by $\mathfrak{p}$, obtained as the
composition of a finite sequence of quadratic transformations. Let $\mathfrak{p_1}$ be the defining ideal
of the desingularization of $\mathfrak{p}$, and set ${\mathcal O}_{W,x}$ and ${\mathcal O}_{W',y}$ as the
local rings of $\mathfrak{p}$ and  $\mathfrak{p_1}$ respectively.  Then the lemma follows immediately
from the fact that ${\mathcal O}_{W,x}$ and ${\mathcal O}_{W',y}$ are isomorphic, because in this case,
we can identify the localizations of $\mbox{Der}_{\mathbf{k}} (W_1)$ at $y$ with that of
$\mbox{Der}_{\mathbf{k}} (W)$ at $x$, and use the fact that the operator $ \Delta $, is defined in terms
of $\mbox{Der}_{\mathbf{k}} (W)$ and that $\Delta$ is a well defined operator on the class of  coherent
ideals over $ W $ (this can be checked at the completion at closed points).
\end{proof}



\section{Non-embedded desingularization.}
\label{nonembeddeddes}

If $X\subset W$ is a closed subscheme  then an embedded
desingularization of $X$, $$\begin{array}{ccc} W_r & \to & W \\
\cup & & \cup \\ X_r & & X,
\end{array}$$
defines a non-embedded desingularization, in the sense that it defines
        a proper birational morphism
$$X_{r}\longrightarrow X_{0}$$ such that $X_{r}$ is regular and the morphism is an isomorphism on \(
\Reg(X_{0}) \).

\

Our procedure of desingularization will also define a non-embedded desingularization of schemes which can
be locally embedded in smooth schemes.  This is not a  restriction at all if we consider Noetherian
separated schemes \( X\) of finite type over a field \( \mathbf{k} \). To prove that
Theorem~\ref{classical} extends to this class of schemes, we only have to prove that given  two different
embeddings of \( X \) we obtain the same non-embedded desingularization.

\begin{theorem}
\label{ThExtLocEmb}
        Let \( X\) be a Noetherian separated scheme of finite type over
        a field of characteristic zero, \( \mathbf{k} \). Consider two
        different closed immersions \( X\subset W\) and \( X\subset
        W'\), where \( W \) and \( W'\) are pure dimensional smooth
        schemes over \( \mathbf{k} \). Let
\begin{equation}
\label{nonem}
\begin{array}{ccccccc}
W & \leftarrow  & W_1 & \ldots & W_{r-1} &\leftarrow & W_r \\
\cup & & \cup & & \cup & & \cup \\ X &  & X_1 & \ldots & X_{r-1}
& & X_r,
\end{array}
\end{equation}
and
\begin{equation}
\label{nonem1}
\begin{array}{ccccccc}
W^{\prime} & \leftarrow  & W_1^{\prime} & \ldots & W_{r-1}^{\prime} &\leftarrow & W_r^{\prime} \\
\cup & & \cup & & \cup & & \cup \\ X &  & X_1^{\prime} & \ldots &
X_{r-1}^{\prime} & & X_r^{\prime},
\end{array}
\end{equation}
be the embedded desingularizations provided by Theorem \ref{classical}.  Then sequences (\ref{nonem}) and
(\ref{nonem1}) define the same non-embedded desingularization  of $X$, i.e., if $$\varphi:
X_{r}\longrightarrow X, \qquad \mbox{ and } \qquad
        \varphi':X'_{r'}\longrightarrow X_{0},$$ are the two induced
        non-embedded desingularizations, then
$$X_{r}=X'_{r'} \qquad \mbox{ and } \qquad \varphi=\varphi'.$$ Moreover the number of blowing ups also
coincide, i.e. \( r=r' \).
\end{theorem}

We refer  to \cite{EncVil99} for the proof of Theorem
\ref{ThExtLocEmb} (see also Remark \ref{renoinm}), but let us
mention here some auxiliary results which contribute to the proof
of this theorem. These results are interesting by themselves,
because  they describe the behavior of the algorithm of resolution
of singularities provided by Theorem \ref{classical} under regular
extensions.

\begin{proposition}
\label{EtaleExtEmbRes} Let \( W\) be a pure dimensional smooth scheme over a field  \(
        \mathbf{k} \) and let \( X\subset W\) be a closed subscheme.
        Then:
\begin{enumerate}
\item[(i)] If \( W'\longrightarrow W\) is  an \'etale morphism and \(
      X'\subset W' \) is the pull-back of \( X \), then the
      desingularization  of $X'\subset W'$ provided by Theorem
      \ref{classical} in \ref{existence}, say
$$(X'_{r}\subset W'_{r})\longrightarrow (X'\subset W')$$ is the
      pull-back of the desingularization of $X\subset W$  provided by
      Theorem \ref{classical}, say
$$(X_{r}\subset W_{r})\longrightarrow (X\subset W).$$
\item[(ii)] If \( W'\longrightarrow W \) is defined by an arbitrary
      extension of the base field, then the desingularization of the pull-back of $X$ in $W'$,
      $X'\subset W'$ provided by Theorem \ref{classical},
$$(X'_{r}\subset W'_{r})\longrightarrow (X'\subset W')$$ is the
      pull-back of the desingularization of $X\subset W$ provided by
      Theorem \ref{classical},
$$(X_{r}\subset W_{r})\longrightarrow (X\subset W).$$
\item[(iii)] More generally, if \( W'\longrightarrow W \) is a regular map, and both $W$ and
$W'$ are schemes in ${\mathcal S}_0$, then the desingularization of
   the pull-back of $X$ in $W'$,     $X'\subset W'$ provided by Theorem \ref{classical},
$$(X'_{r}\subset W'_{r})\longrightarrow (X'\subset W')$$ is the
      pull-back of the desingularization of $X\subset W$ provided by
      Theorem \ref{classical},
$$(X_{r}\subset W_{r})\longrightarrow (X\subset W).$$
        \end{enumerate}
\end{proposition}

Since the embedded desingularization stated in Theorem \ref{classical}
    is defined in terms of the resolution of a suitable defined basic
    object (see \ref{existence}), Proposition \ref{EtaleExtEmbRes} follows
    immediately from the next lemma:

\begin{lemma}
\label{EtaleExtAlgRes} Let \( (W,(J,b),E) \) be a basic object.
Then the algorithm of resolution of basic objects defined in
\ref{AlgResol} and constructed in \ref{Theoremd}, satisfies the
following  additional property: If \( W'\longrightarrow W \) is a
regular map,
        both $W$ and $W'$ are in ${\mathcal S}_0$,  and \( J'\subset\O_{W'} \) and \( E'\)
        are respectively the pull-backs  of \( J\) and \( E\) in $W'$,
        then the resolution of \( (W',(J',b),E') \) defined
in \ref{AlgResol} is also the pull-back of the resolution of
        \( (W,(J,b),E) \) defined in  \ref{AlgResol}.
\end{lemma}

\begin{proof} It is enough to observe that  the proof of
Theorem~\ref{Theoremd} relies on the functions \( \ord \) and $n$ introduced
 in Definition \ref{ordn},
and that these functions are naturally compatible with  regular
maps between schemes in ${\mathcal S}_0$:

\

Let $x'\in W'$, let $x=\pi(x')$, let $J\subset {\mathcal O}_W$ be
a non-zero sheaf of ideals, and let  $y\in W$ be a closed point
contained in  $\overline{\{x\}}$. Consider the diagram:
$$\begin{array}{ccc} {\mathcal O}_{W,x} & \to & {\mathcal
O}_{W,x'}\\ \uparrow & & \\ {\mathcal O}_{W,y} & &
\end{array}$$

If $\ord_xJ=b$, then $\Delta_y^{b-1}(J)\subset {\mathcal O}_{W,y}$
 is a proper ideal. To check that $$\ord_xJ=\ord_{x'}J{\mathcal O}_{W'},$$ note that
the diagram $$\begin{array}{rcl} {\mathcal O}_{W,x} &
\hookrightarrow & {\mathcal O}_{W',x'} \\ \downarrow & &
\downarrow \\ k_1[|x_1,\ldots,x_n|]=\hat{\mathcal O}_{W,x} &
\hookrightarrow & \hat{\mathcal
O}_{W',x'}=k_2[|x_1,\ldots,x_n,x_{n+1},\ldots,x_m|]\end{array}$$
commutes, where $k_1\subset k_2$ since the characteristic is zero,
and that all the maps are faithfully flat.

\end{proof}

An important result in the study of equivariance for
non-embedded schemes is the following lemma:

\begin{lemma} \label{LemaEtale}
         Let \( W\) and \( W'\) be pure dimensional schemes, smooth
         over \( \mathbf{k} \), with \( \dim{W}=\dim{W'} \), and let \(
         J\subset\O_{W} \) and \( J'_{0}\subset\O_{W'} \), be two
         sheaves of ideals.  Assume that for two points \( \xi\in W \)
         and \( \xi'\in W' \) there is an isomorphism \(
         \Theta:\hat{\O}_{W,\xi}\longrightarrow \hat{\O}_{W',\xi'} \)
         such that \( \Theta(\hat{J})=\hat{J}'\) (where \(
         \hat{\O}_{W,\xi} \) and \( \hat{\O}_{W',\xi'} \) denote   the
         completions of $\calo_{W,\xi}$ and $ \calo_{W',\xi'}$
         respectively, and \( \hat{J}_{0}=J_{0}\hat{\O}_{W_{0},\xi_{0}}
         \), \( \hat{J}'_{0}=J'_{0}\hat{\O}_{W'_{0},\xi'_{0}} \)). Then
         there is a common \'etale neighborhood \(
         \tilde{\xi}_{0}\in\widetilde{W}_{0} \) of both \( \xi_{0} \)
         and \( \xi'_{0} \), and and an ideal \(
         \widetilde{J}_{0}\subset\O_{\widetilde{W}_{0}} \) such that
         \begin{equation*}
        \widetilde{J}_{0}=J_{0}\O_{\widetilde{W}_{0}}=
             J'_{0}\O_{\widetilde{W}_{0}}.
         \end{equation*}
\end{lemma}
\begin{proof}
If $ \Theta $ arises from an isomorphism $\Theta: W_0 \to W'_0$ mapping \( \xi_{0}\in W_{0} \) to
\(\xi'_{0}\in W'_{0} \), and $\Theta (J_0)= J'_0$, then we it is enough to take $\widetilde{W}_{0}=W_0$
(note that an isomorphism is an \'etale map). We claim that this is the case in general, at least
replacing $W$ and $W'$ by suitable \'etale neighborhoods in the given points: Since the local rings
$\calo_{W,\xi}/J_{\xi}$ and $\calo_{W',\xi'}/J'_{\xi'}$ are formally isomorphic, then their
henselizations are isomorphic (see \cite[2.6]{Artin}). Now, an isomorphism of these henselizations can be
lifted to an isomorphism of the henselizations of the regular local rings, say $\Gamma:
\overline{\calo}_{W,\xi} \to \overline{\calo}_{W',\xi'} $, mapping $J  \overline{\calo}_{W,\xi}$ to $J'
\overline{\calo}_{W',\xi'}$. Since both henselizations are direct limits of \'etale neighborhoods,
$\Gamma$ also defines an isomorphism, as indicated above, at suitable \'etale neighborhoods.
\end{proof}




\section{Equiresolution. Families of schemes.}
\label{families} Once we fix a constructive algorithm of
desingularization (i.e.  once we attach to each variety a
particular desingularization), it makes sense to ask whether one
can classify varieties in accordance to their (algorithmic)
desingularization. The  formulation of this question already
requires some clarification; but if we are to think of an embedded
variety as a structure defined in terms of equations, it is
conceivable to ask how the algorithm of desingularization will
behave when the coefficients of these equations vary.
\medskip

The questions about classification of varieties lead, quite
naturally, to the notion of {\em families of schemes}.  In this
section we present briefly some results proved in
\cite{EncinasNobileVillamayor} about the good behavior of the
algorithm of desingularization in families and the notion of
equiresolution.

\

Given a smooth morphism  \( W\rightarrow T \),  we will denote by $ W^{(t)}$, the fiber at the point  \(
t\in T \), i.e., $$ W^{(t)}=W\times_{T}\Spec(k(t))$$  where \( k(t) \) is the residue field of \(T\) at
\(t\).  Note that \( W^{(t)} \) is smooth over the field \( k(t) \).

\begin{definition}
\label{familyschemes} Given a smooth morphism  \( W\rightarrow T \),  a {\em family of embedded schemes}
is   a closed subscheme $X \subset W$ such that induced morphism $ X \to T$ is  flat.
\end{definition}

\begin{example}\label{ejemplo1}
Let $C \subset {\mathbb A}^2_k$ be a plane curve smooth over the field $k$ of characteristic zero.  For
an arbitrary scheme $T$ (say for instance a non reduced curve), $C\times T \subset {\mathbb A}^2_k \times
T$ defines a flat family over $T$.
\end{example}

\begin{example}\label{ejemplo2}
Let $ \mathbb{A}^2_k \to {\mathbb A}^1_k$ be the projection on first the coordinate. This is a smooth,
and hence flat morphism. If we now take $W= {\mathbb A}^2_k =\mbox{Spec}(k[x_1,x_2])$, $T= {\mathbb
A}^1_k $, and $X=V(x_1.x_2)\subset W$, then this is not a family since $ X \to T$ is not flat.
\end{example}

\begin{example}\label{ejemplo3}
Let $ {\mathbb A}^4_k \to {\mathbb A}^1_k$ be the projection on the first coordinate (i.e. defined by
$k[x_1]\to k[x_1,x_2,x_3,x_4]$ the natural inclusion).  Set $W= {\mathbb A}^4_k
=\mbox{Spec}(k[x_1,x_2,x_3,x_4])$,  $T= {\mathbb A}^1_k =\mbox{Spec}(k[x_1])$,  and $X=V(x_2^2+x_1x_3^3,
x_4) \subset W$. This example defines a family; the members of this family consist of irreducible
singular curves, except at the the origin of ${\mathbb A}^1_k$ where the fiber is a double line.
\end{example}

So given  a smooth morphism \( W\rightarrow T \), we view $X\subset W$ as a structure defined, at least
locally, by equations with {\em coefficients} in the scheme $T$. The requirement of flatness on $X \to T$
somehow avoids abrupt changes among the different fibers, a fact illustrated in Example \ref{ejemplo2}.

\

There are various contexts where the notion of  families of schemes
  appear, like for instance when dealing  with Hilbert Schemes
  (see for instance \cite{Stromme}):  Fix a projective space ${\mathbb
P}^n_k$ and a
  Hilbert polynomial  $F(t)$; then there is a scheme $S$, together
  with a setting  \( X \subset W={\mathbb P}^n_k \times S \rightarrow S
  \) defining an embedded family such that the fibers parameterize {\em
  all possible} subschemes $Y\subset\mathbb{P}^n_k$ with Hilbert
  polynomial $F(t)$.

\

If we are to accept the notion of families as the right setting where to study the behavior of
desingularization, a new problem arises, as illustrated in Example \ref{ejemplo3}:  It makes sense to
formulate resolution of singularities for the general fibers, consisting of irreducible singular curves,
but not for the particular fiber at the origin, consisting of a double line. Since we prove here that
embedded desingularization is a byproduct of embedded principalization (see \ref{existence}), we overcome
this problem by simply replacing the use of embedded desingularization by principalization of the defining
ideal. This   will be explained with more detail in the  next point of this section (see
``Equisolvability I. Families parameterized by a smooth scheme $T$'' below). The results in
desingularization will  hold for a suitable class of families of schemes: The class of  {\em equisolvable
families}. Thus,  our objectives are:

\begin{itemize}
\item[1.]  To introduce the notion of  {\em equisolvable families}.
The definition of equisolvability  should be stable by fiber products, namely: If  \( X \subset W
\rightarrow T \), is {\em equisolvable}, and if \( X_1 \subset W_1 \rightarrow T_1 \), is defined by
taking fiber products with $ T_1 \to T$, then the later should be equisolvable.

\item[2.] Given  an arbitrary family \( X \subset W \rightarrow T \),
we will want to find a partition of algebraic nature on the underlying topological space of $T$, (a
partition into locally closed sets $T_{\alpha}$)  such that each restricted family, say \( X_{\alpha}
\subset W_{\alpha} \rightarrow T_{\alpha} \), is equisolvable.
\end{itemize}

Note that 2 would provide a formal setting to state that constructive resolution of singularities is {\em
algebraic on the coefficients}, in fact we can take (2)  as a definition, namely, that for {\em any}
family \( X \subset W \rightarrow T \), there is a partition, of algebraic nature, on the underlying
topological space of $T$, such that each restricted family is equisolvable.

\begin{center}
\textbf{Equisolvability I. Families parameterized by a smooth scheme $T$.}
\end{center}

In this part we are going to restrict our attention to  the case of
  families \( X \subset W \rightarrow T \), where $T$ is smooth over a
  field of characteristic zero. Let
  $\mathcal{I}=\mathcal{I}(X)\subset\calo_{W}$ be the defining ideal
  sheaf of $X$. We want to introduce conditions that insure that all
  the different members of the {\em family of ideals}
\begin{equation}\label{famofideal}
{\mathcal G}= (\pi:W \to T, {\mathcal I})
\end{equation}
can be simultaneously principalized by using the algorithm. When this happens, we shall say that the
family is {\it equisolvable}.

\

We propose two conditions (\ref{DefCondAE} and \ref{DefCondTau}), which turn out to be equivalent (see
Theorem \ref{tau}). Condition~\ref{DefCondAE} does not explicitly involve the fibers of the family, but
it requires, essentially, that the centers $C_i$ that appear in the principalization sequence that the
algorithm associates to ${\mathcal I} \subset \calo_W $, $(W,E)$ be ``evenly spread'' over the parameter
space $T$. In condition~\ref{DefCondTau}, a ``numerical'' invariant is associated to the different points
$t \in T$ (this invariant in defined in terms of the principalization sequences of the fibers); it is
required that it be constant along $T$. Both approaches have their advantages, depending on the
situation. One of our aims is to show that given an arbitrary family of ideals, it is possible to
naturally stratify the parameter space $T$ as  a union of locally closed sets so that, along each one, the
restriction of the family is equisolvable.

\

Finally recall that the given embedded  principalization algorithm induces an {\it associated
desingularization algorithm} for couples $X \subset W$ (see \ref{existence}).  We shall say that a family
of embedded schemes \( X \subset W \rightarrow T \) is {\it equisolvable} (relative to the algorithm), if
the associated family of ideals (\ref{famofideal}) is equisolvable. If the family of embedded schemes is
such that all the fibers $X_t$ are reduced, then the desingularization sequence that the associated
algorithm assigns to  $X \subset W$ induces, on each fiber, the resolution sequence that corresponds to
that fiber; that is, it has the property to be expected of a good notion of simultaneous resolution.

\

But this definition of equisolvability applies also to the case where some (or all the) fibers $X^{(t)}$
are non-reduced. So, we have a (we hope, reasonable) notion of equiresolution for families of embedded
schemes where some fibers may be non-reduced. Since in many geometric problems it is unavoidable the
presence of non-reduced fibers in families of schemes, it is an important feature.

\begin{parrafo}
\label{DefCondAE} {\bf Condition of algorithmic equiresolution (AE).}   With the same notation as above,
consider the basic object
\begin{equation} \label{EqBasOFam}
      (W_0,(\mathcal{I},1),E_0=\emptyset),
\end{equation}
and let
\begin{equation}\label{prinseq}
(W_0,E_0) \longleftarrow (W_1,E_1) \longleftarrow ..\longleftarrow (W_r,E_r)
\end{equation}
be the embedded principalization sequence defined by the
algorithm, where    each $W_i \longleftarrow W_{i+1}$ is  defined
by blowing up at $C_i=\Max(f_i)$, which is determined by  the
functions defining the algorithm (see Theorem
\ref{principalization}, Definition \ref{AlgResol} and
\ref{proofprincipal}). Note that we get, by composition, smooth
morphisms $$\pi_{i}: W_i \to T,$$ which induce morphisms $\rho_i:
C_i \to T$, $i=0,\ldots,r-1$.  We will say that the family of
ideals ${\mathcal G}$ satisfies  the {\em condition of algorithmic
equiresolution}  if the morphism
\begin{equation}
\label{dondeesta} \rho _i: C_i \to T
\end{equation}
  is smooth, proper and surjective, for $i=0, \ldots, r$.
\end{parrafo}

\begin{parrafo}
\label{DefCondTau} {\bf Condition $\tau$.}  Given a family of
ideals ${\mathcal G}$ we are going to define a function $\tau
_{{\mathcal G}}$ from the parameter space $T$ into a certain
totally ordered set ${\Lambda ^{(d)}}$ (which depends on $d = $
dimension of the fibers of $\pi: W \to T $ only):
\begin{enumerate}
\item[(i)] The values of $\tau_{\mathcal{G}}$ are sequences whose entries are
either in $\mathbb{Z}$, or in the set $I^{(d)}$  introduced in
Definition \ref{AlgResol},  or (possibly) $\infty$. This sequences
will be ordered lexicographically. Therefore if we take $I_d\times
{\mathbb Z}$ ordered lexicographically and set $L=I_d\times
{\mathbb Z}\bigsqcup \{\infty\}$, we have that
$${\Lambda ^{(d)}}=L\times L\times \ldots.$$

\item[(ii)] For each point $t \in T$, let $
(W^{(t)}_0,(\mathcal{I}^{(t)},1),E^{(t)}_0)$ be the fiber of
$(W_0,({\mathcal I}, 1),E_0=\emptyset)$ at $t\in T$, and now
consider the sequence of transformations of pairs
\begin{equation}\label{tprinseq}
(W^{(t)}_0,E^{(t)}_0) \longleftarrow (W^{(t)}_1,E^{(t)}_1) \longleftarrow\cdots\longleftarrow
(W^{(t)}_{r_{t}},E^{(t)}_{r_{t}})
\end{equation}
which provides an embedded  principalization of the fiber of the
ideal ${\mathcal I}$ over $t$. Let $c^{(t)}_i$ be the number of
connected components of $C^{(\bar{t})}_i := \Max f^{(\bar t)}_i$.
Here  $C_{i}^{(\bar{t})}=C_{i}\times_{T}\overline{k(t)}$ is the
geometric fiber ($\overline{k(t)}$ the algebraic closure of
$k(t)$).  Then we define
\begin{equation*}
      \tau_{{\mathcal G}}(t) = ( \max{f^{(t)}_0},c^{(t)}_0,
      \max{f^{(t)}_1},c^{(t)}_1, \ldots,
      \max{f^{(t)}_{r_t}},c^{(t)}_{r_t}, \infty,\infty, \ldots)\in {\Lambda
^{(d)}}.
\end{equation*}
\end{enumerate}
We will say that the family ${\mathcal G}$ satisfies the {\em
condition $\tau$} if $\tau _{{\mathcal G}} (t)$ is constant for
all $t \in T$.
\end{parrafo}

\begin{theorem}\label{tau} \cite[Theorem~2.3]{EncinasNobileVillamayor}
Let $T$ be an integral and smooth scheme, and  let  ${\mathcal G}=(\pi:W \to T, \mathcal{I})$ be   a
family of ideals. Assume that for all $i=0,\ldots,r-1$, the morphism $\rho _i:C_i \to T$ described in
(\ref{dondeesta})  is proper.  Then ${\mathcal G}$ satisfies condition (AE) if and only if it satisfies
condition~$\tau$,  and in either case, the principalization sequence (\ref{prinseq}) of $ ( W , {\mathcal
I}, E) $ induces, by taking fibers, the principalization sequence of ${\mathcal G}(t):=(W^{(t)},{\mathcal
I}^{(t)})$, for all $t \in T$.
\end{theorem}

The last part of the statement means the following. The length $r$ of the principalization sequence of $
( W , {\mathcal I}, E) $ agrees with the length $r_t$ of the principalization sequence of any fiber
${\mathcal G}(t):=(W^{(t)}, {\mathcal I}^{(t)})$; there  is a natural identification of $\pi_i ^{-1}(t)$
and $W_i ^{(t)}$ (\ref{tprinseq}), for any $t \in T$, and via this identification, the restriction of
$f_i$ to $(\pi _i ^{-1}(t))$ coincides with $f^{(t)}_i$, for $i=0, \ldots, r$.

\begin{definition}
\cite[definition~1.10]{EncinasNobileVillamayor} \label{DefEquiIdeal} A family of ideals ${\mathcal G}=
(\pi:W \to T, {\mathcal I})$ is {\it
  equisolvable} if any of the equivalent conditions of Theorem
  \ref{tau} hold.
\end{definition}

Note that the first hypothesis in Theorem \ref{tau} are
automatically fulfilled if the morphism $\pi: W \to T$ is proper,
in particular if it is projective.

\begin{remark}
Let ${\mathcal G}= (\pi:W \to T, {\mathcal I}, E)$ be a family of ideals. Let $T^{*} \to T$ a smooth
morphism of smooth schemes, and set ${\mathcal G}^{*}= (\pi^{*}:W^{*} \to T^{*}, {\mathcal I}^{*}, E^{*})$
by taking fiber products. If ${\mathcal G}$ fulfills the conditions of the Theorem, then so does
${\mathcal G}^{*}$.
\end{remark}

We mention now some results about families of embedded schemes.

\begin{definition}
We say that a family of embedded schemes ${\mathcal F}=(X\subset W\to T)$ is {\it equisolvable} if the
associated family of ideals is equisolvable.
\end{definition}

\begin{proposition}
Assume that ${\mathcal F}=(X,W,\pi)$ is an equisolvable family of embedded schemes, where for all $t \in
T$ the fiber $X^{(t)}$ is reduced . Then, for all $t\in T$, the resolution sequence (\ref{prinseq})
induces the resolution sequence of $(X^{(t)},W^{(t)})$.
\end{proposition}

\begin{proof}
By definition, the equisolvability of ${\mathcal F}$ means that the associated family of ideals
${\mathcal G}=(\pi:W \to T, \mathcal{I}(X), \emptyset)$ is equisolvable. Letting $W_0=W$,
$\mathcal{I}=\mathcal{I}(X),~E_0=\emptyset$, consider the principalization corresponding to  ${\mathcal
T}=(W_0, {\mathcal I}, \emptyset)$.  Equisolvability of ${\mathcal G}$ asserts that the principalization
sequence of ${\mathcal G}$, induces that of the fiber ${\mathcal G}(t):=(W_0^{(t)}, {\mathcal I}_0^{(t)},
E_0^{(t)})$, in the following sense. We have $C_0^{(t)}=C_0 \cap W_0^{(t)},~f_0^{(t)}=f_0|W_0^{(t)}$ and
the strict transform of $W_0^{(t)}$ in $W_1$ can be identified with $W_1^{(t)}$. Via this isomorphism,
$f_1^{(t)}=f_1|W_1^{(t)}$, $C_1^{(t)}=C_1 \cap W_1^{(t)}$, and so on. Eventually, after $r$ steps, both
the principalization sequence of ${\mathcal T}$ and ${\mathcal G}(t)$ simultaneously stop.

Now $I(X){\mathcal O}_{W_0^{(t)}}=I(X_0^{(t)})$,  and  to obtain a desingularization sequence for
$(X_0^{(t)},W_0^{(t)})$ we use the principalization sequence of the fiber ${\mathcal
G}(t)=(W_0^{(t)},{\mathcal I}_0^{(t)}, \emptyset)$. It is clear that, via the identification of
$W_i^{(t)}$ with a suitable subscheme of $W_i$, for all possible index $i$, $X_i^{(t)}$ corresponds to the
strict transform of $X_0^{(t)}$ to $W_i$. Moreover, the length $s_t$ of the resolution sequence of
$(X_0^{(t)},W_0^{(t)})$ is constant, equal to the length $s$ of the resolution sequence of $(X,W)$. In
fact
  $s$ (resp. $s_t$) is the unique index such that the
proper transform $X_{s} \subset W_{s}$ (resp. $X_{s}^{(t)} \subset W_{s}^{(t)})$ has the same codimension
as the center $C_{s}$ (resp. $C_{s}^{(t)}$) (see \ref{existence}). But $\codim
(X_i^{(t)},W_i^{(t)})=\codim(X_i,W_i)$ and, since $\mathcal G$ is equisolvable,
$\codim(C_i^{(t)},W_i^{(t)})=\codim(C_i,W_i)$.  This proves our contention about the indices. Since the
desingularization functions are defined by restriction, the proposition is proved.
\end{proof}
\medskip

\begin{center}
\textbf{Equisolvability II. Stratification of Hilbert schemes.}
\end{center}

\begin{definition} \label{DefGenFam}
A {\it  general family of embedded schemes} is a pair
\begin{equation*} \label{EqGenFamSch}
      {\mathcal F} = (j: X \to W, \pi:W\to T)
\end{equation*}
where $ \pi:W \to T $ is a smooth morphism of equidimensional Noetherian schemes over a field
${\mathbf{k}}$ of characteristic zero, $j$ is a closed immersion and $p:= \pi\circ j:X \to T$, is flat
(but we do not assume that $W $ or $ T $ are regular).
\end{definition}

\begin{definition}
A {\it  general family of ideals} is a pair
\begin{equation*} \label{EqGenFamIde}
      {\mathcal G}=(\pi:W \to T, {\mathcal I})
\end{equation*}
where $ \pi:W \to T $ is a smooth morphism  (but we do not require $ W $ or $ T $ to be regular or
irreducible).  All fibers ${\mathcal G}(t):=(W^{(t)},{\mathcal I}^{(t)})$ give ideals
$\mathcal{I}^{(t)}\subset\O_{W^{(t)}}$.
\end{definition}

\begin{remark}
Note that the function $\tau _{{\mathcal G}}: T \to \Lambda^{(d)}$
introduced in \ref{DefCondTau}  still can be defined for any
general family. If $ T $ admits a desingularization $ T_1 \to T $
then $ {\mathcal G}_1 $ (obtained from ${\mathcal G}$ by base
change to $T_1$) is a family in the sense of
Definition~\ref{famofideal}. Hence Theorem~\ref{tau} applies to
${\mathcal G}_1$.  Since  $\tau _{{\mathcal G}}$ is defined in
terms of the fibers only, from the properness of $ T_1 \to T $ we
see that the conclusion of Theorem \ref{tau} is also valid for the
family ${\mathcal G}$.
\end{remark}

Given a general family of embedded  schemes, ${\mathcal F} = (j: X \to W, \pi:W \to T)$,  we can
associate to it a  general family of ideals (as we did in (\ref{famofideal})),  say ${\mathcal G}=(\pi:W
\to T, {\mathcal I})$ (\ref{EqGenFamIde}), where ${\mathcal I}=I(X)$.  For simplicity, in the sequel we
denote by $\tau _{{\mathcal F}}$ the function $\tau _{{\mathcal G}}$ corresponding to the family of ideals
${\mathcal G}$.

\begin{parrafo} \label{SetHilbFam}
Now fix the following objects:
\begin{itemize}
\item[(i)] A graded algebra $S$ over a field $k$ of characteristic
zero, finitely generated by elements of degree one, such that
$W={\rm Proj}(S)$ is smooth of finite type over $k$, of dimension
$n$ (which  belongs to the class ${\mathcal S}_{0}$).
\item[(ii)] An element $\alpha\in {\Lambda
^{(d)}}$.
\item[(iii)] A polynomial $Q$ with rational coefficients.
\end{itemize}

Let ${\mathcal H}(W, \alpha, Q)$ be the class of all general families of embedded and reduced schemes
${\mathcal F}_{T}$, of the form $( j : X_T \to W_T, \pi:W_{T}\to T)$ where $\pi$ is obtained from $ W \to
\Spec(k) $ by base change, and $ X_T \subset W_T $ is a closed subscheme such that the induced projection
$ X_T \to T $ is flat ($j $ being the inclusion). We require also for
$\mathcal{F}_{T}\in\mathcal{H}(W,\alpha,Q)$ that:
\begin{itemize}
\item[(0)] $T$ is a reduced scheme of finite type over $k$.

\item[(1)] For all $t\in T$, $\tau _{{\mathcal F},T}(t)= \alpha $
(hence the pull-back of this general family via $T' \to T$, with $T'$ regular is {\it equisolvable} since
it satisfies condition $ \tau $).

\item[(2)] If $( X^{(t)},W)$ is the couple induced by $(X_{T}, W_{T})$
over $t \in T$, then $Q$ is the Hilbert polynomial (relative to the line bundle corresponding to $S(-1)$)
of the embedded scheme $ X^{(t)} \subset W $.
\end{itemize}
\end{parrafo}

Now we can state the following theorem:
\begin{theorem} \label{ThHilb} \cite[4.10]{EncinasNobileVillamayor}
Under the conditions stated in  \ref{SetHilbFam} (and letting $
H:={\mathcal H}(W, \alpha,Q)$) there is a universal object in the
class of general families in $H$ (which will  be called a
universal $(\alpha,Q)$-equisolvable family). That is, there is a
general family, ${\mathcal F}_{H(\alpha,Q)}$ defined by $ (
X_{H(\alpha, Q)}\subset W_{H(\alpha, Q)}) $, $W_{H(\alpha, Q)} \to
H(\alpha, Q)$) such that for each general family ${\mathcal
F}_{T}$ in $ H $, there is a unique morphism $ T \to H(\alpha, Q)$
so that ${\mathcal F}_{T}$ is the pull-back of ${\mathcal
F}_{H(\alpha ,Q)}$.
\end{theorem}

\begin{proof}
Consider the Hilbert scheme $ H(Q) $, parameterizing subschemes of the projective variety $ W $ having
Hilbert polynomial $ Q $. We refer here to \cite[p. 21 (c)]{Mumford}  for a summary of results on Hilbert
schemes.  Let $ X(Q) \subset W \times H(Q)$ be the universal family; note that $ (X(Q), W \times H(Q))$
together with the projection $W \times H(Q) \to H(Q)$ defines a general family, say ${\mathcal F}(W,Q)$.
Then  \cite[Theorem~4.8]{EncinasNobileVillamayor} says that $\tau_{{\mathcal F}(W,Q)}: H(Q) \to
\Lambda^{(n)}$ is an LC-function. Hence its fibers define a partition of $ H(Q) $ into a disjoint union
of locally closed subsets. If $\alpha \in \Lambda^{(n)}$, let $ H(\alpha, Q):= [\tau _{{\mathcal
F}(W,Q)}]^{-1}(\alpha)$. Given a general family, say ${\mathcal F}_T $,  in $ H $, since the natural
morphism $ T \to H(Q) $ obtained by universality of Hilbert schemes has constant value $\alpha$, it
becomes clear that it factors though $ H(\alpha, Q) $ and vice versa, proving the theorem.
\end{proof}

\begin{definition}
An {\it LC-partition} of a scheme $T$, is an expression of $T$ as a disjoint union of locally closed
subsets (or subschemes, with the reduced structure).
\end{definition}

\begin{remark}
In Theorem~\ref{ThHilb} we have introduced an LC-partition in the Hilbert scheme $H(Q)$  which is
naturally related to equiresolution, by means of the (reduced) subschemes $ H(\alpha, Q) $. We might call
this a  {\it stratification} of $H(Q)$ but, in the literature, this term is often used in a more
restricted sense. For instance, in \cite{LipmanTirol} the notion of {\it stratification} of $T$ is
introduced  as an LC-partition (each set thereof is called a stratum) with the following properties:
\begin{itemize}
\item[(i)] The boundary of each stratum is a union of strata.

\item[(ii)] The singular locus of the closure of each stratum is a
union of stratum.

\item[(iii)] Each stratum is smooth.
\end{itemize}

In  \cite[2.5 (c)]{LipmanTirol}  it is remarked that given an LC-partition of an algebraic variety T, one
can attach to this partition a {\it coarsest} stratification of $T$, with the property each locally
closed subset of the partition is a union of strata. Here {\it coarsest} means that any other
stratification with this property is a refinement of the first.  So, if we want to obtain a
stratification  of the Hilbert scheme $H(Q)$, in the sense of \cite{LipmanTirol}, and  naturally related
to equiresolution, all what we need it to take the coarsest stratification associated to the LC-partition
described in Theorem~\ref{ThHilb}.
\end{remark}




\section{Bodn\'ar-Schicho's computer implementation.}
\label{Computer}

In \cite[Theorem 7.3.]{Villa92} and \cite[Theorem
7.13.]{EncVil97:Tirol} give  an algorithm of resolution of basic
objects. Based on this construction,  G.~Bodn\'ar and J.~Schicho
  wrote  a computer program which resolves
basic objects: Given a basic object, \( (W,(J,b),E) \)  the program provides the resolution defined in
terms of functions \( f_{i}^{d} \), as indicated in Definition \ref{AlgResol}, (c)
  (see
\cite{GabSch98} and \cite{GabSch99}). Their program is available at
\begin{center}
       \texttt{http://www.risc.uni-linz.ac.at/projects/basic/adjoints/blowup}.
\end{center}

As a consequence, the authors provided a computer program which
computes the resolution of singularities of a given hypersurface
(cf.  \cite{GabSch98} and \cite{GabSch99}). The main obstacle to
generalize this process to varieties of arbitrary codimension was
that the classical approach for desingularization was based on
reducing to the hypersurface case, and in this process, algorithms
of division and strict transforms of ideals are involved.

\

In the new approach of desingularization presented in this paper
all this is avoided. The  short proofs of embedded resolution of
singularities and embedded principalization of ideals presented in
\ref{existence} and \ref{proofprincipal} are only based on the
existence of an algorithm of resolution of basic objects. As a
consequence, based on these ideas,  G.~Bodn\'ar and J.~Schicho
have been able to adapt their program to produce one that resolves
the singularities of varieties of arbitrary codimension.
  Also, the program
can compute embedded principalization of ideals. We refer to Section \ref{nuevanueva} for some hints on
how their computer program works.

\

About the performance of the program, the complexity of computations seems to be very high but, by now,
we do not know any bound for complexity. In practice ``easy'' examples for \( \dim{W}\leq 4 \) may be
computed. For instance if the input is an ideal \( J\subset\mathbb{Q}[X_{1},\ldots,X_{n}] \) then the
resolution of the variety defined by \( J \) can be computed if \( n\leq 4 \) and the degree of the
generators of \( J \) is not bigger than 7 or 8. The first versions of the program were implemented in
\emph{Maple}, but the most recent version use \emph{Singular}, which performs much better all standard
basis computations.  For example, the resolution of the Whitney umbrella takes \( 1.230 \) seconds of CPU
time in a PowerPC G3 at 700 Mhz.  But sometimes one could find an example which seems ``easy'' but it
requires hard middle computations, since coefficients ideals are complicated. The new version
\cite{BodEncSch} has the option to use the feature of \emph{Singular} for computing in some localizations
rings which improves the performance.




\begin{center}

\part{Preliminaries for constructive resolutions.}
\label{nature}

\end{center}

We will prove resolution of basic objects by induction on the dimension of the ambient space:  To find a
resolution of a \( d \)-dimensional basic object \((W,(J,b),E)\) we will associate to it, at least
locally at each point of \(W\), a \((d-1)\)-dimensional basic object $(Z,({\mathcal A},e),G)$ in such a
way that a resolution of $(Z,({\mathcal A},e),G),$ is equivalent to a  resolution of \((W,(J,b),E)\) in
some sense that we will make precise in the forthcoming sections. This will force us to introduce some
machinery:
\begin{enumerate}
\item The notions of equivalence, inclusion and intersection of basic
objects (see Section \ref{describe}).

\item The notion and properties of the $\Delta$ operator, which will
provide us the way to find $(Z,({\mathcal A},e),G)$ given \((W,(J,b),E)\) (see Definition \ref{DefCoeff}).

\item The notion of {\em general basic object}: Given a  \( d
\)-dimensional basic object \((W,(J,b),E)\), in some situations we will associate to it {\bf locally}, at
each point of $W$, a \((d-1)\)-dimensional basic object. In other words, we will associate  to
\((W,(J,b),E)\), a collection of \((d-1)\)-dimensional basic objects, \((\widetilde{W}_{\alpha},
(I_{\alpha},a_{\alpha}), E_{\alpha})\). This  family of \((d-1)\)-dimensional basic objects will be what
we call a {\em general basic object} (see Section \ref{sectiongbo}).
\end{enumerate}



\section{Intersections,
inclusions and equivalence of basic objects.} \label{describe}
\begin{parrafo} \label{describe1}
There are two types of information encoded in the notion of basic object:  A pair \((W,E)\), and a closed
subset \(F=\Sing(J,b)\). Moreover,  any sequence of transformation of basic objects
\begin{equation}
\label{ressolbo} (W,(J,b),E)=(W_0,(J_0,b),E_0)\longleftarrow (W_1,(J_1,b),E_1) \longleftarrow \ldots
\longleftarrow (W_k,(J_k,b),E_k),
\end{equation}
induces a sequence of transformations of pairs with the same centers,
\begin{equation}
\label{resol2000}
\begin{array}{ccccccc}
(W_{0},E_{0}) & \longleftarrow & (W_{1},E_{1})& \longleftarrow & \cdots & \longleftarrow & (W_{k}, E_{k})
\end{array}
\end{equation}
together with a collection of closed sets \(F_i=\Sing(J_i,b)\) included in \(W_i\), for
\(i=0,1,\ldots,k\).
\end{parrafo}

\begin{parrafo}
\label{inclusions}  {\bf Inclusions of basic objects.} We say  that \((W,(L,b),E)\) {\em is contained} in
\((W,(J,b),E)\),
\begin{equation*}
(W,(L,s),E) \subset (W,(J,b),E)
\end{equation*}
if $\Sing(L,s)\subset \Sing(J,b)$, and
\begin{enumerate}
\item[(i)] {\em For any}   sequence of  transformations  of the basic
object \((W,(L,s),E)\),
\begin{equation*}
(W,(L,s),E) \longleftarrow (W_1,(L_1,s),E_1) \longleftarrow \ldots \longleftarrow (W_k,(L_k,s),E_k),
\end{equation*}
the corresponding induced sequence of transformations of pairs as in (\ref{resol2000}) defines a
permissible sequence
\begin{equation*}
(W,(J,b),E) \longleftarrow (W_1,(J_1,b),E_1) \longleftarrow \ldots \longleftarrow (W_k,(J_k,b),E_k),
\end{equation*}
with the additional property that for $ i=1,\ldots,k$,
\begin{equation*}
       \Sing(L_i,s) \subset \Sing(J_i,s).
\end{equation*}

\item[(ii)] If \(U \subset W\) is an open set then we require the
previous property to hold for the restrictions \((W,(L,s),E)_U\) and \((W,(J,b),E)_U\) for any \(U\) (see
\ref{propalg}(2)).
\end{enumerate}
\end{parrafo}

\begin{example} \label{yqs}
Fix a pair \((W,E)\) and two sheaves of ideals \(J \subset J^{\prime} \subset {\mathcal O}_{W}.\) Then
clearly $$(W,(J^{\prime},b),E) \subset (W,(J,b),E).$$
\end{example}

\begin{remark} \label{idea}
     We have just introduced a notion of inclusion of basic objects.
It may also happen that there is a double inclusion, namely that \((W,(L,s),E)\) is contained in
\((W,(J,b),E)\) and \((W,(J,b),E)\) is contained in \((W,(L,s),E)\) (i.e. $id: (W,E) \to (W,E) $ defines
an isomorphism of both basic objects as in Definition \ref{isoBO}). Then we will say that they are {\em
equivalent}, and will use the notation
     \begin{equation*}
(W,(J,b),E)\cong (W,(L,s),E).
\end{equation*}
\end{remark}

\begin{remark} \label{idea2}
Note that we  view a basic object as a pair, and a way of defining transformations of pairs and closed
sets: Two equivalent basic objects $(W,(J,b),E)\cong(W,(L,s),E)$ define exactly the same transformations
of pairs and the same closed sets. The invariants involved in the constructive resolution of these basic
objects (the functions $f_i$ introduced in Definition  \ref{AlgResol}), will be defined in terms of the
closed subsets described by the basic objects. So the invariants attached to \((W,(J,b),E)\) will
coincide with those attached to \((W,(L,s),E)\). In particular the constructive resolution of both basic
objects given by Definition  \ref{AlgResol} will be the same (i.e. will be defined by the same sequence of
transformations).
\end{remark}

\begin{example}\label{exampleequiv}
If we set \((L,s)=(J^{2},2b)\), then
\begin{equation*}
(W,(J,b),E)\cong (W,(L,s),E),
\end{equation*}
(see also Example \ref{exampleequiv1}).

\end{example}

\begin{parrafo}
\label{intersection}  {\bf Intersections of basic objects.} Given two basic objects, $$(W,(J,b),E) \mbox{
} \mbox{ } \mbox{ } \mbox{and}
  \mbox{ } \mbox{ } \mbox{ } (W,(I,c),E)$$
we define its {\em intersection},
        $(W,(J,b),E) \cap (W,(I,c),E)$,
as the   basic object \((W,(K,e),E)\) where \(K=J^c+I^b\) and \(e=bc\).

The intersection \((W,(K,e),E)\) is contained in both \((W,(J,b),E)\) and \((W,(I,c),E)\), in the sense
defined in \ref{inclusions}. So any transform $(W,(K,e),E) \longleftarrow (W_1,(K_1,e),E_1) $ induces
two  transformations $(W,(J,b),E) \longleftarrow (W_1,(J_1,b),E_1) $ and $(W,(I,c),E) \longleftarrow
(W_1,(I_1,c),E_1) $, and it can be  checked that $$(W_1,(K_1,e),E_1) =(W_1,(J_1,b),E_1)  \cap
(W_1,(I_1,c), E_1)$$ so that the definition of intersection is {\em compatible} with transformations.

\end{parrafo}

\begin{example}
\label{particular} Let \((W_0,(J_0,b),E_0)\) be a basic object, let \(H_0\subset W_0\) be a regular
hypersurface and consider the basic object \((W_0,(\cali(H_0),1),E_0)\).  Then, $$(W_0,(J_0,b),E_0)\cap
(W_0,(\cali(H_0),1),E_0)= (W_0,(J_0+\cali(H_0)^b,b),E_0)=(W_0,(K_0,b),E_0).$$ Note that a resolution of
\((W_0,(K_0,b),E_0)\),
\begin{equation}\label{ka}
     (W_0,(K_0,b),E_0)\longleftarrow \ldots \longleftarrow
(W_{n-1},(K_{n-1},b),E_{n-1})\longleftarrow (W_n,(K_n,b),E_n)
\end{equation}
induces a sequence of transformations of basic objects for both $(W_0,(J_0,b),E_0)$
\begin{equation*}
(W_0,(J_0,b),E_0)\longleftarrow \ldots \longleftarrow (W_{n-1},(J_{n-1},b),E_{n-1})\longleftarrow
(W_n,(J_n,b),E_n)
\end{equation*}
and $(W_0,(\cali(H_0),b),E_0)$, $$ (W_0,(\cali(H_0),b),E_0)\longleftarrow \ldots \longleftarrow
(W_{n-1},(\cali(H_{0})_{n-1},b),E_{n-1})\longleftarrow (W_n,(\cali(H_{0})_n,b),E_n). $$ Since \(H_{0}\)
is a regular hypersurface, \(\cali(H_{0})_n\) is the sheaf of ideals defining the strict transform of
\(H_0\), say \(H_{0,n} \subset W_n\) (see Remark \ref{so} (2)), and since (\ref{ka}) is a resolution of
\((W_0,(K_0,b),E_0)\), then
\begin{equation*}
\Sing(J_n,b)\cap H_{0,n}=\emptyset.
\end{equation*}
\end{example}

\begin{parrafo} {\bf Basic objects defined over different pairs.}
\label{incluyrestric} We have defined operations of inclusions and intersections of two basic objects
over a same pair.  Now let \((W,E=\{H_1,..,H_r\})\) be a pair and let  \( Z \subset W\) be a closed
smooth subscheme. Assume that each \(H_i\) intersects \(Z\) transversally at \(\widetilde{H}_i\) for
$i=1,\ldots,r$, and that \((Z, \widetilde{E}=\{\widetilde{H_1},...,\widetilde{H_r}\})\) is a pair (for
instance, this will be the case if \(E=\widetilde{E}=\emptyset\)).

Note that a closed subset in \(Z\) is a closed subset in \(W\) and that any sequence of transformation of
the pair \((Z,\widetilde{E})\) induces one over \((W,E)\) (see Definition \ref{normalc}). Note also that
\begin{equation*}
(Z,(A,e),\widetilde{E}) \subset (W,(\mathcal{I}(Z),1),E),
\end{equation*}
since any sequence of transformations over the first induces one
over the second, with inclusion of the corresponding closed sets.
In this context we can also have two basic objects with a double
inclusion as in Remark \ref{idea} (an inclusion of the
corresponding closed sets), and then we will say they are {\em
equivalent}, and denote it by
\begin{equation*}
(Z,(A,e),\widetilde{E}) \cong (W,(J,b),E).
\end{equation*}
If $\mbox{dim}(Z)=d'$, we will  abuse notation  by saying that
$(W,(J,b),E)$ has a structure of a $d'$-dimensional basic object.
\end{parrafo}

\begin{example} \label{curvintcom}
Let \(W\) be smooth of dimension \(d\), and let \(X \subset W\) be a curve which is defined by \(d-1\)
global sections \(f_1,...,f_{d-1}\) in \(\calo_W\). Assume that \(Z=V(<f_1>)\) is smooth. Then
\begin{equation*}
(W,(J,1), \emptyset) \cong (Z,(A,1),\emptyset)
\end{equation*}
where \(J=<f_1,...,f_{d-1}> \subset \calo_W\) and \(A=J \calo_Z\).
\end{example}




\section{The differential operator  and the inductive nature of the
algorithm.} \label{SeccDelta} The objective of this section is to state Proposition \ref{ordenuno},
    which is a key step in the inductive proof of Theorem
    \ref{resolucionob}.  Proposition \ref{ordenuno} says that  we can
    associate (at least locally) a $d-1$ dimensional basic object to
    some class of $d-$dimensional basic objects (we will refer to them
    as  {\em basic objects within the case $\ord_{\calB}\cong 1$} -see
    Definition \ref{withincase1}).  The rest of this section is devoted
    to giving a flavor  of the idea behind of Proposition \ref{ordenuno}
    as well as to explain the role of the sheaf of differentials in
    this process.  The proof of Proposition \ref{ordenuno} will be
    postponed to Section \ref{proposition}. Note that all arguments
exposed here
    extend to the class of schemes described in Definition  \ref{laclasescero}.

\begin{definition}
\label{def1} Given a sheaf of ideals \(J \subset {\mathcal O}_W\) we will denote by \(R(1)(F)\)  the
closed subset of points where
     \(F=V(J)\) has codimension 1 in \(W\).
\end{definition}
\begin{definition}
\label{withincase1} Let ${\calB}=(W,(J,b),E)$ be a  basic object, and consider the equivariant function
introduced in Section \ref{examplesfamilies},
\begin{equation*}
\begin{array}{cccc}
\ord_{\calB}: & F=\Sing(J,b) & \to  &  {\mathbb Q} \\ & x & \to &
      \frac{\nu_J(x)}{b}.
\end{array}
\end{equation*}
We will say that ${\calB}$ is within the {\bf case} \(\ord_ {\calB}\cong 1\) if the maximum value of
\(\ord_{\calB}: F \to {\mathbb Q}\) is 1.
\end{definition}

\begin{proposition}
\label{ordenuno} Let \({\calB}=(W,(J,b),E)\) be a d-dimensional basic object, assume that \(E=\emptyset\)
and that \( \ord_{ B}\cong 1\).  Then:
\begin{enumerate}
\item[a)] There is an open covering \(\{U_{\alpha} \}_{\alpha\in
    {\Lambda}}\), and for each index \(\alpha\in {\Lambda} \) a closed
    and smooth hypersurface \(\widetilde{W}_{\alpha}\subset
    U_{\alpha}\), such that if \((U_{\alpha},(J_{\alpha},b),
    \emptyset)\) is the restriction of \((W, (J,b),\emptyset)\) to
    \(U_{\alpha}\) as in Remark \ref{DefRestricBO},  then
\begin{equation*}
(U_{\alpha},(J_{\alpha},b), \emptyset) \subset (U_{\alpha},(\mathcal{I}(\widetilde{W}_{\alpha}),1),
\emptyset).
\end{equation*}
\item[b)] If \(U_{\alpha} \cap R(1)(\Sing (J,b))=\emptyset \) then
$(U_{\alpha},(J_{\alpha},b), \emptyset)$ has  structure of $(d-1)$-dimensional basic object, i.e. there
is a $(d-1)$-dimensional basic object \((\widetilde{W}_{\alpha},({\mathcal A}_{\alpha},e_{\alpha}),
\emptyset)\) such that $$(U_{\alpha},(J_{\alpha},b), \emptyset)\cong (\widetilde{W}_{\alpha},({\mathcal
A}_{\alpha},e_{\alpha}), \emptyset).$$
\end{enumerate}
\end{proposition}

\begin{remark}
\label{enoempty}

This proposition is the key point of our inductive argument. We are using the notation of
\ref{incluyrestric} (here \(E=\emptyset\)). In general, when \(E \neq \emptyset\), our constructive
resolution will not follow only from induction on the dimension of ambient space \(W\), but rather from
induction on pairs \((W,E)\). The function \( n(x) \)  introduced in \ref{nnn}  will allow us to treat
the case \( (W,E) \) when \( E\neq\emptyset \).
\end{remark}

\begin{definition}
\label{delta1} Let  \(W\) be  a smooth scheme over a field \({\bf k}\) and let $\Omega_{W/{\bf k}}$ be
the (locally free) sheaf of differentials of $W$ over \({\bf k}\).  Given a closed point \( \xi\in W \)
and a regular system of parameters \( \{x_{1},\ldots,x_{n}\}\subset\mathcal{O}_{W,\xi} \), let
\begin{equation*}
         \left\{\frac{\partial\ }{\partial x_{1}},\ldots,
         \frac{\partial\ }{\partial x_{n}}\right \}
\end{equation*}
be a basis of the  dual of the \( \mathcal{O}_{W,\xi} \)-module, \( (\Omega_{W/{\bf k}}^{1})_{\xi} \).
We will express the differential operator \( \Delta \) in terms of this basis although it will be
independent of any choice of parameters.
\end{definition}

\begin{definition}
\label{delta} (See also Section \ref{completions}).  Let  \(W\) be  a smooth scheme over a field \({\bf
k}\), and let $J\subset {\mathcal O}_{W}$ be a non-zero sheaf of ideals. Given a closed point \( \xi\in W
\) and a regular system of parameters \( \{x_{1},\ldots,x_{n}\}\subset\mathcal{O}_{W,\xi} \), we define \(
\Delta(J)\subset {\mathcal O}_{W}\) as the ideal  locally generated  by
\begin{equation*}
       \left\{f_{i}\right\}_{i=1}^{r}\cup \left\{\frac{\partial
       f_{i}}{\partial
       x_{j}}\right\}_{\genfrac{}{}{0pt}{}{i=1,\ldots,r}{j=1,\ldots,n}}
\end{equation*}
  where
$J_{\xi}=\langle f_1,\ldots,f_r\rangle$. We also define by induction $\Delta^i(J)=\Delta
\left(\Delta^{i-1}(J)\right)$, for $i\geq 1$, where $\Delta^0(J)=J$ and $\Delta^1(J)=\Delta(J)$. Note
that $\Delta(J)$ is independent of the choice of parameters.
\end{definition}

\begin{parrafo}
\label{PropOrdDelta} {\bf Properties of the \(\Delta\) operator.}
      If \(J\) is a non-zero sheaf of ideals  at any irreducible
component of the smooth scheme \(W\) of characteristic zero,  then the ideal \(\Delta(J)\) satisfies  the
following properties:
\begin{enumerate}
\item[{\bf P1.}]  \( J\subset\Delta(J)\subset\cdots\subset
         \Delta^{b}(J)=\calo_{W} \) for some \( b \in \nat\).
\item[{\bf P2.}]  For any point \(\xi \in W\), \(
         \nu_{\calo_{W,\xi}}(J_{\xi})=b>1 \) if and only if  \(
         \nu_{\calo_{W,\xi}}(\Delta(J)_{\xi})=b-1 \).  In particular \(
         \nu_{\calo_{W,\xi}}(J)=b>0 \) if and only if  \(
         \nu_{\calo_{W,\xi}}({\Delta}^{b-1}(J))=1 \).
\item[{\bf P3.}] For any  point \(\xi\in W\),
\(\nu_{\calo_{W,\xi}}(J)\geq b>1\) if and only if  \(\xi\in
V(\Delta^{b-1}(J))\).
\end{enumerate}
\end{parrafo}

\begin{parrafo}
\label{deinduc} {\bf Basic objects and the $\Delta$ operator.} Let \( (W,(J,b),E)\) be a basic object
over a field \({\bf k}\) of characteristic zero, then by looking at the Taylor expansion of the
generators of $J$ we have that
\begin{equation}
\label{zero} V(\Delta^{b-1}(J))=\Sing(J,b).
\end{equation}
See also \ref{propertyd}.
\end{parrafo}

\begin{remark}
\label{withincase1oh} Given a basic object ${\calB}=(W,(J,b),E)$, set $F=\Sing(J,b)$, and note that the
following statements are equivalent:
\begin{enumerate}
\item[(a)] The maximum value of \(\ord_{\calB}: F \to {\mathbb Q}\) is
1 (i.e. $B$ is within the case $\ord_{\calB}\cong 1$).
\item[(b)] For every \(x \in F\),   \(\ord_{\calB}(x)=1\).
\item[(c)] The order of \(J_x\) at \({\mathcal O}_{W,x}\) is \(b\) for
every \(x\in F\).
\item[(d)] The order of \(\Delta^{b-1}(J)\) is zero or one at every
point   of \(W\) and the order is one at points of \( F \).
\end{enumerate}
\end{remark}

\begin{definition} \label{DefCoeff}
Let  \( (W,(J,b),E) \) be a basic object, and let \( Z\subset W \) be a closed and smooth subscheme such
that
\begin{enumerate}
       \item  For any \( H\in E \), the closed subschemes \( Z \) and
       \( H \) intersect transversally.

       \item The hypersurfaces of \( Z \), \( E\cap Z=\{H\cap Z\mid
       H\in E\} \), have only normal crossings.
\end{enumerate}
Then we define the {\em coefficient ideal}
        \begin{equation}
\label{coeficientes}
      \Coeff_{Z}(J)= \sum_{i=0}^{b-1}\left(\Delta^{i}(J)\mathcal{O}_{Z}
      \right)^{\frac{b!}{b-i}}\subset \calo_Z.
        \end{equation}
If the coefficient  ideal is non-zero, we also define the basic object
        \begin{equation}
\label{bocoeficientes}
      (Z,(\Coeff_{Z}(J),b!),E\cap Z).
        \end{equation}
\end{definition}

\begin{remark}
\label{Eempty} Note that the conditions required in Definition
\ref{DefCoeff} are trivially satisfied if \( E=\emptyset \).  We
will show that the ideals
\(\mathcal{A}_{\alpha}=\Coeff_{\widetilde{W}_{\alpha}}(J_{\alpha})\)
will satisfy the conditions in Proposition~\ref{ordenuno}, with
$e=b!$.

\end{remark}

\begin{example}
\label{exchar} Let \( W=\mathbb{A}_{\mathbf{k}}^{3} \) and consider the ideal \(J=\langle
z^2+x^3y^3\rangle\), generated by an equation expressed  in a Tschirnhausen form (following Abhyankar).
Then $$ \Delta(J)=\left\langle z^{2}+x^{3}y^{3}, 2z, 3x^{2}y^{3}, 3x^{3}y^{2}\right\rangle,$$ and if the
characteristic of \( \mathbf{k} \) is not two, then \( \Delta^{2}(J)=\mathcal{O}_{W} \).  Note that
\begin{equation*}
         \Sing(J,2)=V(\Delta(J))=V(x,z)\cup V(y,z),
\end{equation*}
and that  \(\ord_{\calB} \cong 1\). In fact, if $Z=V(\langle z\rangle )$ and ${\mathcal A}=\langle
(x^3y^3)\rangle$, then
\begin{equation}
\label{believe} (W,(J,2),\emptyset)\cong(Z, ({\mathcal A},2),\emptyset),
\end{equation}
since the equation has a Tschirnhausen form. So
$(W,(J,2),\emptyset)$ has a structure of two-dimensional basic
object (note that ${\mathcal A}=\Coeff_{Z}(J)$). Although we will
not prove this equivalence here (see Proposition
\ref{PropBasInduc}),  we already point out that this is the key
step in the induction argument behind the proof of Proposition
\ref{ordenuno} and hence in the proof of Theorem
\ref{resolucionob} (we refer the reader to Section
\ref{proposition} for a complete proof that an equivalence of the
form  (\ref{believe}) holds in general).
\end{example}

\begin{remark}
\label{despuesde} Equivalence (\ref{believe}) can be established
thanks to the fact that $Z\subset W$ is a smooth hypersurface
defined by an ideal $I(Z) \subset \Delta^{b-1}(J)$. In general,
given a basic object ${\calB}=(W,(J,b),E)$ such that $\ord_{\calB}
\cong 1$, we will not be able to find a smooth hypersurface
$\widetilde{W}\subset W$ such that globally,
\begin{equation}  \label{forminduc}
\mathcal{I}(\widetilde{W})\subset \Delta^{b-1}(J).
\end{equation}
We shall prove  that a smooth hypersurface \(\widetilde{W}\) that
fulfills  condition (\ref{forminduc}), also has the property that
\begin{equation} \label{forminduccc}
(W,(J,b),E) \subset (W, (\mathcal{I}(\widetilde{W}),1),E),
\end{equation}
(see \ref{corgir}).  Furthermore, we shall proof that if
\({\mathcal B}= (W,(J,b),E) \longleftarrow {\mathcal
B}_1=(W_1,(J_1,b),E_1)\) is a transformation as in
\ref{transformation} or as in \ref{projections}, then
\(\mathcal{I}(\widetilde{W})_1\subset \Delta^{b-1}(J_1)\). Here
\(\mathcal{I}(\widetilde{W})_1=\mathcal{I}(\widetilde{W_1})\)
where \(\widetilde{W_1}\) is the strict transform of
\(\widetilde{W}\) in the first case (see Remark \ref{so} (2)), and
the pull-back in the second case. So, in either case \(
\Delta^{b-1}(J_1)\) has order at most one, so
$\Sing(J.b)=\emptyset$ or \(\ord_{{\mathcal B}_1}\cong 1 \).

The key point for the proof of Proposition \ref{ordenuno} is the
inclusion (\ref{forminduccc}). The role of (\ref{forminduc}) is
technical, but provides a way of finding hypersurfaces for which
\ref{forminduccc} holds. This fact naturally leads us to the
notion of {\em general basic object}, which will be introduced in
Section~\ref{sectiongbo}. The following simple example will
illustrate this notion, with $b=1$, and hence $\Delta^{b-1}(J)=J$.
\end{remark}

\begin{example}\label{ejemcurva}
Let \(W\) be smooth of dimension \(d\), and let \(X \subset W\) be a smooth curve. One can always define
an open covering of \(W\), \(\{U_{\alpha}\}_{\alpha\in {\Lambda}}\), such that the restriction of the
curve, \(X_{\alpha}=X \cap U_{\alpha}\), is a complete intersection on each open subset $U_{\alpha}$,
i.e., $$ \cali(X_{\alpha})=\langle f_1,\ldots,f_{d-1}\rangle=J_{\alpha}. $$ Let $E=\emptyset$ and let
\((U_{\alpha},(J_{\alpha},1),E_{\alpha})\) be the  restriction of \((W,(J=\cali (X),1),E)\) to
\(U_{\alpha}\) and   define \(Z_{\alpha}=V(f_1)\). Then
\begin{equation*}
(U_{\alpha},(J_{\alpha},1),E_{\alpha}) \subset (U_{\alpha},(f_{1},1),E_{\alpha}),
\end{equation*}
and
\begin{equation*}
(U_{\alpha},(J_{\alpha},1), \emptyset) \cong (Z_{\alpha},(\Coeff(J_{\alpha}),1),\emptyset).
\end{equation*}
As in Example \ref{curvintcom}, both basic objects describe the same
    closed set after any sequence of transformations.  In particular, a
    resolution of the basic object
    \((U_{\alpha},(J_{\alpha},1),\emptyset )\) is equivalent to a
    resolution of \(( Z_{\alpha},(\Coeff(J_{\alpha}),1),\emptyset)\).
    The advantage of working with the second basic object is that it is
    a \((d-1)\)-dimensional basic object.  This already illustrates
    that our form of induction on the dimension is of local nature, and also
    leads to the consideration of closed sets which are
    {\em locally} described by basic objects (see Section
    \ref{sectiongbo} where the notion of general basic objects is
    introduced).
\end{example}

\begin{remark} \label{constraint}
Note that in Example~\ref{ejemcurva} the closed sets \( \Sing(\Coeff(J_{\alpha}),1) \) patch to the curve
\( X=\Sing(J,1) \).

Example \ref{exchar} illustrates why we require that the characteristic of the underlying field be zero,
at least if we want condition (d) of Remark \ref{withincase1oh} among the equivalent conditions in case
\(\ord_{\calB}\cong 1\).  It is actually this point which imposes the constraint on the characteristic
for the resolution of singularities of algebraic varieties, because it is not possible, in general, to
use the argument of reduction in the dimension that we are using here.
\end{remark}

\begin{remark}\label{deladesig}
Set $ W \longleftarrow W_1$ a monoidal transformation with center
$C$, let $H_1$ be the exceptional locus, and let $J \subset
\calo_W$ be a sheaf of ideals with order $b$ along points in $C$.
Then $J \calo_{W_1}=I(H_1)^b\cdot J_1$ and the order of $J_1$ is
{\em at most} $b$  at the exceptional points of $W_1$ (in $H_1
\subset W_1$). This is quite a general result, but in our context
(over a field of characteristic zero) it also follows from Remark
\ref{despuesde}: $\Delta^{b-1}(J)$ has order one along the points
in $C$, and  hence, as indicated in Remark  \ref{despuesde},  $
\Delta^{b-1}(J_1)$ has at most order one along exceptional points;
so the claim on $J_1$ follows from \ref{PropOrdDelta}, P2).
  \end{remark}



\section{General basic objects and the main inductive theorem:
   Theorem \ref{Theoremd}.} \label{sectiongbo} In the previous section
we have already observed the need to generalize the notion of basic object in order to use an inductive
argument for the proof of Theorem \ref{resolucionob}. This leads us to the notion of {\em general basic
   objects} which will be developed in this section. Inclusion,
intersection and equivalence of basic objects will naturally extend to this new context, as well as the
notion of equivariant functions (see Proposition \ref{wordngbo}). At the end of this section we extend
Theorem \ref{resolucionob}, in terms of general basic objects, with a {\em constructive} formulation
(\ref{remarkresol}), in Theorem \ref{Theoremd}. The proof of Theorem \ref{Theoremd} will be developed in
Sections \ref{proofbo} and \ref{proofofgbo}. Finally, in Section \ref{SeccConsFunc}, we introduce the
functions \( f_{i}^{d} \) so that the resolution of basic objects, provided by Theorem \ref{Theoremd},
coincides with that defined by the algorithm described in Definition \ref{AlgResol}.

\begin{definition}
\label{defbo} A {\em \(d\)-dimensional general basic object over a
   pair \((W,E)\)} consists of an open covering of \( W \),
\(\{U_{\alpha}\}_{\alpha \in \Lambda} \) such that if \((U_{\alpha},E_{\alpha})\) is the restriction of
\((W,E)\) to \(U_{\alpha}\), then we have:
\begin{enumerate}

\item[(i)] {\bf A collection of basic objects.} For
   every \(\alpha \in \Lambda\) there is a closed and smooth
   \(d\)-dimensional subscheme \(\widetilde{W}_{\alpha} \subset
   U_{\alpha}\), which intersects transversally all hypersurfaces of
   \(E_{\alpha}\) defining a pair \((\widetilde{W}_{\alpha} ,
   \widetilde{E}_{\alpha})\),
   and a basic object
\begin{equation*}
(\widetilde{W}_{\alpha}, (B_{\alpha},d_{\alpha}),\widetilde{E}_{\alpha}).
\end{equation*}
Obviously, for each index \(\alpha\) the closed set
\begin{equation*}
\Sing(B_{\alpha},d_{\alpha})\subset U_{\alpha}
\end{equation*}
     is locally
closed in  \(W\) (compare to conditions in \ref{DefCoeff}).

\item[(ii)] {\bf A patching condition.}
There is a closed subset \(F \subset W\) such that
\begin{equation*}
F \cap U_{\alpha}=\Sing(B_{\alpha},d_{\alpha})
\end{equation*}
for every \(\alpha\in \Lambda\).

\item[(iii)] {\bf Stability of patching (I).} Let
\begin{equation*}
(W,E) \longleftarrow (W_1,E_1)
\end{equation*}
be  a permissible transformation   with center \(Y\subset F\), let \(\{U_{\alpha,_1}\}\) be  the pullback
of   \(\{U_{\alpha}\}_{\alpha\in  \Lambda}\)  to  \(W_1\), and for each \(\alpha\in  \Lambda\) let
\begin{equation*}
(\widetilde{W}_{\alpha}, (B_{\alpha},d_{\alpha}), \widetilde{E}_{\alpha})\longleftarrow
(\widetilde{W}_{\alpha,1}, (B_{\alpha,1},d_{\alpha}), \widetilde{E}_{\alpha,1}).
\end{equation*}
be the corresponding transformation of basic objects. Then there is a closed set \(F_1 \subset W_1\) so
that
\begin{equation*}
F_1 \cap U_{\alpha,1}=\Sing(B_{\alpha,1},d_{\alpha})
\end{equation*}
for each index \(\alpha \in
  \Lambda\).

\item[(iv)] {\bf Stability of patching (II).}
Let \(W \leftarrow W_1\) be a projection and let
\begin{equation*}
(W,E) \longleftarrow (W_1,E_1)
\end{equation*}
be the  corresponding transformation of pairs as defined in \ref{projections}. Let \(\{U_{\alpha,1}\}\)
be  the pullback of   \(\{U_{\alpha}\}_{\alpha\in I}\)  to  \(W_1\), and for each \(\alpha\in \Lambda\)
let
\begin{equation*}
(\widetilde{W}_{\alpha}, (B_{\alpha},d_{\alpha}), \widetilde{E}_{\alpha})\longleftarrow
(\widetilde{W}_{\alpha,1}, (B_{\alpha,1},d_{\alpha}), \widetilde{E}_{\alpha,1}).
\end{equation*}
be the transformation of basic objects as defined in \ref{projections}. Then there is a closed set \(F_1
\subset W_1\) such that
\begin{equation*}
F_1 \cap U_{\alpha,1}=\Sing(B_{\alpha,1},d_{\alpha})
\end{equation*}
for each index \(\alpha \in  \Lambda\).

\item[(v)] {\bf Stability of patching (III). } The   patching condition
defined in (iii) and (iv)   holds after any sequence of transformations: Given a sequence of
transformations of pairs, $$
\begin{array}{ccccccccc}
(W_0,E_0) &\longleftarrow & (W_1,E_1) & \longleftarrow & \ldots & \longleftarrow & (W_r,E_r) &
\longleftarrow  & (W_{r+1},E_{r+1}) \\ \cup & & \cup & & & & \cup & & \\ F_0 & & F_1 & & & & F_r &
\end{array}
$$ where for $i=0,1,\ldots,r$,  $W_{i+1}\to W_i$ is defined either by:
\begin{enumerate}
\item[(1)] blowing up  at centers \(Y_i \), permissible for the pair
\((W_i,E_i)\), and \(Y_i\) included in the inductively defined closed sets \(F_i \subset W_i\), or
\item[(2)] a projection \(p: W_{i+1} \to W_i\),
\end{enumerate}

\noindent there is an  open covering \(\{U_{\alpha,r+1}\}\) of \(W_{r+1}\) (the pull back of
\(\{U_{\alpha}\}\)), a sequence of transformations of basic objects,
\begin{multline} \label{transalpha}
(\widetilde{W}_{\alpha}, (B_{\alpha},d_{\alpha}), \widetilde{E}_{\alpha})\longleftarrow
(\widetilde{W}_{\alpha,1}, (B_{\alpha,1},d_{\alpha}),\widetilde{E}_{\alpha,1}) \longleftarrow\cdots \\
\cdots\longleftarrow (\widetilde{W}_{\alpha,r+1},(B_{\alpha,r+1},d_{\alpha}), \widetilde{E}_{\alpha,r+1}),
\end{multline}
    and a closed set
$F_{r+1}\subset W_{r+1}$,  such that for each $\alpha\in \Lambda$,
\begin{equation*}
F_{r+1} \cap U_{\alpha,r+1}=\Sing(B_{\alpha,r+1},d_{\alpha}).
\end{equation*}

\item[(vi)] \textbf{Restriction to open sets.}
If \( V\subset W \) is an open set, consider the restriction of all data to \( V \): The open covering \(
\{U_{\alpha}\cap V\}_{\alpha\in\Lambda} \), the basic objects \( (\widetilde{W}_{\alpha},
(B_{\alpha},d_{\alpha}),\widetilde{E}_{\alpha})_{V} \)
  and the closed set \( F_{V}=F\cap V \). Then
we require that all properties (i), (ii), (iii) (iv) and (v) hold for the restriction (see Remark
\ref{DefRestricBO}).
\end{enumerate}
A general basic object will be denoted by \(({\mathcal F},(W,E))\), the restriction to an open set \( V
\) will be
  denoted by
\( (\mathcal{F}_{V},(V,E_{V})) \), and we will write a  sequence of transformations and projections as
\begin{equation}
\label{trangbo}
\begin{array}{ccccccc}
({{\mathcal F}}_0,(W_0,E_0)) & \longleftarrow & \ldots & \longleftarrow & ({{\mathcal
F}}_{r},(W_{r},E_{r})) & \longleftarrow & ({{\mathcal F}}_{r+1},(W_{r+1},E_{r+1}))\\ \cup & & & & \cup &
& \cup\\ F_0 & & & & F_{r} & & F_{r+1}.
\end{array}
\end{equation}
\end{definition}

\begin{remark} \label{notaopen}
If \( (\mathcal{F},(W,E)) \) is a general \( d \)-dimensional
basic object, note that \( d \) can be different from \( \dim{W}
\). A basic object \( (W,(J,b),E) \) defines a general basic
object \( (\mathcal{F},(W,E)) \), with the trivial open covering
and in this case \( d=\dim{W} \). A general basic object can be
described by giving two different open coverings. What is
important here are the closed sets \(F\) that it defines. That is
why in the notation for general basic objects \(({\mathcal
F},(W,E))\) there is no reference to the open covering which
appears in the definition (see Definition \ref{isoGBO}, where the
notion of isomorphism of general basic objects is stated).
\end{remark}

\begin{example}\label{examplegbo1}
Example \ref{ejemcurva} illustrates what a  \(d\)-dimensional general basic object can be. In this case,
the collection of locally closed basic objects is given by $$(U_{\alpha},(J_{\alpha},1),E_{\alpha}),$$
and the closed subset $F\subset W$ of condition (ii) in Definition \ref{defbo} is $X$. Now, consider a
sequence of transformations of pairs, $$
\begin{array}{ccccccccc}
     (W,E) & \longleftarrow &
(W_1,E_1) & \longleftarrow & \ldots & \longleftarrow & (W_r,E_r) & \longleftarrow & (W_{r+1},E_{r+1})\\ F
& & F_1 & & & & F_r & & F_{r+1},
\end{array}
$$ where each transformation is either as in Definition \ref{defbo} (iii)  or as in Definition
\ref{defbo} (iv). Then $F_{r+1}$ is either the  strict transform of \(F_{r}\), in the first case (this is
a particular feature discussed in \ref{so} (2)),  or the pull-back of \(F_{r}\) in $W_{r+1}$, in the
second case.
\end{example}

Most definitions on basic objects extend naturally to general basic objects, and we extend only some of
them. In Proposition \ref{wordngbo} we extend the  two equivariant  functions introduced in Section
\ref{examplesfamilies}.

\begin{definition}\label{resgbo}
A resolution of a general basic object \(({{\mathcal F}}_0,(W_0,E_0))\) is a sequence of transformations
as in (\ref{trangbo}) which fulfills the following two conditions:
\begin{enumerate}
\item[(i)] The sequence involves only transformations as the ones
   introduced  in Definition
\ref{defbo} (iii).
\item[(ii)] The closed set \(F_{r+1}\) is empty.
\end{enumerate}

Note that if $\{U_{\alpha}\}$ is an open covering of $W$ as in Definition \ref{defbo}, then for any \(
\alpha \)  we obtain a resolution of the basic object \(
(\widetilde{W}_{\alpha},(B_{\alpha},d_{\alpha}),\widetilde{E}_{\alpha})\)
  in the sense of Definition \ref{defResol}.
\end{definition}

\begin{definition}\label{isoGBO}
Let \(({{\mathcal F}_0},(W_0,E_0))\) and \(({{\mathcal F}_0}^{\prime}, (W_0^{\prime},E_0^{\prime}))\) be
two general basic objects and let \(\theta: (W, E)\to (W^{\prime},E^{\prime})\) be an isomorphism of pairs
(\ref{DefActPar}). We will say that \(\theta\) induces an {\em isomorphism of general basic objects},
\begin{equation*}
\theta: ({{\mathcal F}_0},(W_0,E_0))\to ({{\mathcal F}_0}^{\prime}, (W_0^{\prime},E_0^{\prime})),
\end{equation*}
if the following conditions hold:
\begin{enumerate}
\item[(i)] The isomorphism \(\theta: W\to W^{\prime}\) induces an
isomorphism of the closed subsets defined by the two  general basic objects
\begin{equation*}
\theta:F \cong F^{\prime}.
\end{equation*}
\item[(ii)] If
\begin{equation*}
\begin{array}{ccccccc}
({{\mathcal F}}_0,(W_0,E_0)) & \longleftarrow & \ldots & \longleftarrow & ({{\mathcal
F}}_{r-1},(W_{r-1},E_{r-1})) & \longleftarrow & ({{\mathcal F}}_{r},(W_{r},E_{r}));\\ \cup & & & & \cup &
& \cup\\ F_0 & & & & F_{r-1} & & F_r
\end{array}
\end{equation*}
is sequence of transformations of general basic objects, and
\begin{equation*}
\begin{array}{ccccccc}
(W,E) & \longleftarrow & (W_1, E_1)& \longleftarrow & \ldots & \longleftarrow & (W_k, E_k)
\end{array}
\end{equation*}
is the corresponding sequence of transformation of pairs, then:
\begin{enumerate}
\item[(a)] The  sequence induced by $\theta$ as in Remark \ref{lifting1},
say
\begin{equation*}
\begin{array}{ccccccc}
(W^{\prime},E^{\prime}) & \longleftarrow & (W_1^{\prime}, E_1^{\prime})& \longleftarrow & \ldots &
\longleftarrow & (W_k^{\prime}, E_k^{\prime}),
\end{array}
\end{equation*}
defines a sequence of transformation of general basic objects,
\begin{equation*}
\begin{array}{ccccccc}
({{\mathcal F}}_0^{\prime},(W_0^{\prime},E_0^{\prime})) & \longleftarrow & \ldots & \longleftarrow &
({{\mathcal F}}_{r-1}^{\prime},(W_{r-1}^{\prime},E_{r-1}^{\prime})) & \longleftarrow & ({{\mathcal
F}}_{r}^{\prime},(W_{r}^{\prime},E_{r}^{\prime}))\\ \cup & & & & \cup & & \cup\\ F_0^{\prime} & & & &
F_{r-1}^{\prime} & & F_r^{\prime}.
\end{array}
\end{equation*}
\item[(b)] The isomorphisms
\(\theta_i: (W_i,E_i)\longrightarrow (W_i^{\prime},E_i^{\prime})\), presented  in Definition
\ref{lifting1}, induce an isomorphism of the closed subsets defined by the general basic objects,  i.e.,
     \(\theta_i(F_i)=F_i^{\prime} \mbox{ \ for \ } i=0,1,\ldots,r+1\).
\end{enumerate}

\item[(iii)] For any open set \( U\subset W_{0} \), set \( U'=\theta(U) \) and
consider the restrictions \( (\mathcal{F}_{0,U},(U,E_{0,U})) \) and \(
(\mathcal{F}'_{0,U'},(U',E'_{0,U'})) \) (see Definition \ref{defbo} (vi)). Then we require that
properties (i) and (ii) hold for the restrictions.
\end{enumerate}
\end{definition}

\begin{definition}
\label{familyequivbo} Assume that for any general  basic object \(({\mathcal F},(W,E))\) there is an
upper semi-continuous function
\begin{equation*}
f_{\mathcal{F}}: F\to (T, \geq)
\end{equation*}
associated to it, this is what we call a family of functions with values on \( T \).  We will say that
the family of functions \(f_{\mathcal{F}}\) is {\em equivariant} if:
\begin{enumerate}
\item[(i)] For any isomorphism of general basic objects
\(\Theta:({\mathcal F},(W,E))\to (\mathcal{F}^{\prime},(W^{\prime},E^{\prime}))\) we have that
\begin{equation*}
f_{\mathcal{F}^{\prime}}\circ\Theta=f_{\mathcal{F}}.
\end{equation*}

\item[(ii)] If $U\subset W$ is an opens set, let \(
(\mathcal{F}_{U},(U,E_{U})) \) be the restriction as in Definition \ref{defbo} (vi). Then the function
\(f_{\mathcal{F}_U }\) is the restriction of \(f_{\mathcal{F}}\).
\end{enumerate}
\end{definition}

\begin{proposition}
\label{wordngbo} Let \(({\mathcal F},(W,E))\) be a general basic object, let $\{U_{\alpha}\}_{\alpha\in
\Lambda}$ be the corresponding open covering of $W$ and let
$(\widetilde{W}_{\alpha},(B_{\alpha},d_{\alpha}),\widetilde{E}_{\alpha})$ be the collection of
$d$-dimensional basic objects  associated to \(({\mathcal F},(W,E))\). Then  the functions
\begin{equation*}
\ord^{d}_{\alpha}: \Sing(B_{\alpha},d_{\alpha}) \to {\mathbb Q} \mbox{ } \mbox{ } \mbox{ and }\mbox{
}\mbox{ } n^{d}_{\alpha}: \Sing(B_{\alpha},d_{\alpha}) \to {\mathbb Z}
\end{equation*}
patch so as to define equivariant functions
\begin{equation*}
\ord^{d}_{\mathcal{F}}:F \to {\mathbb Q} \qquad\text{and}\qquad n^{d}_{\mathcal{F}}: F \to{\mathbb Z}
\end{equation*}
which verify the requirements of Definition~\ref{familyequivbo}.
\end{proposition}

We sketch the proof below, which parallels the proof of Lemma \ref{propertiesnord} for basic objects, and
we refer to Section~\ref{Hironaka} for technical details.

\begin{proof}
Because of the  way they are defined, any transformation of general basic objects
\begin{equation}
\label{trangbo1}
\begin{array}{ccccccc}
({{\mathcal F}}_0,(W_0,E_0)) & \longleftarrow & \ldots & \longleftarrow & ({{\mathcal
F}}_{r},(W_{r},E_{r})) & \longleftarrow & ({{\mathcal F}}_{r+1},(W_{r+1},E_{r+1}))\\ \cup & & & & \cup &
& \cup\\ F_0 & & & & F_{r} & & F_{r+1},
\end{array}
\end{equation}
induces, for each index \(\alpha\),  a sequence of transformations of basic objects
\begin{multline}
\label{transalpha1} (\widetilde{W}_{\alpha}, (B_{\alpha},d_{\alpha}),
\widetilde{E}_{\alpha})\longleftarrow ((\widetilde{W}_{\alpha})_1, ((B_{\alpha})_1,d_{\alpha}),
(\widetilde{E}_{\alpha})_1)\longleftarrow\cdots \\
\cdots\longleftarrow ((\widetilde{W}_{\alpha})_{r+1}, ((B_{\alpha})_{r+1},d_{\alpha}),
(\widetilde{E}_{\alpha})_{r+1}).
\end{multline}
In general, an arbitrary sequence of transformations over
\((\widetilde{W}_{\alpha},
(B_{\alpha},d_{\alpha}),\widetilde{E}_{\alpha})\), will not give
rise to a sequence of transformations of the general basic object.
However, this will be so if we fix \({\alpha}\) together with a
point \(x_0 \in \Sing (B_{\alpha},d_{\alpha})\), and take the
sequence (\ref{transalpha1}) to be an arbitrary \(x_0\)-extendable
sequence (\ref{exextendable}).
\medskip

So assume now that the sequence (\ref{trangbo1})  induces an
\(x_0\)-extendable sequence (\ref{transalpha1}). In this case, for
any index $j=0,1,\ldots,r$, there is an identification of fibers,
\begin{equation*}
(\Sing ((B_{\alpha,j}),d_{\alpha}))_{x_0}=(F_j)_{x_0}.
\end{equation*}

The characterization of the rational number \(\ord^{d}_{\alpha}(x_0)\) (d=dim($\widetilde{W}_{\alpha}$))
in terms of the dimensions of these closed sets (see Lemma \ref{propertiesnord}) shows that, for \(x_0\in
F_0 \cap U_{\alpha} \cap U_{\beta}\), \(\ord^{d}_{\alpha}(x_0)=\ord^{d}_{\beta}(x_0).\) In particular
\(\ord^{d}_{\mathcal{F}}\) is well defined (see also Proposition \ref{HiroTrick}). The property of
equivariance of the function \(\ord^{d}_{\mathcal{F}}\) follows now from (\ref{equifibra}).
\medskip

As for the function \(n^{d}_{\mathcal{F}}\), note that it is well defined by the very definition of
general basic object; the same can be said about the property of equivariance. In fact any
    isomorphism
$$ \theta: ({{\mathcal F}_0},(W_0,E_0))\to ({{\mathcal F}_0}^{\prime}, (W_0^{\prime},E_0^{\prime})) $$
defines an isomorphism $$\theta_0: (W_0, E_0)\to (W^{\prime}_0,E^{\prime}_0)$$ and the functions
\(n_{\mathcal{F}}\) and \(n_{\mathcal{F}^{\prime}}\) are defined in terms of these pairs.
\end{proof}

\begin{theorem}
\label{Theoremd} {\em (Theorem (d))} To each  \( d \)-dimensional general basic object
\((\mathcal{F}_0,(W_0,E_0))\) we can attach a resolution \(R_{\mathcal{F}_0}\) in the sense of Definition
\ref{resgbo},
\begin{equation}
\label{trangbo2}
\begin{array}{ccccccc}
({{\mathcal F}}_0,(W_0,E_0)) & \longleftarrow & \ldots & \longleftarrow & ({{\mathcal
F}}_{r},(W_{r},E_{r})) & \longleftarrow & ({{\mathcal F}}_{r+1},(W_{r+1},E_{r+1}))\\ \cup & & & & \cup &
& \cup\\ F_0 & & & & F_{r} & & F_{r+1}=\emptyset,
\end{array}
\end{equation}
which has the following equivariance property: If
$\Theta:({{\mathcal F}_0},(W_0,E_0))\to ({\mathcal F}_0^{\prime},
(W_0^{\prime},E_0^{\prime}))$ is an isomorphism of general basic
objects, and if
\begin{equation*}
\begin{array}{ccccccc}
R_{\mathcal{F}_0}:\ ({{\mathcal F}}_0,(W_0,E_0)) & \longleftarrow & \ldots & \longleftarrow & ({{\mathcal
F}}_{r-1},(W_{r-1},E_{r-1})) & \longleftarrow & ({{\mathcal
F}}_{r}=\emptyset,(W_{r},E_{r}))\\
\cup & & & & \cup & & \cup \\
  F_0 & & & &
F_{r-1} & & F_r=\emptyset,
\end{array}
\end{equation*}
and
\begin{equation*}
\begin{array}{ccccccc}
R_{\mathcal{F}_0^{\prime}}:\ ({{\mathcal F}}_0^{\prime},(W_0^{\prime},E_0^{\prime})) & \longleftarrow &
\ldots & \longleftarrow & ({{\mathcal F}}_{r-1}^{\prime},(W_{r-1}^{\prime},E_{r-1}^{\prime})) &
\longleftarrow & ({{\mathcal F}}_{r}^{\prime},(W_{r}^{\prime},E_{r}^{\prime}))\\ \cup & & & & \cup & &
\cup\\ F_0^{\prime} & & & & F_{{r^{\prime}}-1}^{\prime} & & F_{r^{\prime}}^{\prime}=\emptyset,
\end{array}
\end{equation*}
are the corresponding resolutions associated to \(\mathcal{F}_0\) and \(\mathcal{F}_0^{\prime}\), then
there is a lifting of \(\Theta\) to isomorphisms
\begin{equation*}
\Theta_i: ({{\mathcal F}_i},(W_i,E_i))\to ({\mathcal F}_i^{\prime}, (W_i^{\prime},E_i^{\prime})), \qquad
\text{(so }\Theta_i(F_i)=F_i^{\prime}) \qquad \text{for}\quad i=1,\ldots, r,
\end{equation*}
or, in other words,  both resolutions are related by an isomorphism
  as in Definition  \ref{isoGBO} (ii).
\end{theorem}

We will show later, that the resolutions \( R_{\mathcal{F}_{0}} \), to be constructed here, will also be
that defined  by the functions \( f_{i}^{d} \) that appear in Definition  \ref{AlgResol}. The fact that
the definition of the functions \( f_{i}^{d} \) is given later, in section~\ref{SeccConsFunc}, is due to
the fact that the invariants involved in the definition of such functions is largely motivated by our
proof of \ref{Theoremd}.

The proof of Theorem \ref{Theoremd} will be presented through the
next two sections: In Section \ref{proofbo} the theorem will be
proven for the case of basic objects, and in Section
\ref{proofofgbo} the general case will be treated. Since the proof
is based in an inductive argument, we already show why Theorem
\ref{Theoremd} holds for $0-$dimensional general basic objects:
namely, in that in that case, each
$(\widetilde{W}_{\alpha},(B_{\alpha},d_{\alpha}),E_{\alpha})$ is
zero dimensional, so we can assume that each
$\widetilde{W}_{\alpha}$ is a point, and hence, each $B_{\alpha}$
is a non-zero ideal in a field. Therefore,
$\Sing(B_{\alpha},d_{\alpha})=\emptyset$, and hence,
$F_0=\emptyset$.




\section{Proof of Theorem \ref{Theoremd} for basic objects.}
\label{proofbo} Giving a basic object ${\mathcal B}=(W,(J,b),E)$ we will distinguish three cases:
\begin{itemize}
\item Case 1:  \(\ord_{\calB} \cong 1\) and $E=\emptyset$
\item Case 2:  \(\ord_{\calB} \cong 1\) and $E\neq\emptyset$.
\item Case 3: The general case: Proof of Theorem \ref{Theoremd} for
basic objects.
\end{itemize}
Cases  1 and 2 are not
  special  cases: They appear as intermediate steps
when proving the general case.
\medskip

\begin{center}
{ Case 1.}
\end{center}
\medskip
We will first study the case  \({ B}=(W,(J,b),\emptyset)\) such that \(\ord_{B} \cong 1\), or a transform
of this case. Induction in this case parallels the treatment of Tschirnhausen transformations
(\ref{exchar}). This is the simplest situation to be considered, but at the same time the most
significant; in fact it shows why centers defined by our procedure are regular, and why the algorithm of
resolution of basic objects is equivariant. The key point in this case is  the fact that the closed set
\(\Sing(J,b)=V(\Delta^{b-1}(J))\) is defined by an ideal of order at most one.
\medskip

\begin{proposition}
\label{case1} Theorem \ref{Theoremd} holds in  Case~1.
\end{proposition}

To prove this statement  we will need some auxiliary results: Lemmas \ref{case12} and \ref{case112}. The
proof of Proposition \ref{case1} will be given in \ref{proofofcase1}. To illustrate the general idea
behind this case, let us consider the following example:

\begin{example}
\label{examplec11} Let \(W= \mbox{Spec}(k[x,y,z])\), let \(J=< z^2+y^3x^2> \subset {\mathcal O}_W =
k[x,y,z]\), \(E= \emptyset\) and  consider the basic object \((W,(J,2),E)\). In this case  \(F=
\Sing(J,2)\) is a union of two lines, and therefore \(R(1)(F)=\emptyset\) (see Definition \ref{def1}).
  Note that \(\langle z\rangle\in\Delta(J)\) (see (\ref{forminduc})).
\medskip

Let \( \widetilde{W}=V(\langle z\rangle) \subset W\). By Proposition \ref{PropBasInduc} we may consider
the two dimensional basic object
\begin{equation*}
    (\widetilde{W},{\mathcal A},\emptyset)=(\widetilde{W},(y^3x^2,2),\emptyset).
\end{equation*}
     This   basic object defines the
2-dimensional structure associated to \((W, (J,2),E)\) in the
sense of \ref{incluyrestric}. In this case note that
\begin{equation*}
     \Sing(<z^2+y^3x^2>, 2)=\Sing (y^3x^2,2),
\end{equation*}
and  furthermore the equality holds after any sequence of
transformations. For example, blow-up the origin in both cases and
compare the closed sets defined by $$ (W,(J,2),
\emptyset)\longleftarrow (W_{1},(J_{1},2),E_1) \mbox{ } \mbox{ }
\mbox{ } \mbox{ and } \mbox{ } \mbox{ } \mbox{ } (\widetilde{W},(
{\mathcal A},2),\emptyset) \longleftarrow
(\widetilde{W}_{1},({\mathcal A}_{1},2),E_1).
$$ Summarizing we have that:

\begin{enumerate}
\item[\(\bullet\)] A resolution of \((W,(J,2),\emptyset)\)  induces a
resolution of \(( \widetilde{W},( {\mathcal A}, 2), \emptyset)\).
\item [\(\bullet\)] A resolution of \((W,(J_1,2),E_1)\) induces  a resolution
of \(( \widetilde{W},({\mathcal A}_1,2), \widetilde{E}_1)\).
\end{enumerate}
\end{example}

\begin{remark} \label{sobrecaso1} Note that:
\begin{enumerate}
\item[(1)]   If \((W,(J,b),E)\) is such that \( \ord\cong 1 \),
so is any transform (see Remark \ref{despuesde}).
\item[(2)] If there is an inclusion of basic objects,
\((W,(K,d),E)\subset (W,(J,b),E)\) and  \((W,(J,b),E)\) is such that \( \ord\cong 1 \), then same holds
for \((W,(K,d),E)\). This follows from \ref{intersection}, $$(W,(K,d),E)\cap (W,(J,b),E)\cong
(W,(K,d),E),$$
  and the fact that the function \(
\ord \) is equivariant (\ref{propertiesnord}).
\end{enumerate}
\end{remark}

\begin{lemma}
\label{case12} Let \(\calB=(W,(J,b),\emptyset)\) be a basic object
such that \(\ord_{{\calB}} \cong 1\). Set
\(F=\Sing(J,b)=V(\Delta^{b-1}(J))\) and let
\(G=\Sing(J,b)\setminus R(1)(F)\) be the union  of the irreducible
components of $\Sing(J,b)$ of codimension greater than or equal to
2. Then:
\begin{enumerate}
\item[(a)] \(R(1)(F)\) and \(G\) are both open and closed in \(F\).
\item[(b)] \(R(1)(F)\) is either empty or regular.
\item[(c)] If \(R(1)(F)\) is not empty, then it is a permissible center for
\((W,(J,b),\emptyset)\) and if
\begin{equation*}
(W,(J,b),\emptyset) \longleftarrow (W_1,(J_1,b),E_1)
\end{equation*}
     is the transformation with center \(R(1)(F)\),  then
\(W_1=W\) and \(R(1)(F_1)= \emptyset \).
\item[(d)] If \(R(1)(F)= \emptyset
\)  then the basic object \((W,(J,b),\emptyset)\) has a structure of \((d-1)\)-dimensional general basic
object.
\end{enumerate}
\end{lemma}

\begin{proof}
Note that if \(x \in R(1)\), then \(J_x= <l^b> \subset \calo_{W,x}\), where \(l \in \calo_{W,x}\) is an
element of order one. This reduces to two observations:
\begin{enumerate}
\item[(i)] If an ideal of order one at  \(\calo_{W,x}\) defines a
hypersurface, the ideal must be principal and generated by an element of order one.
\item[(ii)] If \(Z\) is a smooth irreducible hypersurface included in
\(\Sing(J,b)\), then the order of \(J\) at the generic point if \(Z\) is at least b, and
\(J=\mathcal{I}(Z)^b\cdot J^{\prime} \subset \calo_W\).
\end{enumerate}

These observations prove statements (a), (b) and (c). Statement (d) follows from
Proposition~\ref{ordenuno}, which we reformulate below for a sequence of transformations. Note that if \(
R(1)(F)=\emptyset \) then \( \Coeff_{Z}(J)\neq 0 \) (see Definition \ref{DefCoeff}).
\end{proof}

\begin{lemma}\label{case112}
Let \({\calB}_0=(W_0,(J_0,b),\emptyset)\) be a basic object  such that \( ord_{{\calB}_0} \cong 1\) and
let
     \begin{equation}
\label{equiseq0} {B}_0=(W_0,(J_0,b),E_0=\emptyset) \longleftarrow \ldots \longleftarrow {B}_k=(W_k ,
     (J_k,b),E_k)
\end{equation}
be a sequence of transformations of basic objects. Then:
\begin{enumerate}
\item[(a)] \(R(1)(\Sing (J_k,b))\) is the strict transform of  \(R(1)(\Sing
(J_0,b))\).
\item[(b)] If  \((W_0,(J_0,b),E_0=\emptyset) \subset
(W_0, (\mathcal{I}(\widetilde{W_0}),1),E_0=\emptyset )\) for some smooth hypersurface
$\widetilde{W_0}\subset W_0$, and
\begin{equation*}
(W_0, (\mathcal{I}(\widetilde{W_0}),1),E_0= \emptyset ) \longleftarrow \cdots \longleftarrow (W_k,
(\mathcal{I}(\widetilde{W_k}),1),E_k)
\end{equation*}
is the   sequence induced by (\ref{equiseq0}), then
\begin{equation}\label{agregadotrento00}
(W_k , (J_k,b),E_k)\subset (W_k, (\mathcal{I}(\widetilde{W_k}),1),E_k),
\end{equation}
  and
\(\mathcal{I}(\widetilde{W_k})\) defines a smooth hypersurface \(\widetilde{W_k}\subset W_k\) transversal
to all hypersurfaces of \(E_k\).
\item[(c)] If $U_{0,\alpha}$ and $\widetilde{W_0}_{\alpha}$ denote,
respectively, the open covering of $W_0$ and the corresponding smooth hypersurfaces stated in
Proposition \ref{ordenuno}, then for each index \(\alpha\), sequence (\ref{equiseq0}) induces a sequence
\begin{equation*}
(U_{0,\alpha},(J_{0,\alpha},b),E_{0,\alpha}=\emptyset)\longleftarrow \cdots\longleftarrow
((U_{k,\alpha}),(J_{k,\alpha},b),(E_{k,\alpha}))
\end{equation*}
the  open subsets \((U_{k,\alpha})\) define an open cover of \(W_k\), and for each index \({\alpha}\),
this sequence can be replaced by one in dimension d-1 (\ref{ordenuno} (b)) since
\begin{equation} \label{agregadotrento01}
(U_{k,\alpha},(J_{k,\alpha},b),E_{k,\alpha}) \cong (\widetilde{W}_{k,\alpha},({\mathcal
A}_{\alpha},e_{\alpha}), {\overline E}_{k,\alpha}).
\end{equation}

\end{enumerate}

\end{lemma}

\begin{proof}
Part (b) follows from \ref{inclusions}, and (c) follows from
Proposition~\ref{ordenuno}. Part (a) follows from the proof of
Lemma~\ref{case12}.
\end{proof}
\begin{remark}\label{agregadotrento1}
If only monoidal transformations arise in sequence
(\ref{equiseq0}), then:
\begin{enumerate}
\item[1)] If $R(1)(\Sing(J_k,b))=\emptyset$, then
${\calB}_k=(W_k,(J_k,b), E_k)$ has a $(d-1)$-dimensional structure (Lemma \ref{case12} (d))).
\item[2)] If $R(1)(\Sing(J_k,b))=\emptyset$ and
${\calB}_k=(W_k,(J_k,b),E_k)\cong {\calB}'_k=(W'_k,(J'_k,b'), E_k)$, then ${\calB}'_k$ has a
$(d-1)$-dimensional structure (Lemma \ref{case12} (d)).

\item[3)] Let $(W_k,({\mathcal C},s), E_k)\subset {\calB}_k$.   If $R(1)(\Sing({\mathcal C},s))=\emptyset$ and
$ (W_k,({\mathcal C},s), E_k)\subset {\calB}_k=(W_k,(J_k,b),E_k)$, then $(W_k,({\mathcal C},s), E_k)$ has
a $(d-1)$-dimensional structure (Lemma \ref{case12} (d)).
\end{enumerate}
All three observations follow from the fact that isomorphisms stated in (\ref{agregadotrento01}) require
only  inclusion (\ref{agregadotrento00}) to hold.

\end{remark}

\begin{parrafo}
\label{proofofcase1} {\em Proof of Proposition \ref{case1}:} Assume that ${\calB}=(W,(J,b),E=\emptyset)$
is within Case 1. If \(R(1)(\Sing (J,b))=\emptyset\), then we apply Lemma \ref{case12} (d), and
  the statement follows by induction on
\(d\); in fact we are assuming Theorem~\ref{Theoremd} for \( (d-1) \)-dimensional general basic objects.
\medskip

If \(R(1)(\Sing (J,b))\neq \emptyset\), then by Lemma \ref{case12} (c),  after blowing-up \(R(1)(F)\) we
may assume that the basic object has a structure of \((d-1)\)-dimensional general basic object. Again  we
proceed by induction.
\medskip

If \( \Theta:{\calB}=(W,(J,b),E=\emptyset) \cong {\calB}^{\prime}= (W',(J',b'),E'=\emptyset)\) is an
isomorphism of basic objects, then \({B}^{\prime}\) is within case 1, and there is an isomorphism
$\Theta: W \to W',$ mapping \(F=\Sing(J,b)\) isomorphically to \(F'= \Sing(J',b')\). It is clear that for
{\em any} such isomorphism \(\theta\), $\theta (R(1)(F))=R(1)(F')$. Note also that if \(
R(1)(F)=\emptyset=R(1)(F^{\prime}) \) then \( \Theta \) induces an isomorphism of the corresponding \(
(d-1) \)-dimensional general basic objects (see Lemma \ref{case12} (d)).

By  induction  there is a well defined equivariant resolution
\({\mathcal R}_{{\calB}}\) for any \({\calB}=(W,(J,b),E)\) within
Case 1. Hence  Theorem \ref{Theoremd} holds within this case. \qed
\end{parrafo}

\begin{remark}\label{caso1prima}
Fix  \({\calB}_0=(W_0,(J_0,b),\emptyset)\) within Case 1, and now consider any sequence of permissible
monoidal transformations as in (\ref{equiseq0}),
  not necessarily related to the resolution
\({\mathcal R}_{{\calB}_0}\) defined in \ref{proofofcase1}. We claim now that the basic object
\({\calB}_k=(W_k,(J_k,b),E_k)\) can also be considered and treated as in Case 1 (although $E_k$ is no
longer empty), in the sense that we can attach a resolution to ${\calB}_k$ by the same arguments. In
fact, $F=\Sing(J_0,b)$ is a disjoint union of two closed sets $F=R(1)(F) \cup G$ (Lemma \ref{case12}),
and ${\calB}_0$ has a structure of $d-1$-dimensional general basic object in an open neighborhood of $G$.
Now set $F_k=\Sing(J_k,b)$ and the disjoint union $ F_k= R(1)(F_k)\cup G_k $ as before. By Lemma
\ref{case112} (a), $R(1) (F_k)$ is the strict transform of $R(1)(F)$; and locally at points of $G_k$,
${\calB}_k$ inherits the $d-1$-dimensional structure of ${\calB}_0$ locally at $G$. So the whole point is
to note that $R(1) (F_k)$ is either empty or a permissible center (it has normal crossings with $E_k$),
which follows from the fact that it is the strict transform of $R(1)(F)$, and all irreducible components
of the centers are either disjoint or included in the intermediate strict transforms of this hypersurface.
\end{remark}
\medskip

\begin{center}
{\bf  The general pattern of  Cases 2 and 3.}
\end{center}
\medskip

\begin{parrafo}
\label{pattern} Cases 2 and 3 follow a general pattern which we describe now: For each $d-$dimensional
basic object \({\mathcal B}=(W,(J,b),E)\) we will define  an upper semi-continuous of function with
values in a totally ordered set (a family in the sense of Definition \ref{familyequiv}),
\begin{equation}
\label{hache} h=h_{\mathcal{B}}:\Sing(J,b)\to (T, \leq).
\end{equation}
Note that $\Max(h_{\mathcal{B}})\subset \Sing(J,b)$.

The function \( h_{\mathcal{B}} \) is such that there will be a
$d$-dimensional basic object $(W,(A,c),E)$ associated to
$\Max(h)$. In this case we will say that $h$ is {\em associated to
}  $\equiv (W,(A,c),E)$, and it has the following properties:
\begin{enumerate}
\item[{\bf P1.}] There is an  inclusion of basic objects
as in (\ref{inclusions})
\begin{equation*}
(W,(A,c),E)\subset (W,(J,b),E).
\end{equation*}
Hence, any sequence of transformations of basic objects,
\begin{equation}
\label{alpha} (W_k,(A_k,c),E_k)\longleftarrow \ldots\longleftarrow (W_N,(A_N,c),E_N),
\end{equation}
induces a sequence of transformations with the same centers
\begin{equation}
\label{beta} (W_k,(J_k,b),E_k)\longleftarrow \ldots\longleftarrow (W_N,(J_N,b),E_N)
\end{equation}
and \( \Sing(A_s,c) \subset \Sing(J_s,b) \) for \( k \leq s\leq N \).

\item[{\bf P2.}] If the basic object \((W_k,(J_k,b),E_k)\)
is within Case $i$, then \((W_k,(A_k,c),E_k)\) will be within Case \(i-1\) (for \(i=2,3\)), which is
simpler.

\item[{\bf P3.}] If sequence (\ref{beta}) follows from sequence
(\ref{alpha}) as in  (P1), then we have that:
\begin{enumerate}
\item[{\bf (a)}] If \(\Sing(A_N,c) \neq \emptyset\), then \( \max h_k =\max
h_{k+1}=...= \max h_{N-1}= \max h_N \), and
\begin{equation*}
\Sing(A_j,c)=\Max(h_j) \mbox{ for } \mbox{ \ } k \leq j \leq N.
\end{equation*}

\item[{\bf (b)}] If \(\Sing(A_N,c)= \emptyset\), then
\( \max h_k =\max h_{k+1}=\ldots= \max h_{N-1}> \max h_N \), and
\begin{equation*}
\Sing(A_j,c)= \Max (h_j) \mbox{ for} \mbox { \ }  k \leq j \leq N-1.
\end{equation*}
\end{enumerate}

\item[{\bf P4.}] If \(\Theta:  (W_k,(J_k,b),E_k) \cong
(W'_k,(J'_k,b),E'_k)\) is an isomorphism of basic objects and if \(h_k\), \(h'_k\) are the corresponding
upper semi-continuous functions
\begin{equation*}
h_k:\Sing(J_k,b)\longrightarrow (T,\geq) \qquad\text{and}\qquad h_k': \Sing(J_k',b)\longrightarrow
(T,\geq),
\end{equation*}
then $h_k'(\theta(\xi))=h_k(\xi)$ for all $\xi \in \Sing(J_k,b).$
\item[{\bf P5.}] Both the functions \(h_i\), and the basic objects
\((W_i,(A_i,c),E_i)\) attached to it, are compatible with open restrictions of \((W_i,(J_i,b),E_i)\) (see
Definition \ref{DefRestricBO}).
\end{enumerate}
\end{parrafo}

\begin{remark}\label{equifunchi}
Note that:
\begin{enumerate}
\item[{\bf A)}] A basic object provides a way of describing closed
sets. We may think of \((W_i,(A_i,c),E_i)\) as a basic object
attached to the value \(\max h_i\); and property  (P3) expresses
this fact properly. It indicates, in particular, that if the
sequence (\ref{alpha}) is a {\em resolution} of $(W,(A,c),E)$,
then in  the sequence (\ref{beta}) we have that
\begin{equation*}
\max h_k =\max h_{k+1}=...= \max h_{N-1}>\max h_N,
\end{equation*}
and \( \Sing(A_j,c)= \Max (h_j) \) for \( k \leq j \leq N-1\). In other words, lowering \(\max h_k\) is
guaranteed by a resolution of the basic object \((W_k,(A_k,c),E_k)\), attached to this maximum value.

\item[{\bf B)}] If a sequence of transformations is defined (as in
\ref{transformation} and \ref{projections}),
\begin{equation}
\label{equiseq} {\mathcal B}=(W,(J,b),E) \longleftarrow \ldots \longleftarrow {\mathcal B}_k=(W_k,
(J_k,b),E_k),
\end{equation}
if
\begin{equation*}
\Theta: {\mathcal B}=(W,(J,b),E) \to {\mathcal B}^{\prime}=
(W^{\prime},(J^{\prime},b^{\prime}),E^{\prime})
\end{equation*}
is  an isomorphism of basic objects, and if
\begin{equation*}
{\mathcal B}^{\prime}=(W^{\prime},(J^{\prime},b^{\prime}),E^{\prime}) \longleftarrow \ldots
\longleftarrow {\mathcal B}_k^{\prime}=(W^{\prime}_k , (J^{\prime}_k,b^{\prime}),E^{\prime}_k)
\end{equation*}
is the corresponding induced sequence,  so that there are isomorphisms $$\Theta_j : (W_j,(J_j,b),E_j)
\cong (W_j^{\prime},(J_j^{\prime},b^{\prime}),E_j^{\prime}),$$ for \( 0 \leq j \leq k\) as in Definition
\ref{isoBO} (ii), then (P3) and  (P4) assert that \(\Theta_k : (W_k,(J_k,b),E_k) \cong
(W_k^{\prime},(J_k^{\prime},b^{\prime}),E_k^{\prime})\) induces an isomorphism \((W_k,(A_k,c),E_k)\cong
(W_k^{\prime},(A_k^{\prime},c^{\prime}),E_k^{\prime})\).

\item[{\bf C)}] \emph{The analog for general basic objects:}
If \({\mathcal B}=({\mathcal F},(W,E))\) is now a general basic object and $\{U_{\alpha}\}_{\alpha\in I}$
is an open covering of $W$ as in Definition \ref{defbo} (i), then for each index \(\alpha\) and each
basic object $(\widetilde{W}_{\alpha}, (B_{\alpha},d_{\alpha}),\widetilde{E}_{\alpha})$ there is an upper
semi-continuous function
\begin{equation*}
h_{\alpha}: \Sing(B_{\alpha},d_{\alpha})\longrightarrow(T,\leq).
\end{equation*}
If the collection of functions $\{h_{\alpha}\}$ patch so as to define a function,
\begin{equation*}
        h : F \to (T,\leq)
\end{equation*}
then  the  argument exhibit in {\bf(A)} and {\bf(B)}
   show that  the different basic objects
$$(\widetilde{W}_{\alpha},(A_{\alpha},c_{\alpha}),\widetilde{E}_{\alpha})$$ (defined in terms of
\(\alpha\)), define a general basic object attached to the value \(\max h_i\).

\item[{\bf D)}] In our approach, the functions
\(h_i\) will be defined entirely in terms of the functions \(\ord\) and \(n\). So
Proposition~\ref{wordngbo} guarantees that they define functions on general basic objects as required in
\textbf{(C)}.
\end{enumerate}
\end{remark}
\medskip

\begin{center}
{\bf Case 2.}
\end{center}
\medskip

Now we consider the case of a basic object \( (W,(J,b),E)\) such that the function \(\ord_{\mathcal B}
\cong 1\), but where $E$ may not be empty.

\begin{proposition}\label{case2bo}
Theorem \ref{Theoremd} holds for d-dimensional basic objects which are within Case 2.
\end{proposition}

The proof of this result, which will be stated in
\ref{proofofcase2}, is based on the construction of suitable
functions as indicated in \ref{pattern}. To illustrate the
philosophy behind this case we present the following example:

\begin{example}
\label{examplec2} Let \( W=\mathbb{A}_{\mathbf{k}}^{2} \). Set \( J=\left\langle x-y\right\rangle \) the
ideal of the diagonal $D$, and set \( E=\{H_{1},H_{2}\} \), where \( H_{1} \) is the \( x \)-axis and \(
H_{2} \) the \( y \)-axis. The resolution of the basic object \( (W,(J,1),E) \) consists of two
transformations:
\begin{equation*}
(W,(J,1),E)\longleftarrow(W_{1},(J_{1},1),E_{1}) \longleftarrow(W_{2},(J_{2},1),E_{2})
\end{equation*}
The first transformation is the blowing-up at the origin and the second is the blowing-up at \( D' \),
where \( D' \) is the strict transform of the diagonal.

Now consider the automorphism \( \Theta \) defined by the symmetry
about  the diagonal. The resolution is equivariant and the
automorphism \( \Theta \) lifts to automorphisms \( \Theta_{1} \)
in \( (W_{1},E_{1}) \) and \( \Theta_{2} \) in \( (W_{2},E_{2})
\).
\end{example}

\begin{parrafo}
\label{nbo} Let  ${\mathcal B}=(W,(J,b),E)$ be a basic object such that $\ord_{{\mathcal B}} \cong 1$.
Remark
    \ref{sobrecaso1} asserts that for any sequence of transformations
of basic objects
\begin{equation}
\label{alpha2} {\mathcal B}=(W,(J,b),E) \longleftarrow {\mathcal B}_1=(W_1 ,
     (J_1,b),E_1) \longleftarrow
\ldots\longleftarrow {\mathcal B}_k=(W_k ,
     (J_k,b),E_k),
\end{equation}
then $\Sing(J_k, b)=\emptyset$ or $\ord_{{\mathcal B}_i} \cong 1$
for \(1 \leq i \leq k\). Now consider the disjoint union
\begin{equation}
\label{number} E_k =E^+_k \cup E^-_k
\end{equation}
where \(E^-_k\) denotes the strict transform of hypersurfaces of
     \(E\). Note that sequence (\ref{alpha2})
     induces a sequence of
transformations
\begin{equation}
\label{alpha22} (W,(J,b),\emptyset) \longleftarrow (W_1 ,
    (J_1,b),E_1^+)\longleftarrow \ldots \longleftarrow (W_k ,
     (J_k,b),E^+_k).
\end{equation}
\end{parrafo}

\begin{definition}\label{ngbo}
With the notation introduced above, we define the upper semi-continuous function
\begin{equation*}
n_k: \Sing(J_k,b) \to {\mathbb Z}
\end{equation*}
\begin{equation*}
n_k(x)= \sharp\{ H_i \in E^-_k \ : \  x \in H_i \}.
\end{equation*}

If \(E^-_k=\{H_1,H_2,\ldots,H_s \}\), and \(\max n_k=m\), we set
\begin{equation} \label{EqDefIdealP}
       \mathcal{P}_{k}=\mathcal{P}_{k}(m)=\prod_{i_{1}<\cdots<i_{m}}
       \left(
       \sum_{j=1}^{m}\mathcal{I}(H_{i_{j}})     \right),
\end{equation}
and define the basic objects
\begin{equation}\label{funchic2}
        (W_k,(A_k,c),E_k^+)=
      \left\{
       \begin{array}{lcl}
       (W_k,(J_k,b),E_k^+) \cap (W_k,(\mathcal{P}_{k},1),E_k^+) &
        \text{if} & m>0  \\
(W_k,(J_k,b),E_k^+) & \text{if} & m=0,
         \end{array}
\right.
\end{equation}
and
\begin{equation}\label{funchic2000oh}
(W_k,(A_k,c),E_k)= \left\{
\begin{array}{lcl}
(W_k,(J_k,b),E_k) \cap (W_k,(\mathcal{P}_{k},1),E_k) &
\text{if} & m>0  \\
(W_k,(J_k,b),E_k) & \text{if} & m=0.
\end{array}
\right.
\end{equation}
\end{definition}

\begin{remark}\label{nngbo}
Note that:
\begin{enumerate}
\item[(a)] By definition
\begin{equation} \label{inclucaso2}
        (W_k,(A_k,c),E_k^+)
         \subset (W_k,(J_k,b),E_k^+).
\end{equation}

\item[(b)] The function \(n_k \)is  associated to $(W_k,(A_k,c),E_k)$ (linked by
the five properties stated in \ref{pattern}). These properties will be discussed below.

\item[(c)]For any \( \xi\in W_{k} \),
\begin{equation*}
        (\mathcal{P}_{k})_{\xi}=\calo_{W_k} \quad \Longleftrightarrow \quad
         \xi\not\in\Max{n_{k}}.
\end{equation*}
In fact if \( \xi\in\Max{n_{k}} \) then
\begin{equation*}
        (\mathcal{P}_{k})_{\xi}=\sum_{j=1}^{m={\mbox{\tiny max}(n_k)}}
        \mathcal{I}(H_{i_{j}})_{\xi}
        \qquad\text{ and }\qquad
        \xi\in\bigcap_{j=1}^{m}H_{i_{j}}\qquad H_{i_{j}}\in E^-_k .
\end{equation*}
  So $\Max n_k=\Sing(A_k,c)$. Note now that the expressions in
(\ref{funchic2000oh}) and (\ref{funchic2}), as intersection of basic objects, are compatible with
transformations in the sense of \ref{intersection}. This fact together with (a) prove that Property (P1)
stated  in \ref{pattern} holds for the functions $n_k$.

\item[(d)] If \( C \) is a permissible center of \(
(W_{k},(\mathcal{P}_{k},1),E_{k}) \), for any hypersurface \( H\in E_{k}^{-} \), and any irreducible
component \( C' \) of \( C \), either \( C'\subset H \) or \( C'\cap H=\emptyset \). This last
observation and the following lemma will indicate that
\begin{equation}\label{congcase2}
(W_k,(A_k,c),E_k)\cong (W_k,(A_k,c),E_k^{+}),
\end{equation}
in the sense of \ref{isoBO}, despite the fact that \( E_{k}\neq
E_{k}^{+} \). In fact the following lemma proves that both  basic
objects  have the same sequences of transformations.
\end{enumerate}
\end{remark}
\begin{lemma}
\label{case22} Let \((W,E)\) be a pair, assume that \(E\) is a
disjoint union  \( E=E^+ \cup E^-\). If \(C \subset W\) is a
permissible center for \((W,E^+)\), and each irreducible component
of \(C\) is either disjoint from or included in hypersurfaces of
\(E^-\), then \(C\) is permissible for \((W,E)\).
\end{lemma}

\begin{proof}
This follows from the definition of normal crossing (see Definition \ref{normalc}).
\end{proof}

\begin{parrafo}
\label{proofofcase2} {\em Proof of Proposition \ref{case2bo}:} First note that given  a sequence of
transformation of basic objects
\begin{equation}\label{alpha2oh}
{\mathcal B}=(W_0,(J_0,b),E_0) \longleftarrow \ldots \longleftarrow {\mathcal B}_k=(W_k,
    (J_k,b),E_k),
    \end{equation}
and isomorphisms $$\Theta: {\mathcal B}=(W,(J,b),E) \to {\mathcal B}^{\prime}=
(W^{\prime},(J^{\prime},b^{\prime}),E^{\prime}),$$ sequence (\ref{alpha2oh}) induces a sequence
\begin{equation}\label{ecuiniso1}
{\mathcal B}^{\prime}=(W^{\prime},(J^{\prime},b^{\prime}),E^{\prime}) \longleftarrow \ldots
\longleftarrow {\mathcal B}_k^{\prime}=(W^{\prime}_k ,
     (J^{\prime}_k,b^{\prime}),E^{\prime}_k).
     \end{equation}
In particular, for $0 \leq j \leq k$, we have the  isomorphisms
\begin{equation*}
       \Theta_j : (W_j,E_j) \cong (W_j^{\prime},E_j^{\prime}),
\end{equation*}
as in Definition \ref{isoBO} (b). Now set
\begin{equation}\label{factorization}
E_k=E_k^+\cup E_k^- \qquad\text{and}\qquad E_k^{\prime}=E_k^{\prime+}\cup E_k^{\prime-},
\end{equation}
and consider the  functions $n_k^{\prime}: \Sing(J_k^{\prime},b)\to {\mathbb Z}$,  and the corresponding
basic objects $$(W_k^{\prime},(A_k^{\prime},c^{\prime}),E_k^{\prime}).$$

Clearly the functions \( n_{k} \) are compatible with the
expressions in (\ref{factorization}) defined in terms of \(E_k\)
and \(E_k^{\prime}\). Hence, $$n_k(x)=n_k^{\prime}(\theta_k(x)).$$
Note that \(\theta_k\) defines an isomorphism $$(W_k,(A_k,c),E_k)
\cong (W_k^{\prime},(A_k^{\prime},c^{\prime}),E_k^{\prime})$$ as
 can be checked directly from the construction of these basic
object.
\medskip

We are now ready to give a  proof of Proposition \ref{case2bo}. We
begin by making some inductive assumptions on sequence
(\ref{alpha2oh}): Assume first that sequence (\ref{alpha2oh}) is
equivariant, and assume also that
\begin{equation*}
        \max n_0 \geq \max n_1 \geq\ldots\geq \max n_k.
\end{equation*}
Let \(k_0\) be the index such that \( \max n_{k_0-1}> \max n_{k_0}=\ldots= \max n_k \), and set
$(W_{k_0},(A_{k_0},c),E_{k_0})$ as in (\ref{funchic2000oh}).

By (\ref{congcase2})  a resolution of the basic object \((W_{k_0},(A_{k_0},c),E_{k_0}) \) is equivalent
to a resolution of \((W_{k_0},(A_{k_0},c), E_{k_0}^{+} ) \). We will show that
\((W_{k_0},(A_{k_0},c),E_{k_0}^{+} )\) is within Case 1, in the extended form indicated in Remark
\ref{caso1prima} (see  \ref{pattern} (P2)). To check this, note that \((W_0,(J_0,b),\emptyset)\) in
(\ref{alpha22})  is within Case 1, so the assertion follows from Remark \ref{agregadotrento1}.

\medskip

Assume now, by induction on \(k\), that the last \(k-k_0\) transformations in sequence (\ref{alpha2oh})
are defined by the first \(k-k_0\) steps of the equivariant resolution of
\((W_{k_0},(A_{k_0},c),E^+_{k_0}) \). We finally enlarge sequence (\ref{alpha2oh}) by the sequence
induced by the equivariant resolution of \((W_{k_0},(A_{k_0},c),E^+_{k_0}) \), as was indicated in
\ref{pattern} (P1). This defines for some index \(k^{\prime}>k\) an equivariant enlargement of sequence
(\ref{alpha2oh}). Note that if \(\max n_k=0\) this extension is already an equivariant resolution of
\((W,(J,b),E)\). If not
\begin{equation*}
\max n_{k_0}=...=\max n_k=...=\max n_{k^{\prime}-1}>\max n_{k^{\prime}}.
\end{equation*}
Repeating this procedure, ultimately \(\max n_{k^{\prime}}=0\), and hence, we have constructed an
equivariant resolution within this case, say
\begin{equation*}
R_{\mathcal B}=(W,(J,b),E) \longleftarrow (W_1 ,
     (J_1,b),E_1) \longleftarrow
\ldots\longleftarrow (W_N ,
     (J_N,b),E_N)
\end{equation*}
together with  a sequence \(0=k_0<k_1<k_2....<k_s \leq N\), and
\begin{equation*}
     \max n_{k_{i-1}}>  \max n_{k_{i}}= \max n_{k_{i}+1}=\ldots=\max
n_{k_{i+1}-1}>\max n_{k_{i+1}}\ldots.
\end{equation*}
This proves Proposition \ref{case2bo}. \qed
\end{parrafo}

Note that in the Example~\ref{examplec11},  \( \max{n_{0}}=2>\max{n_{1}}=1>\max{n_{2}}=0 \).
\medskip

\begin{center}
{\bf  Case 3: Proof of Theorem \ref{Theoremd} for basic objects.}
\end{center}
\medskip

Our next objective is the following statement:

\begin{proposition}\label{case3bo}
Theorem \ref{Theoremd} holds for $d$-dimensional basic objects.
\end{proposition}

We will start with an example to illustrate the invariants involved in the general case, and then we will
introduce some auxiliary definitions and results that will be used in the proof of the proposition, which
will be detailed in \ref{olvidada}.

\begin{example}\label{examplecase3bo}
Let \(W\) be the affine plane \({\mathbb A}^2_k\), let
     \(J=<x^3,y^4> \subset k[x,y]\),  set
\(b=2\), and consider the basic object \({\mathcal B}=(W,(J,2),E=\emptyset)\). In this case
\(\Sing(J,2)\) is the origin in \({\mathbb A}^2_k\), and $\ord_{{\mathcal B}}=3/2$. Now blow-up at the
origin, \((W,(J,2), \emptyset ) \longleftarrow (W_1,(J_1,2), E_1=\{H_1\}).\)

Note that \(J_1 \subset {\mathcal O}_{W_1}\) is defined so that \(J {\mathcal
O}_{W_1}=\mathcal{I}(H_{1})^2 J_1\). However in this case, since the function \( \ord \) is not equal to
1, there is a new factorization of the form
\begin{equation*}
J_1=\cali(H_1)\overline{J}_1.
\end{equation*}
     This
factorization of \(J_1\) does not hold for a basic object for which \(\ord \cong 1\), in fact \(J_1=
\overline{J}_1\) if \(\ord \cong 1\).
\end{example}

\begin{parrafo}\label{weakagain}
Given a sequence of transformation of basic objects
\begin{equation}
\label{alpha3} (W_0,(J_0,b),E_0) \longleftarrow \ldots \longleftarrow (W_k ,
     (J_k,b),E_k),
\end{equation}
consider the  expression
\begin{equation}
\label{expresionj} J_k=\cali(H_1)^{\alpha (1)} \cali(H_2)^{\alpha (2)}\ldots \cali(H_k)^{\alpha
(k)}\cdot\overline{J}_k
\end{equation}
which is unique, if we require that \( \overline{J}_k\) does not vanish along (any component of) any
\(H_i\), arising from a previous transformation.
\end{parrafo}

\begin{remark} \label{RemMejorbk}
Let \( b_{k}=\max\nu_{\bar{J}_{k}} \), where $\nu_{\bar{J}}$  is the numerator of the function introduced
in \ref{nnn}. Note that \( \Sing(\bar{J}_{k},b_{k})=\Max{\nu_{\bar{J}_{k}}} \). If the transformation is
such that \( C_{k}\subset\Sing(\bar{J}_{k},b_{k}) \) then \( (W_{k+1},(\bar{J}_{k+1},b_{k}),E_{k+1}) \)
is the transform of the basic object \( (W_{k},(\bar{J}_{k},b_{k}),E_{k}) \). Moreover if \(
b_{k+1}=\max{\nu_{\bar{J}_{k+1}}} \) then \( b_{k+1}\leq b_{k} \) (see Remark \ref{deladesig}).
\end{remark}

\begin{definition}\label{word}
\label{defworder} Given a sequence of transformations of basic objects as in (\ref{alpha3}), and
expressions as in (\ref{expresionj}), we define the  upper semi-continuous function:
\begin{equation*}
\begin{array}{rccl} \word_k: & \Sing(J_k,b) & \to & {\mathbb Q
}\\
     & x & \to & \word_k(x)= \nu_x( \overline{J}_k)/b.
\end{array}
\end{equation*}
\end{definition}

\begin{remark}
\label{remarkc3} Given a sequence of transformations as in (\ref{alpha3}) note that:
\begin{enumerate}
\item[(a)] If sequence (\ref{alpha3}) is defined by blowing up centers
\(C_i \subset \Max\word_i \), then by Remark \ref{RemMejorbk}
\begin{equation*}
     \max \word_0 \geq \max
\word_1 \geq \cdots\geq \max \word_k.
\end{equation*}

\item[(b)] If \(\max \word_k=0\) then
\begin{equation*}
     J_k=\mathcal{I}(H_1)^{\alpha
(1)} \mathcal{I}(H_2)^{\alpha (2)}\cdot\ldots\cdot \mathcal{I}(H_k)^{\alpha (k)}.
\end{equation*}

\item[(c)] If \(\max \word_k > 0\) then
the functions \( \word_k \) will play the role of the functions \(h_k\) in  (\ref{hache}). In fact if
\(\max \word_k =b_{k}/b >0\) for \( 0<b_{k} \) then we define:
\begin{equation}\label{funchic3}
(W_k,(A_k,c),E_k)=(W_k,(\overline{J}_k,b_{k}),E_k) \cap (W_k,(J_k,b),E_k).
\end{equation}

\item[(d)] The basic object
\((W_k,(\overline{J}_k,b_k),E_k) \) is such that \(\ord\cong 1\), so the same holds for the basic objects
\( (W_k,(A_k,c),E_k) \) (\ref{sobrecaso1} (2)), and \(\Sing(\overline{J}_k,b_k)=\Max \word_k\). By (a)
\(b_i/b=\max \word_i \geq b_{i+1}/b= \max \word_{i+1}\). Note that if equality holds, then
$$(W_{i+1},(\overline{J}_{i+1},b_{i+1}),E_{i+1})$$ is the transform of
\((W_{i},(\overline{J}_{i},b_{i}),E_{i}) \).

\item[(e)] The function \(\word_k\) is associated to $(W_k,(A_k,c),E_k)$,
and hence it satisfies    the five properties stated in
\ref{pattern}.

\end{enumerate}
\end{remark}

\begin{lemma}\label{casomono}
Let $\Theta:{\mathcal B}=(W_0,(J_0,b),E_0)\cong {\mathcal
B}^{\prime}=(W^{\prime}_0,(J^{\prime}_0,b^{\prime}),E^{\prime}_0)$
be an isomorphism of  basic objects. Consider a sequence of
transformations over ${\mathcal B}$,
\begin{equation}
(W_0,(J_0,b),E_0) \longleftarrow \ldots \longleftarrow (W_k ,
     (J_k,b),E_k),
\end{equation}
together with the induced sequence over ${\mathcal B}^{\prime}$ (as in Definition \ref{isoBO} (ii)),
\begin{equation}
{\mathcal B}^{\prime}=(W^{\prime},(J^{\prime},b^{\prime}),E^{\prime}) \longleftarrow \ldots
\longleftarrow {\mathcal B}_k^{\prime}=(W^{\prime}_k ,
     (J^{\prime}_k,b^{\prime}),E^{\prime}_k),
     \end{equation}
and let
\begin{equation*}
     J_k=\mathcal{I}(H_1)^{\alpha
(1)} \mathcal{I}(H_2)^{\alpha (2)}\cdot\ldots\cdot \mathcal{I}(H_k)^{\alpha (k)}\cdot \overline{J}_k,
\end{equation*}
and
\begin{equation*}
     J_k^{\prime}=\mathcal{I}(H_1^{\prime})^{\alpha^{\prime}
(1)} \mathcal{I}(H_2^{\prime})^{\alpha^{\prime} (2)}\cdot\ldots\cdot
\mathcal{I}(H_k^{\prime})^{\alpha^{\prime} (k)}\cdot \overline{J}_k^{\prime},
\end{equation*}
be the factorizations of the ideals $J_k$ and $J_k^{\prime}$ as in
\ref{weakagain}. Set $\xi_{k}\in
Sing(\overline{J}_k,b_{k})=\Max\word_{k}$ and
$\theta(\xi_{k})=\xi_{k}^{\prime}$. Then:
\begin{enumerate}
\item[(i)] For each \(\xi_k \in \Sing(J_k,b)\) the rational numbers
\(\word_k(\xi_{k})\), and the $\{\alpha(i) /b\}_{\{i=1,\ldots,k\}}$ corresponding to hypersurfaces
containing $\xi_k$, can be expressed in terms of
    the rational numbers
$$\{\ord_{\mathcal B_0} (\xi_0),\ord_{\mathcal B_1}(\xi_1),\ldots,\ord_{\mathcal B_k}(\xi_k)\},$$ where
$\xi_i$ denotes the image of $\xi_k$ in $W_i$.
\item[(ii)] For each $i=1,\ldots,k$,  \(\xi_{k}\in H_i \in E_k\) if and only
if \(\xi^{\prime}_{k}\in \theta_k(H_i) \in E_k^{\prime}\), and in
that case \(\alpha(i)/b=\alpha^{\prime}(i)/b^{\prime}\).
\item[(iii)] The function $\word_k$ is equivariant, i.e.,
\(\word_k(\xi_k)=\word_k^{\prime}(\xi_k^{\prime})\).
\end{enumerate}
\end{lemma}

\begin{proof}
We first indicate how to prove (i)  by induction on \(k\). If \(k=0\), \(\word_0=\ord_{\mathcal B_0}\) so
it follows from \ref{propertiesnord}.
\medskip

If \(k=1\), and \(\xi_1 \notin H_1 \), then \(\xi_1\) can be identified with \(\xi_0 \in W_0\), and
\(\word_1(\xi_1)=\word_0(\xi_0)=\ord_0(\xi_0)\). If \(k=1\) and \(\xi_1 \in H_1 \), note that (as in
Example~\ref{examplecase3bo}) \(\alpha_1\) is the order of \(J\) along the first center of
transformation. Hence \(\alpha_1/b=\ord_{{\mathcal B}_0}(\xi_0)\), and then \(\ord_{{\mathcal
B}_1}(\xi_1)=\alpha_1/b+ \word_1(\xi_1)\) and (i) holds in this case.
\medskip

Assume that (i)  holds for \(i=k-1\), and we want to prove the assertion for \(i=k\). Then we argue as
above distinguishing the case when  \(\xi_k \in  H_k\), or \(\xi_k\notin H_k\). A general expression of
the invariants mentioned in  (i)  in terms of \(\{ \ord_{\mathcal B_0}(\xi_0), \ord_{\mathcal
B_1}(\xi_1),\ldots, \ord_{\mathcal B_k}(\xi_k)\}\) can be found in \cite[Theorem 7.6]{EncVil97:Tirol}.
\medskip

Parts (ii) and (iii) follow from the proof of part (i).
\end{proof}

\begin{parrafo} \label{olvidada}
{\em Proof of Proposition \ref{case3bo}:} Let \((W,(J,b),E)\) be a basic object, and assume  by induction
on \(k\), that we have defined an equivariant sequence of transformations
\begin{equation}
\label{alpha33} (W_0,(J_0,b),E_0) \longleftarrow \ldots \longleftarrow (W_k ,
     (J_k,b),E_k),
\end{equation}
at centers \(C_i \subset \Max\word_i \) for $i=0,1,\ldots,k-1$. Hence, if \(\Theta: {\mathcal
B}=(W,(J,b),E) \to {\mathcal B}^{\prime}= (W^{\prime},(J^{\prime},b^{\prime}),E^{\prime})\) is an
isomorphism, then there is an induced sequence of transformations
\begin{equation}
\label{ecuiniso33} {\mathcal B}^{\prime}=(W^{\prime}_0,(J^{\prime}_0,b^{\prime}),E^{\prime}_0)
\longleftarrow \ldots \longleftarrow {\mathcal B}_k^{\prime}=(W^{\prime}_k ,
    (J^{\prime}_k,b^{\prime}),E^{\prime}_k)
     \end{equation}
at centers \(\Theta(C_i)=C_i^{\prime} \subset \Max\word_i \) for $i=0,1,\ldots,k-1$ (see Definition
\ref{isoBO}).
\medskip

Let  \(k_0\) be  the smallest index  in sequence (\ref{alpha33}) (or in (\ref{ecuiniso33})) such that
$$\max \word_{k_0}=\max\word_{k}.$$ We distinguish two cases:
\begin{itemize}
\item {Case \(\max \word_{k_0}>0\).}
\end{itemize}
Recall that in this case, if we set  $\max \word_{k_0}=b'/b$, the basic object
$$(W_{k_0},(A_{k_0},c),E_{k_0})=(W_{k_0},(\overline{J_{k_0}},b'),E_{k_0}) \cap
(W_{k_0},(J_{k_0},b),E_{k_0})$$ is within the previously treated Case 2; moreover, \(
\Sing(A_{k_{0}},c)=\Max\word_{k_{0}} \), and the basic object $(W_{k_0},(A_{k_0},c),E_{k_0})$ fulfills
all conditions stated in  \ref{pattern}. As indicated in Remark~\ref{equifunchi} (A), a resolution of
$(W_{k_0},(A_{k_0},c),E_{k_0})$, say
\begin{equation}
\label{alpha33oh} (W_{k_0},(A_{k_0},c),E_{k_0})\longleftarrow \ldots\longleftarrow (W_N,(A_N,c),E_N),
\end{equation}
induces a sequence
\begin{equation}
\label{beta33} (W_{k_0},(J_{k_0},b),E_{k_0})\longleftarrow \ldots\longleftarrow (W_N,(J_N,b),E_N),
\end{equation}
such that $$\max \word_{k_0}=\max \word_{{k_0}+1}=..=\max \word_{N-1}>\max \word_N.$$

Assume that the last \(k-k_0\) steps of sequence (\ref{alpha33}) are the first \(k-k_0\) steps of the
sequence (\ref{beta33}) (induced by the resolution (\ref{alpha33oh})), and now extend sequence
(\ref{alpha33}), as above, to a sequence of length \(N\). Note that \(N>k\) and that all the hypothesis
on sequence (\ref{alpha33}) also hold on this enlarged sequence. In particular, if sequence
(\ref{ecuiniso33}) is defined in terms of sequence (\ref{alpha33}), then by Lemma \ref{casomono} and the
equivariance of Case 2,  sequence (\ref{ecuiniso33}) can also be extended to a sequence of length \(N\),
such that \(\word_{k_0} > \max\word_N\), and the  isomorphisms \(\theta_j\), \(j=0,1,\ldots,k\) can be
lifted  for $j=k+1,\ldots,N$.
\medskip

If \(\max \word_N>0\) we repeat the argument. But this can happen only finitely many times. In fact, if
     \(\max\word_k =b'/b\),
then this positive rational number can drop at most \(b^{\prime}\) times. This leads us ultimately to the
case \(\max \word =0\).
\medskip

\begin{itemize}
\item {Case \(\max \word_{k_0} =0\).}
\end{itemize}
In this case we have that
\begin{equation}
\label{remarkc33b}
    J_{k_0}=\mathcal{I}(H_1)^{\alpha
(1)} \mathcal{I}(H_2)^{\alpha (2)}\cdot\ldots\cdot \mathcal{I}(H_{k_0})^{\alpha ({k_0})},
\end{equation}
and hence also that
\begin{equation}
\label{remarkc33bprime}
    J_{k_0}^{\prime}=\mathcal{I}(H_1^{\prime})^{\alpha
(1)} \mathcal{I}(H_2^{\prime})^{\alpha (2)}\cdot\ldots\cdot \mathcal{I}(H_{k_0}^{\prime})^{\alpha
({k_0})}.
\end{equation}
This is what we call the {\em monomial case} which will be treated in Section \ref{SeccMonom}. We will
only mention here that in this case  it is very simple to extend sequence (\ref{alpha33}) (and hence
sequence (\ref{ecuiniso33})) to a resolution. This can be done in various ways, and simply by considering
the exponents \(\alpha(i)\): One can define a  new upper semi-continuous equivariant function whose
value at \( \xi \in \Sing (J_k,b)\) depends entirely on  \( \{ \alpha(i)/b : \xi \in H_i \in E_k\}\).
These functions are equivariant by Lemma \ref{casomono} (ii), and reach their maximum value at a smooth
permissible center. An extension of the first \(k_0\) steps of the equivariant sequence (\ref{alpha33})
(of (\ref{ecuiniso33})) to an equivariant resolution  is achieved by repeating this procedure (see
Section \ref{SeccMonom}).
     \qed
\end{parrafo}




\section{Proof of Theorem \ref{Theoremd} for  general basic objects.}
\label{proofofgbo} We now address the proof of Theorem \ref{Theoremd} for general basic objects
\({\mathcal B}=({\mathcal F},(W,E))\) of dimension \(d\).  Fix an open covering \(\{U_{\alpha}\}_{\alpha
\in \Lambda}\), a collection of smooth \(d\)-dimensional subschemes \(\widetilde{W}_{\alpha} \subset
U_{\alpha}\) and \(\mathcal{B_{\alpha}}=(\widetilde{W}_{\alpha}, (B_{\alpha},d_{\alpha})
\widetilde{E}_{\alpha})\) with the properties stated in Definition \ref{defbo}. Recall that by
Proposition \ref{wordngbo}, there are well defined equivariant functions, $$\ord^d: F\to {\mathbb Q}
\mbox{ }\mbox{ } \mbox{ and } n^d: F\to {\mathbb Z},$$ obtained by patching the functions $\ord_{
\mathcal{B_{\alpha}}}^d$ and $n_{\mathcal{B_{\alpha}}}^d$.

\

As in the case of basic objects we will  consider three cases:
\begin{itemize}
\item Case 1:  $\ord_{\calB}=1$ and $E=\emptyset$.
\item Case 2:  $\ord_{\calB}=1$ and $E\neq\emptyset$.
\item Case 3. The general case: Proof of Theorem~\ref{Theoremd}.
\end{itemize}
Again, we stress here that Cases 1 and 2 are not special cases,
and that they appear as intermediate steps when proving the
general case.

\

Recall that Theorem~\ref{Theoremd} holds for \( d=0 \), so we will argue by induction on $d$.

\

\begin{itemize}
\item {\bf Case 1.}
\end{itemize}
\medskip
Since \({\mathcal{F}}\) is of dimension \(d\), the closed set
\(F(\subset W)\), has  dimension at most \(d-1\). Lemma
\ref{case12} asserts that the \((d-1)\)-dimensional components of
\(F\) are both open and closed in \(F\), defining a smooth
permissible center; and furthermore, by blowing up such center the
transform has a structure of \((d-1)\)-dimensional general basic
object. A similar argument as the one used in the proof of
Proposition \ref{case1} shows that this particular transformation
is equivariant. This proves the theorem for \(d\)-dimensional
general basic objects within Case 1.
\begin{remark}\label{Caso1Prima}
We may also repeat the outcome of Remark \ref{caso1prima} in the context of general basic objects, namely
that the construction of the resolution indicated above, also applies for transform of a general basic
object within Case 1.
\end{remark}
\medskip

\begin{itemize}
\item {\bf Case 2.}
\end{itemize}
\medskip

Let \({\mathcal F}_{0}\) be a general basic object
    within Case 2. As  in Proposition
\ref{wordngbo}, set  $$n_{{\mathcal F}_0} : F_0 \to {\mathbb N}.$$ Given a sequence of transformation of
general basic objects,
    \begin{equation}
\label{aqui}
        \begin{array}{ccccccc} ({{\mathcal F}}_0,(W_0,E_0)) &
            \longleftarrow & \ldots & \longleftarrow &
            (\mathcal{F}_{r-1},(W_{r-1},E_{r-1})) & \longleftarrow &
            (\mathcal{F}_{r},(W_{r},E_{r})).\\
            \cup & & & & \cup & & \cup \\
            F_0 & & & & F_{r-1} & & F_r
        \end{array}
\end{equation}
we can consider the partition, as  in (\ref{number}), $$E_k=E_k^+\cup E_k^-.$$ Finally define the
functions
\begin{equation*}
n_k: F_k \to {\mathbb Z}
\end{equation*}
\begin{equation*}
n_k(x)= \sharp\{ H_i \in E^-_k \ : \  x \in H_i \}.
\end{equation*}
  In this way we extend, to general basic objects, the previously
defined function on basic objects.
  The same holds for the  \( \mathcal{P}_{k}(m) \) in
(\ref{EqDefIdealP}) as we sketch below.
\medskip
For each $\alpha\in \Lambda$  (as in the covering introduced in \ref{defbo}), define basic objects
\((\widetilde{W}_{\alpha,k}, (A_{\alpha,k},c_{\alpha}), \widetilde{E}_{\alpha,k}^+)\) and
\((\widetilde{W}_{\alpha,k}, (A_{\alpha,k},c_{\alpha}), \widetilde{E}_{\alpha,k})\), with
\((\widetilde{W}_{\alpha,k},(A_{\alpha,k},c_{\alpha}), \widetilde{E}_{\alpha,k})
    \subset (\widetilde{W}_{\alpha,k},(B_{\alpha,k},d_{\alpha}),
\widetilde{E}_{\alpha,k})\) as in (\ref{funchic2}) and (\ref{funchic2000oh}):
\begin{equation*}
(\widetilde{W}_{\alpha,k}, (A_{\alpha,k},c_{\alpha}), \widetilde{E}_{\alpha,k}^+)=
(\widetilde{W}_{\alpha,k},(B_{\alpha,k},d_{\alpha}), \widetilde{E}_{\alpha,k}^+)
         \cap (\widetilde{W}_{\alpha,k},(\mathcal{P}_{\alpha,k},1),E_k^+),
\end{equation*}
and
\begin{equation*}
(\widetilde{W}_{\alpha,k}, (A_{\alpha,k},c) \widetilde{E}_{\alpha,k})=
(\widetilde{W}_{\alpha,k},(B_{\alpha,k},d) \widetilde{E}_{\alpha,k})
         \cap (\widetilde{W}_{\alpha,k},(\mathcal{P}_{\alpha,k},1),E_k).
\end{equation*}

Now note that the  \((\widetilde{W}_{\alpha,k}, (A_{\alpha,k},c) \widetilde{E}_{\alpha,k}) \) define a
new \(d\)-dimensional general basic objects (see Remark~\ref{equifunchi} (C)).  So for \(k=0\) we obtain
a new general basic object, say \(({\mathcal F^{\prime}_0},(W_0,E_0))\), and we repeat the argument used
in Case 2 in the context of basic objects: Since \(({\mathcal F^{\prime}_0},(W_0,E_0))\) is in Case 1
there is an equivariant resolution,
\begin{equation*}
       \begin{array}{ccccccc}
      ({{\mathcal F}}_0^{\prime},(W_0,E_0)) &
      \longleftarrow & \ldots & \longleftarrow &
      (\mathcal{F}_{r-1}^{\prime},(W_{r-1},E_{r-1}))
      &
      \longleftarrow &
      (\mathcal{F}_{r}^{\prime},(W_{r},E_{r}))\\
            \cup & & & & \cup & & \cup\\
      F_0^{\prime} & & & & F_{{r^{\prime}}-1}^{\prime} &
            & F_{r^{\prime}}^{\prime}=\emptyset
       \end{array}
\end{equation*}
which induces a sequence of transformations
\begin{equation*}
       \begin{array}{ccccccc} ({{\mathcal F}}_0,(W_0,E_0)) &
          \longleftarrow & \ldots & \longleftarrow &
            (\mathcal{F}_{r-1},(W_{r-1},E_{r-1})) & \longleftarrow &
            (\mathcal{F}_{r},(W_{r},E_{r}))\\
            \cup & & & & \cup & & \cup \\
            F_0 & & & & F_{r-1} & & F_r,
        \end{array}
\end{equation*}
such that \(\max n_0 = \max n_1 =...=\max n_{r-1}>\max n_r\), and \(F_i^{\prime}=\Max n_i \subset F_i\)
for \(i=0,1,...,r-1\) (see Remark \ref{equifunchi} (A)).
\medskip

If \(\max n_r \neq 0\), a new \(d\)-dimensional general basic object is attached  to the value \(\max
n_r\), as we did before. We repeat this argument up to the case in which \(\max n_r=0\). In this case
\((\mathcal{F}_{r},(W_{r},E_{r}))\) is within Case 1 (in the sense of Remark \ref{Caso1Prima}), and hence
can be extended to a resolution.

\

Assume now that there is an isomorphism of general basic objects, say \( {\mathcal F} \to {\mathcal
G}\).  If \( {\mathcal F}^{\prime}\) and \( {\mathcal G}^{\prime}\) are both defined, as before, in terms
of the value \(\max n\), then by Remark \ref{equifunchi}
  (B), \( {\mathcal F}^{\prime} \cong {\mathcal
G}^{\prime}\) (see also Remark \ref{equifunchi} (D)). Hence these resolutions are equivariant.
\medskip

\begin{itemize}
\item {\bf Case 3: Proof of Theorem \ref{Theoremd}.}
\end{itemize}
\medskip

We may argue as in the previous case, to show that the functions \(\word_i\) can also be defined for
general basic objects, and for a sequence of transformations of general basic objects.  The same
arguments used above, replacing here the functions \(n_i\) by \(\word_i\), show that:
\begin{enumerate}
\item[(i)] For any general basic object \( {\mathcal F}_0 \), we can associate
to it a sequence of transformations,
\begin{equation}
\label{sequencei}
        \begin{array}{ccccccc} ({{\mathcal F}}_0,(W_0,E_0)) &
            \longleftarrow & \ldots & \longleftarrow &
(\mathcal{F}_{N-1},(W_{N-1},E_{N-1})) & \longleftarrow &
           (\mathcal{F}_{N},(W_{N},E_{N}))\\
            \cup & & & & \cup & & \cup \\
            F_0 & & & & F_{N-1} & & F_N,
        \end{array}
\end{equation}
such that \(\max \word_N=0\).

\item[(ii)] Given an isomorphism \(\Theta: {\mathcal F}_0
\to {\mathcal F}_0^{\prime}\), then the sequences defined in (i) for \({\mathcal F}_0\) and \({\mathcal
F}_0^{\prime}\), are both of same length, say \(N\), and linked by isomorphisms \( \theta_i:{\mathcal
F}_i \to {\mathcal F}_i^{\prime}\) for \( 0\leq i \leq N\).
\end{enumerate}

Finally, the discussion in Case \(\word_N=0\) applies also for general basic objects. In fact the
extension of sequence (\ref{sequencei}) to a resolution is achieved by means of a function that depends
only on the invariants  of Lemma \ref{casomono}. Thus, Lemma \ref{casomono}
    together with the equivariance of the functions
\(\ord_{{\mathcal F}_i}\) assert that the extensions of sequence (\ref{sequencei}) to a resolution of \(
\mathcal{F}_{0} \) (of \( \mathcal{F}'_{0} \)) are defined together with linking isomorphisms. This
proves Theorem \ref{Theoremd}. \qed




\part{The Algorithm.}
\label{appendix}




\section{On the definition of  the functions $f_i^d$.}
\label{SeccConsFunc} We address here the explicit description of the functions \(f_i^d\) introduced in
Definition \ref{AlgResol}. These functions provide the elementary proofs of desingularization and
principalization given in \ref{existence} and \ref{proofprincipal}; and the resolution of basic objects
they define turn out being that defined in Theorem \ref{Theoremd}. We summarize some previous results for
self-containment, although some references to Part \ref{nature} will be needed.
\begin{parrafo}
\label{descloc} Let \((W_0,(J_0,b),E_0=\{H_1,\ldots,H_l\})\) be a \(d\)-dimensional basic object. Given a
sequence of transformations,
\begin{equation}
\label{extra11} (W_{0},(J_{0},b),E_{0})
\stackrel{\pi_1}{\longleftarrow}\cdots\stackrel{\pi_r}{\longleftarrow} (W_{r},(J_{r},b),E_{r})
\end{equation}
we define expressions
\begin{equation}
\label{extra22} J_{i}=\mathcal{I}(H_{l+1})^{a_{1}}\cdots \mathcal{I}(H_{l+i})^{a_{r}}\overline{J}_{i},
\end{equation}
where \(H_{l+i}\) denotes the exceptional divisor at that \(i\)-th blowing-up, and  upper semi-continuous
functions
\begin{equation*}
        \begin{array}{cccc}
      \word_{i}^{d}: & \Sing(J_{i},b) & \longrightarrow &
      \frac{1}{b}\mathbb{Z}\subset\mathbb{Q}\\ & \xi & \longrightarrow &
      \frac{\nu_{\overline{J}_{i}}(\xi)}{b}
        \end{array}
\end{equation*}
are defined in terms of (\ref{extra22}) (see (\ref{expresionj}) and
  Definition \ref{word}).
\medskip

Assume now that sequence (\ref{extra11}) is a sequence of transformations at permissible centers \(Y_i\)
such that
\begin{equation}
        \label{adcons}
        Y_{i}\subset\Max\word_{i}^d\subset\Sing(J_{i},b)
\end{equation}
for \(i=0,1,\ldots,r\). In this case:
\begin{equation} \label{IneqwordPt}
        \word_{i-1}^d(\pi_{i}(\xi_{i}))\geq\word_{i}^d(\xi_{i})
\end{equation}
for every \( \xi_{i}\in\Sing(J_{i},b) \), and  equality holds if \(\pi_{i}(\xi_{i})\notin Y_{i-1}\) (see
Remark \ref{deladesig}). In particular
\begin{equation}
        \label{Ineqword}
        \max\word_{0}^d\geq\cdots\geq\max\word_{r}^d.
\end{equation}

\

Pick \(k\in \{0,1,\ldots,r\}\).  If \( \max\word_{k}^d>0 \) let \( k_{0} \) be the smallest index so that
\begin{equation*}
\max\word_{k_{0}-1}^d>\max\word_{k_{0}}^d=\max\word_{k}^d,
\end{equation*}
(\( k_{0}=0 \) if \( \max\word_{0}^d=\cdots=\max\word_{r}^d \)). Write
\begin{equation*}
E_{k}=E_k^{+}\sqcup E_{k}^{-}
\end{equation*}
where \(E_{k}^{-} \) is the set hypersurfaces of \(E_{k} \) which are strict transforms of hypersurfaces
of \( E_{k_{0}} \).  Now define
\begin{equation*}
        n_{k}^d(\xi)=\left\{
        \begin{array}{lll}
      \#\{H\in E_{k}\mid \xi\in H\} & {\rm if} &
      \word_{k}^d(\xi)<\max\word_{k}^d \\
      \#\{H\in E_{k}^{-}\mid \xi\in H\} & {\rm if} &
      \word_{k}^d(\xi)=\max\word_{k}^d.
       \end{array}
        \right.
\end{equation*}
\end{parrafo}

\begin{definition}
\label{Deft1} If condition (\ref{adcons}) holds, and \( \max\word_{r}^d>0 \), we define, for the index
\(r\), a function \(t_r^d\) by setting:
\begin{equation*}
\begin{array}{ccccc}
     t_{r}^d: & \Sing(J_{r},b) & \longrightarrow &
        (\mathbb{Q}\times\mathbb{Z},\leq) &   \\
         & \xi & \longrightarrow & (\word_{r}^d(\xi),n_{r}^d(\xi)), &
\end{array}
\end{equation*}
where $\mathbb{Q}\times\mathbb{Z}$ is ordered lexicographically.
In the same way we define functions \(
t_{r-1}^d,t_{r-2}^d,\ldots,t_{0}^d \). If \(Y_r\subset \Max
t^d_r\) is a permissible center, then \(Y_r\) is said to be a {\em
\(t^d_r\)-permissible center}.
\end{definition}

\begin{parrafo}
\label{PropiedFt} {\bf Properties of the inductive function \( t^d_i \)} (cf.
\cite[4.15]{EncVil97:Tirol}).
\begin{enumerate}
       \item Set \(F^i=\Sing(J_i,b)\). Given \((q,m) \in
        \mathbb{Q}\times\mathbb{Z}\), set \(F^i_{(q,m)}=\{ x \in F^i:
        t^d_i(x) \geq (q,m)\}\).  Note that
        \begin{equation*}
      F_{(q,m)}=\{x \in F^i: \word_i(x)
      > q \} \cup \{x \in F^i: \word_i(x)\geq q, \mbox{ } n_i(x) \geq m \}.
        \end{equation*}
       Here the first term is closed since \(\word\) is upper semi-continuous
       and takes only finitely many different values. Since $n_i$ and
$\word_i$ are  upper semi-continuous functions, it easily follows that  the
        second term is
       also
        closed.  This shows that \(F_{(q,m)} \) is closed, and hence that
        each \( t_{i}^d \) is upper semi-continuous.  Note also that \(
        \Max{t_{i}^d}\subset\Max\word_{i} \).
        \item If the sequence of transformations
        \begin{equation}
      \label{extra01000}
      (W_{0},(J_{0},b),E_{0})
      \stackrel{\pi_1}{\longleftarrow}
     (W_{1},(J_{1},b),E_{1})
      \stackrel{\pi_2}{\longleftarrow}
      \cdots\stackrel{\pi_r}
     {\longleftarrow} (W_{r},(J_{r},b),E_{r})
       \end{equation}
       is \(t^d\)-permissible, namely if  \(
        Y_{i}\subset\Max{t_{i}^d} (\subset\Max\word_{i})\),
then for each index \( i=0,1,\ldots,r\),
\begin{equation*}
        t_{i-1}^d(\pi_{i}(\xi_{i}))\geq t_{i}^d(\xi_{i})
\end{equation*}
for all $\xi_{i}\in\Sing(J_{i},b)$, and equality holds if \( \pi_{i}(\xi_{i})\not\in Y_{i-1} \). In
particular
\begin{equation*}
        \max{t_{0}^d}\geq\cdots\geq\max{t_{r}^d}.
\end{equation*}
\item We say that \( \max{t}^d \) drops at \( i_{0} \) if \(
\max{t_{i_{0}-1}^d}>\max{t_{i_{0}}} \).  If \( \max\word_{0}^d=\dfrac{b'}{b} \) and \( \dim{W_{0}}=d \),
note that \( \max{t_{i}^d}=\left(\dfrac{s}{b},m\right) \), \( 0\leq s\leq b' \), \( 0\leq m\leq d \).  So
it is clear that \( \max{t^d} \) can drop at most \( b'd \) times.
\item The functions \(t^d_i\) are the {\em inductive invariants}.  In fact,
if \(R(1)(\Max t)\neq \emptyset\) (\ref{def1}), then this is our
{\em canonical choice of center}; and after repeatedly blowing up
 this center we may assume that \(R(1)(\Max t)= \emptyset\).  On
the other hand, if \(R(1)(\Max t)= \emptyset\), then via some form
of induction which will be described below, it is possible to
construct a unique enlargement of sequence  (\ref{extra01000}),
\begin{multline}
        \label{nueva2}
        (W_{0},(J_{0},b),E_{0})\longleftarrow\cdots\longleftarrow
         (W_{r},(J_{r},b),E_{r})\longleftarrow\hspace{2cm} \\
        \longleftarrow
        (W_{r+1},(J_{r+1},b),E_{r+1})\longleftarrow\cdots\longleftarrow
        (W_{N},(J_{N},b),E_{N}),
\end{multline}
such  that \( \max{t_{r}^d}=\max{t_{r+1}^d}=\cdots=\max{t_{N-1}^d} \) and either
\begin{enumerate}
        \item  \( \Sing(J_{N},b)=\emptyset \); or
\item  \( \Sing(J_{N},b)\neq\emptyset \);
        and \( \max\word_{N}=0 \); or
        \item  \( \Sing(J_{N},b)\neq\emptyset \),      \( \max\word_{N}>0 \)
      and \( \max{t_{N-1}}>\max{t_{N}} \).
\end{enumerate}
\end{enumerate}

Note that Property 3 says that the function
  \(\max\word^d\) can drop at most finitely
many times, in particular for some index \(N\), either (a) or (b) will hold. We show now why Property (4)
holds, and why the construction of the equivariant resolution \(R_{\mathcal F}\), in Theorem
\ref{Theoremd}, follows essentially from the function $t^i_j$: Roughly speaking,  we {\em attach a new
basic object} to $\max{t_{r}^d}$, say $(W_{r},(J''_{r},b''),E_{r})$, so that $\Max t_{r}^d=
\Sing(J''_{r},b'')$, and a sequence of transformations
\begin{equation}\label{algoextrag}
(W_{r},(J''_{r},b''),E_{r}) \longleftarrow\cdots\longleftarrow (W_{N},(J''_{N},b''),E_{N})
\end{equation}
is naturally associated to sequence  (\ref{nueva2}), so that $\Max t_i^d=\Sing(J''_{r},b'')$
($i=r,r+1,...,N-1$), and such that $\max{t^d}_N<\max{t^d}_{N-1}$ if and only if
$\Sing(J''_{N},b'')=\emptyset$ (if  and only if  sequence (\ref{algoextrag}) is a resolution). The point
is that a resolution of $(W_{r},(J''_{r},b''),E_{r})$ is easy to achieve by induction, essentially as in
\ref{ordenuno}, which we prove in \ref{proposition}. We will do this in two steps, by first attaching a
basic object  to the function $\word^d$ (in the sense described below, see \ref{pattern} for a more
precise description), and finally to the function $t^d$.

\end{parrafo}

\begin{parrafo} {\bf  }
     \label{algomas}
Fix a \(d\)-dimensional basic objects \((W_0,(J_0,b),E_0)\), and let
\begin{equation}\label{algoextra}
(W_{0},(J_{0},b),E_{0}) \longleftarrow\cdots\longleftarrow (W_{r},(J_{r},b),E_{r})
\end{equation}
be  a \(t_i^d\)-permissible sequence (so \(\max \word_i^d\geq \max \word_{i+1}^d\) and \(\max t_i^d\geq
\max t_{i+1}^d\)). For each $i=0,1,\ldots,r$, consider the  factorization
\(\label{masextra}J_i=\cali(H_{l+1})^{a_1}\cdot\ldots\cdot\cali(H_{l+i})^{a_i} \overline{J}_i\) as in
(\ref{extra22}). We will assume that sequence (\ref{algoextra}) has been defined by induction on \(r\),
together with some other added conditions that we will impose below.

\begin{itemize}
\item Case \(\max \word_r>0\).
\end{itemize}
The function \(t_r^d\) was defined only when \(\max \word_r^d>0\),
and in that  case $$\Max t_i^d\subset \Max \word_i^d \subset
V(\overline{J}_i).$$
  Set \(r_0 (\leq r)\) the smallest index such that \(\max
\word_{r_0}= \max \word_{r}\). A basic object
$$(W_{r_0},(J'_{r_0},b'),E_{r_0})$$ can be attached to the value
\(\max \word_{r_0}\), in the sense indicated above (and denoted by
\((W_{r_0},(A_{r_0},c),E_{r_0})\) in (\ref{funchic3}), see Remark
\ref{remarkc3} (e)). Here all the centers \( Y_i \subset \Max
\word_i\),
  and we
will assume that the last \(r-r_0\) steps of (\ref{algoextra}) are defined by a sequence of
transformations
\begin{equation}\label{transword1}
(W_{r_0},(J'_{r_0},b'),E_{r_0}) \longleftarrow (W_{r_0+1},(J'_{r_0+1},b'),E_{r_0+1}) ... \longleftarrow
(W_{r},(J'_{r},b'),E_{r})
\end{equation}
as in the case of  sequence (\ref{algoextrag}) (or see (P1) in \ref{pattern}). Note that
\((W_{r_0},(J'_{r_0},b'),E_{r_0})\) is closer to the setting in
   \ref{ordenuno}, in fact $\ord \cong 1$ (it is within Case 2, see Remark
\ref{remarkc3} (d)). Consider the partition
     \( E_j=E_j^+ \cup E_j^-\), for each index \( r_0 \leq j \leq
r\), and assume that sequence (\ref{transword1}) is such that
\begin{equation*}
     \max n_{r_0} \geq \max
n_{r_0+1}\geq \cdots\geq \max n_r.
\end{equation*}
Since $n_j$ is the second coordinate of $t_j$ we just view it as a function on $\Max \word_j=
\Sing(J'_j,b')$. Set \(r_1 (\geq r_0)\) the smallest index such that \( \max n_{r_1}=\max n_{r}\) and let
$$(W_{r_1},(J''_{r_1},b''),E_{r_1})$$ be the basic object attached to \(\max n_{r_1}\) in
(\ref{funchic2000oh}) (see also Remark \ref{nngbo} (b)).

All centers \( Y_i \subset \Max n_i(= \Max t_i \subset \Max \word_i)\), and we assume that the last
\(r-r_1\) terms of sequence (\ref{transword1}) are defined by
\begin{equation}\label{transword2}
(W_{r_1},(J''_{r_1},b''),E_{r_1}) \longleftarrow (W_{r_1+1},(J''_{r_1+1},b''),E_{r_1+1}) ...
\longleftarrow (W_{r},(J''_{r},b''),E_{r})
\end{equation}
as in the case of  \ref{algoextrag} (or see (P1) in \ref{pattern}).

The point is that  \((W_{r_1},(J''_{r_1},b''),E_{r_1})\cong
(W_{r_1},(J''_{r_1},b''),E_{r_1}^+) \) (see (\ref{congcase2})),
and the right hand term is essentially in the setting of
\ref{ordenuno} (it is within Case 1, see (\ref{proofofcase2})). It
follows now from \ref{caso1prima} that the set
\(F=\Sing(J'',b'')=\Max t^d\) is a disjoint union \( R(1)(F) \cup
F_1\), where \( R(1)(F)\) is a permissible center, and after
blowing up that smooth hypersurface we may assume that
\(R(1)(F)=\emptyset\), in which case
\((W_{r_1},(J''_{r_1},b''),E_{r_1})\) has a \((d-1)\)-dimensional
structure. Theorem \ref{Theoremd} defines a resolution, which in
this case has first center \(R(1)(F)\neq \emptyset \),  and then
proceeds by induction on \(d-1\).

Assume, by induction, that \((I_{d-1}, \geq)\) is defined,
together with functions \(f_i^{d-1}\) which have  the properties
stated in Definition  \ref{AlgResol} for general basic objects of
dimension \(d-1\). Assume also that the resolution that they
provide, is the same resolution as Theorem (d-1) (Theorem
\ref{Theoremd}). We also assume that \(I_{d-1}\) has a biggest
element, say \( \infty_{d-1} \in I_{d-1}\), and that this value is
never reached by any \(f_i^{d-1}\).

Set \( \overline{f}_{r_1}^{d-1}(x)=\infty_{d-1} \) if \( x \in R(1)(F)\), and \(
\overline{f}_{r_1}^{d-1}(x)= {f}_{0}^{d-1}(x)\) if \(x \in F_1\). So if \(R(1) (F) \neq \emptyset\) the
following properties hold:
\begin{gather}\label{ecdesiguno}
     \Max \overline{f}_{r_1}^{d-1}=R(1)(F); \qquad \qquad
    \max \overline{f}_{r_1}^{d-1}> \max \overline{f}_{r_1+1}^{d-1} \\
     \overline{f}_{r_1}^{d-1}(x)= \overline{f}_{r_1+1}^{d-1}(x')
     \qquad\text{if}\qquad x \notin R(1)(F),
\end{gather}
where \(x' \in \Sing(J''_{r_1+1},b'')\) is, in this case, the point naturally identified with \(x\), and
assume that the first transformation in sequence (\ref{transword2}) is defined by blowing up at \(\Max
\overline{f}_{i}^{d-1}=R(1) (F)\) (see Case 1 in \ref{proofofcase1}).

Fix some index \(i\), \(r_1<i \leq r\). For any \(x \in \Max n_i(= \Max t_i \subset \Max \word_i)\) set \(
\overline{f}_{i}^{d-1}(x)={f}_{i-r_1}^{d-1}(x)\), and assume that sequence (\ref{transword2}) is defined
by blowing up at \(\Max \overline{f}_{i}^{d-1}\). Extend sequence (\ref{transword2}) to a resolution by
means of these functions. Finally extend sequence (\ref{algoextra}) by this resolution in the sense of
\ref{pattern} (P1). Since \(\max t^d \) can drop only finitely many times, eventually we come to the
following case.
\medskip

\begin{itemize}
\item Case \(\max \word_r=0\).
\end{itemize}

\noindent Set \(r_0 \leq r\) the smallest index for which \(\max \word_{r_0}=0\). This is the case in
which \(J_{r_0}\) is in the setting of Remark \ref{remarkc3} (b).  This is a simple case, in which no
form of induction is required. It is not hard to formulate the development in the proof of Case \(\word=
0\) in Proposition \ref{case3bo}, by defining a totally ordered set \((\Gamma,\geq)\) and an upper
semi-continuous function \(h_{r_0}: \Sing (J_{r_0},b) \to \Gamma\) so that, as in that proof, \(\Max
h_{r_0}\) is a permissible center (see Section~\ref{SeccMonom}). Furthermore, by setting inductively
functions \(h_i : \Sing(J_i,b) \to \Gamma \),  a resolution
\begin{equation}\label{algoextra0} (W_{r_0},(J_{r_0},b),E_{r_0})
\longleftarrow\cdots\longleftarrow (W_{N},(J_{N},b),E_{N})
\end{equation}
is defined by blowing up at \(\Max h_i \subset \Sing (J_i,b),E_i)\) for \(i=0,1,\ldots,N-1\), and the
sequence is such that
\begin{equation}\label{ecdesigdos}
\max h_j > \max h_{j+1}; \mbox{ \ } \mbox{  } \mbox{ \ } h_{j}^{d}(x)= h_{j+1}^{d}(x') \mbox{
     if } x \notin \Max h_j,
\end{equation}
where \(x' \in \Sing(J_{j+1},b)\) is, in this case, the point naturally identified with \(x\). Assume
that the last \(r-r_0\) steps of sequence (\ref{algoextra})  are the first \(r-r_0\) steps of sequence
(\ref{algoextra0}). Finally extend sequence (\ref{algoextra})
  to a resolution
of length N.
\end{parrafo}

\begin{definition}\label{thefunction}
Set \(T^d=\{\infty\}\sqcup ({\mathbb Q}\times {\mathbb Z}) \sqcup \Gamma\) where this disjoint union is
totally ordered by setting that \(\infty\) is the biggest element, and that \(\alpha < \beta \) if
\(\beta \in ({\mathbb Q}\times {\mathbb Z})\) and \(\alpha \in \Gamma \). We now set \( I_{d}=T^d \times
I_{d-1} \) ordered lexicographically, and define \( f_r^d: \Sing(J_r,b) \to I_{d}\).
\begin{enumerate}
\item[(i)] If  \(\max \word_r> 0\), and \(x \in \Max t_r (\subset
\Sing(J_r,b))\): \( f_r^d(x)=(\max t_r, \overline{f}_r^d(x))\).
\item[(ii)]If  \(\max \word_r=0\): \(
f_r^d(x)= (h_{r-r_0}(x),\infty_{d-1})\).
\item[(iii)] If \(\max \word_r>0\), \(x \notin \Max t_r^d\) and
\(\word_r(x)>0\): There is a smallest index \(r'>r\) in the resolution, such that \(\max
t^d_{r'}=t^d_r(x)\). Note that \(x\) can be identified with a point \(x' \in \Max t^d_{r'}\). Set
\(f_r^d(x)= f_{r'}^d(x')\).
\item[(iv)] If \(\max \word_r>0\) and \( \word_r(x)=0\):
There is a smallest index \(r'>r\) in the resolution, such that \(\max \word_{r'}=0\). Note that \(x\)
can be identified with a point \( x' \in \Sing(J_{r'},b)\) and we set \(f_r^d(x)= f_{r'}^d(x').\)
\end{enumerate}
\end{definition}

\begin{remark} \label{long}{\bf On Definition \ref{AlgResol}.}
\begin{enumerate}
\item[(1)] The value $(\infty,\infty_{d-1})\in I_{d}$ is never reached,
and up to induction on the dimension $d$, choices of centers (namely, \(\Max f^d\)) are defined either as
\(R(1) (\Max t^d)\) in case \( \word^d
>0\), or as \( \Max h^d \) in case \(\word^d=0\). This shows that
  the centers \(\Max f^d\) are regular and permissible.

\item[(2)] Fix, as before, a basic object \({\mathcal B}=(W,(J,b),E)\)
of dimension \(d\), and let
\begin{equation*}
R_{{\mathcal B}}:(W,(J,b),E) \longleftarrow (W_1,(J_1,b),E_1)\longleftarrow...\longleftarrow
(W_M,(J_M,b),E_M)
\end{equation*}
be the the resolution defined by Theorem \ref{Theoremd}. It follows from (\ref{ecdesiguno}) and
(\ref{ecdesigdos}) that \(\max f_j^d> \max f_{j+1}^d\) and that \(f_j^d(x)= f_{j+1}^d(x')\) if \(x \notin
\Max f_j\), where \(x' \in \Sing(J_{j+1},b)\) is the point naturally identified with \(x\).

\item[(3)] One can check from (2)  that given \(\alpha \in I_{d}\),
\begin{equation*}
F_{\alpha}=\{ x \in \Sing(J_r,b) :f_r^d (x)\geq \alpha\}= \cup \pi^s_r(\Max f_s^d),
\end{equation*}
where the union is taken over each index \(s \geq r\) such that \(\max f_s^d \geq \alpha\), and where
\(\pi^s_r: W_s \to W_r\) is the composite morphism. In particular \(F_{\alpha}\) is closed, so the
functions \(f^d_r\) are upper semi-continuous.
\item[(4)] Conditions B (e) and B (f) of Definition \ref{AlgResol}
hold for the functions \(f^d_r\). In fact, by Lemma \ref{propertiesnord} and
  Lemma \ref{casomono}, it follows that both conditions hold for the
functions \(h_i\) and \( t^d_i\) since they hold for the functions \(\ord\) and \(n\) (see also
Proposition \ref{wordngbo}).

\item[(5)] Note that \(I_{d}=T^d \times\ldots\times T^0 \), and that
\((\infty, \infty,\ldots, \infty) \in I_{d}\) is the biggest element.

\item[(6)] It follows now that  these functions have the properties
described in Definition  \ref{AlgResol}; and since all invariants involved rely on the functions \(\ord\)
and \(n\)), Definition \ref{AlgResol} also applies in the setting of general basic objects.
\end{enumerate}
\end{remark}
\begin{example}\label{examplefunction}
Fix \(W= {\mathbb A}^3_k\), \(E= \emptyset\), \(X \subset W\) a smooth subscheme and set \((W,
(\mathcal{I}(X),1),\emptyset)\). Note that \(\Sing(\mathcal{I}(X),1)=X\) and that \(t^3(x)=(1,0)\) for
any \(x \in X\).

If \(X\) is a smooth hypersurface, then \(X=R(1) (\Max(t^3))\) so \(f^3(x)=((1,0), \infty, \infty) \in
I_3\) for  any \(x \in X\).

If \(X\) is a smooth curve, then \(R(1) (\Max(t^3))= \emptyset \), and \(X=R(1) (\Max(t^2))\). In this
case \(f^3(x)=((1,0),(1,0), \infty) \in I_3\) for any \(x \in X\).
\end{example}

\begin{remark}\label{renoinm} It turns out that given  a subscheme
\(X \subset W\), the algorithmic resolution of  \((W,
(\mathcal{I}(X),1),\emptyset)\) depends, to some extent, only on
$X$ (on $\calo_W/I(X)$). In fact, if $J \subset \calo_W$ and
$J'\subset \calo_{W'}$ are such that $ \calo_W /J$ and
$\calo_{W'}/J'$ are isomorphic, and if $$ R_{{\mathcal
B}}:(W,(J,1),\emptyset) \longleftarrow
(W_1,(J_1,1),E_1)\longleftarrow\ldots\longleftarrow
(W_M,(J_M,1),E_M) $$ and $$R_{{\mathcal B'}}:(W',(J',1),\emptyset)
\longleftarrow (W'_1,(J'_1,1),E'_1)\longleftarrow\ldots
\longleftarrow (W'_{M'},(J'_{M'},1),E_{M'}) $$
  are the resolutions defined by Theorem \ref{Theoremd}, then $M=M'$
  and, for each index $i$, there is an identification of $ V(J_i) (\subset W_i)
$ with $ V(J'_i) (\subset W'_i)$. To check this, note first that if dim $W=d$ $
>  $ dim $W'=d'$, then $R_{{\mathcal B}}:(W,(J,1),\emptyset)$ has a
structure of $d-d'$ general basic object. Hence the first $d-d'$ coordinates of $f_i$ are $(1,0)$ at any
point of $V(J_i)$. Use this argument to reduce to the case in which $d=d'$, and then check, by induction
on $i$, that in that setting, there is an identification of the set  $ V(J_i) (\subset W_i) $ with $
V(J'_i) (\subset W'_i)$ such that the functions $f_i$ and $f'_i$ coincide. See  Section
\ref{nonembeddeddes} for more details.
\end{remark}




\section{Bodn\'ar-Schicho's program for resolution of singularities.}
\label{nuevanueva} In \cite{GabSch98}, G. Bodn\'ar and J. Schicho presented a computer program that
produces resolution of basic objects. Full details can be found  in \cite{GabSch99}. Given a basic
object, \( (W,(J,b),E) \) the program provides its resolution,
\begin{equation*}
       (W,(J,b),E)\longleftarrow (W_{1},(J_{1},b),E_{1})
       \longleftarrow\cdots\longleftarrow (W_{r},(J_{r},b),E_{r}),
\end{equation*}
defined in terms of functions \( f_{i}^{d} \) (see Definition \ref{AlgResol}).

\

Among the various problems to solve, the first one is how to encode the variety \( W \).  This is not too
difficult in case \( W \) is of finite type over \( \mathbf{k} \).  They cover \( W \) by open subsets \(
U_{\alpha} \), such that each \( U_{\alpha} \) is a closed subscheme of the affine space \(
\mathbb{A}_{\mathbf{k}}^{m} \), for some \( m \).  In fact the program presents the resolution of the
basic object as a tree of charts.  By looking at the tree one can follow the affine charts of the
sequence of blowing-ups, and the top level of the tree are the charts of \( W_{r} \).

\

 From Definition \ref{thefunction} we see that the functions \(
f_{i}^{d} \) are defined in terms of the function \( \ord \) (introduced in \ref{ordn}), and $n$
(\ref{nnn})
  of several auxiliary basic objects. Given an ideal \( J \), the
function \( \ord(J) \) is computed in terms of the operator \( \Delta \) (\ref{PropOrdDelta}). Finally to
compute \( \Delta(J) \), we need partial derivatives (see \ref{delta}); and to compute derivatives we
need regular systems of parameters at every point.

\

A global section \( f \) of \( \O_{U_{\alpha}} \) is a restriction of a  polynomial \(
f\in\mathbf{k}[X_{1},\ldots,X_{m}] \) modulo the ideal defining  \( U_{\alpha} \). If \(
f_{1},\ldots,f_{n} \) are global sections of \( \O_{U_{\alpha}} \) such that suitable \( d\times d \)
minors of the jacobian matrix do not vanish at any point of \( U_{\alpha} \), then these global sections
give rise to regular systems of parameters, locally at any point of \( U_{\alpha} \).

\begin{parrafo} \label{DefChartBS}
In the computer program given in  \cite{GabSch99} the variety \( W \) is covered by open charts, say \(
U_{\alpha} \), where each chart is presented as follows:
\begin{itemize}
      \item \( U_{\alpha} \) is a closed smooth variety of \(
      \mathbb{A}_{\mathbf{k}}^{m} \) for some \( m \) (depending on \(
      U_{\alpha} \)).  So that \( U_{\alpha} \) is defined by
      polynomials, called \emph{dependencies}, \( D_{1},\ldots,D_{s}\in
      \mathbf{k}[X_{1},\ldots,X_{m}] \).

      \item Polynomials \(
      P_{1},\ldots,P_{d}\in\mathbf{k}[X_{1},\ldots,X_{m}] \) which give
      rise to regular system of parameters for every point of \(
      U_{\alpha} \).

      \item A \( m\times d \) matrix with entries in \(
      \mathbf{k}[X_{1},\ldots,X_{m}] \), such that the restriction of
      each coefficient to $U_{\alpha}$ is the
      partial derivative of \( \dfrac{\partial X_{i}}{\partial P_{j}}
      \).
\end{itemize}
With this structure they are able to compute derivatives in the chart \( U_{\alpha} \).  In particular,
given generators of an ideal \( J\subset\O_{U_{\alpha}} \) they compute generators of \( \Delta(J) \).
To compute the maximal order of \( J \) in the chart \( U_{\alpha} \), one has to find the smallest index
\( b \) such that \( \Delta^{b}(J)=\O_{U_{\alpha}} \), which is accomplished by Gr\"obner basis
computations.
\end{parrafo}

\begin{remark}
Let \( U_{\alpha} \) be a chart in \ref{DefChartBS} and let \( C_{\alpha}=\Max   f_{i}^{d}\) be the
center defined by the algorithm described in Definition  \ref{AlgResol} (i.e. by the functions \(
f_{i}^{d} \)). The center \( C_{\alpha}\subset U_{\alpha} \) is defined by some equations of \(
\mathbf{k}[X_{1},\ldots,X_{m}] \), say \( \mathcal{I}(C_{\alpha})=\left\langle
f_{1},\ldots,f_{r}\right\rangle \).  The program will arrange matters so that \(
\mathcal{I}(C_{\alpha})=\left\langle P_{1},\ldots,P_{r}\right\rangle \) for some \( r<d \).  Using these
equations one may compute the blowing-up with center \( C_{\alpha} \), \( (U_{\alpha})_{1}\longrightarrow
U_{\alpha} \), and cover \( (U_{\alpha})_{1} \) by \( r \) charts as in \ref{DefChartBS}.
\end{remark}

\begin{parrafo}
\textbf{Computing \( \mathcal{I}(C_{\alpha}) \).} We follow the proof of Theorem~\ref{Theoremd} for basic
objects: If the basic object is within Case~1 of Section \ref{proofbo}, then, and as in
  Proposition \ref{ordenuno}, at every chart \(
U_{\alpha} \) they define a smooth hypersurface \( \widetilde{U}_{\alpha}\subset U_{\alpha} \) such that
\begin{equation*}
       (U_{\alpha},(J_{\alpha},b),E_{\alpha})\subset
       (U_{\alpha},(\mathcal{I}(\widetilde{U}_{\alpha}),1),E_{\alpha}).
\end{equation*}
The above inclusion asserts that \( C_{\alpha}\subset\widetilde{U}_{\alpha} \), and hence, if \(
f\in\mathbf{k}[X_{1},\ldots,X_{m}] \) is the equation defining the hypersurface \( \widetilde{U}_{\alpha}
\), then \( f \) is one of the equations defining \( C_{\alpha} \).

Note that the hypersurface \( \widetilde{U}_{\alpha} \) may not be irreducible.  If \(
R(1)(\Sing(J,b))\neq\emptyset \) then \( R(1)(\Sing(J,b)) \) is a smooth hypersurface and a union of some
connected components of \( \widetilde{U}_{\alpha} \). The program of \cite{GabSch99} covers \( U_{\alpha}
\) by open sets, as in \ref{DefChartBS}, of two types:
\begin{enumerate}
       \item  Charts where \( R(1)(\Sing(J,b))=\widetilde{U}_{\alpha}
\).

       \item  Charts where \( R(1)(\Sing(J,b))=\emptyset \).
\end{enumerate}
For charts of the first type the center is \( C_{\alpha}=\widetilde{U}_{\alpha} \).

For charts in as in (2), they provide the computations of the \( (d-1) \)-dimensional basic object
\begin{equation*}
       (\widetilde{U}_{\alpha},
       (\Coeff_{\widetilde{U}_{\alpha}}(J),b!),
       E_{\alpha}\cap\widetilde{U}_{\alpha}),
\end{equation*}
which requires an expression of \( \widetilde{U}_{\alpha} \) as in \ref{DefChartBS}.  If \( f=P_{j} \)
for some index \( j \), then \( \widetilde{U}_{\alpha} \) fulfills \ref{DefChartBS} by adding \( P_{j} \)
to dependencies, and we have \( d-1 \) parameters.

In general \( f \) will not be a parameter, and the program given in \cite{GabSch99} uses an operation
called \emph{cover and exchange}: it covers \( U_{\alpha} \) by charts as in \ref{DefChartBS} such that \(
f=P_{j} \) for some \( j \) at each new chart.  Following this procedure one reduces the computation of
the center for a \( d \)-dimensional basic object to the computation for a \( (d-1) \)-dimensional basic
object.

Assume now that the basic object is within Case~2 or 3 of Section
\ref{proofbo}.  In these cases new basic objects were defined,
always following the general pattern of Remark \ref{pattern}. In
both cases the basic objects are defined in terms of intersections
of other basic objects, and this is a computation (see
\ref{intersection}) which the program can also handle.
\end{parrafo}




\section{Proof of Proposition \ref{ordenuno}.}
\label{proposition} The objective of this section is to prove Proposition \ref{ordenuno}, which is a key
step in the inductive proof of Theorem \ref{Theoremd}:

\

\noindent {\bf Proposition \ref{ordenuno}} {\em Let \({\calB}=(W,(J,b),E)\) be a basic object, assume that
\(E=\emptyset\) and that \( \ord_{ B}\cong 1\). Then:
\begin{enumerate}
\item[a)] There is an open covering \(\{U_{\alpha} \}_{\alpha\in {\Lambda}}\),
    and for each
index \(\alpha\in {\Lambda} \) a closed and smooth hypersurface \(\widetilde{W}_{\alpha}\subset
U_{\alpha}\), such that if \((U_{\alpha},(J_{\alpha},b), \emptyset)\) is the restriction of \((W,
(J,b),\emptyset)\) to \(U_{\alpha}\)
     then
\begin{equation*}
(U_{\alpha},(J_{\alpha},b), \emptyset) \subset (U_{\alpha},(\mathcal{I}(\widetilde{W}_{\alpha}),1),
\emptyset).
\end{equation*}
\item[b)] If \(U_{\alpha} \cap R(1)(\Sing (J,b))=\emptyset \) then
$(U_{\alpha},(J_{\alpha},b), \emptyset)$ has  structure of $d-1-$dimensional basic object, i.e. there is
a $d-1$ basic object \((\widetilde{W},({\mathcal A},e), \emptyset)\) such that
$$(U_{\alpha},(J_{\alpha},b), \emptyset)\cong (\widetilde{W},({\mathcal A},e), \emptyset).$$
\end{enumerate}
}

To prove this proposition we will need some auxiliary  results. Some of them clearly express the role of
the derivatives in our statement, while the others are related to the study of the restriction of basic
objects to smooth hypersurfaces. The proof of Proposition \ref{ordenuno} will be detailed in
\ref{otraolvidada}.

\

\begin{itemize}
\item {\bf The role of the derivatives (an idea of J. Giraud).}
\end{itemize}

\begin{lemma}
\label{invertible} Let \(J\subset \calo_W\) be an invertible sheaf of ideals  and let \( \delta\) be a
globally defined derivation. Then
\begin{equation*}
        (W,(J,b),E) \cap (W,(\delta(J),b-1),E)\cong(W,(J,b),E),
\end{equation*}
or equivalently,
\begin{equation*}
        (W,(\delta(J),b-1),E) \subset (W,(J,b),E).
\end{equation*}
\end{lemma}

\noindent {\em Proof:} From the properties stated in \ref{PropOrdDelta}, we have that
\begin{equation*}
        \Sing(J,b) \cap \Sing(\delta(J),b-1)= \Sing(J,b).
\end{equation*}
We only have to prove that this equality is stable (or preserved) after a transformation,
\begin{equation}\label{rangle1}
(W,(J,b),E) \longleftarrow (W_1,(J_{1},b),E_1).
\end{equation}
It is simple to check this when  the transformations is as in \ref{projections}. We discuss here the case
when the transformation is as in \ref{transformationop}.

Let \(H\) be the exceptional divisor, and let \(\cali(H) \subset \calo_{W_1}\) be the corresponding
invertible sheaf of ideals. We claim that $\cali(H)\cdot\delta$ is an invertible sheaf of derivations on
\(W_1\). To see this we argue locally: Let \(\xi\in W\) be a closed point and choose a regular system of
parameters \(\{x_1,...x_n\}\subset \calo_{W,\xi}\) so that the center of the monoidal transformation is
locally defined by \(\langle x_1,...,x_s\rangle\). Now consider an affine neighborhood \( U \) of \( \xi
\) such that \( x_{1},\ldots,x_{s} \) are global sections of \( \mathcal{O}_{U} \) and such that \( J \)
is generated by a global section, say  \( f \). For simplicity we may assume that \( U=W \). The scheme
\(W_1\) is defined by patching the affine rings $$A_i=\calo_{W}[x_1/x_i,..,x_s/x_i],\mbox{ } \mbox{
}\mbox{ } \mbox{ }i\in \{1,\ldots,s\},$$ and \(\cali(H)=\langle x_i\rangle \) at \(A_i\). Note also that
\begin{equation*}
        \delta\left(\frac{x_{j}}{x_{i}}\right)=
        \frac{\delta(x_{j})}{x_{i}}-
        \frac{x_{j}}{x_{i}}\frac{\delta(x_{i})}{x_{i}},
\end{equation*}
and that $$\cali(H)\delta|_{\mbox{Spec}(A_i)}=x_i.\delta: A_i\to A_i,$$ and hence \(\cali(H)\cdot
\delta\) is an invertible sheaf of derivations on \(W_1\).

\

Now in  $A_i$ consider the factorization \(f=x_i^b g_i\), so that \((\Spec(A_i),(\langle
g_i\rangle,b),E_{i,1})\) is the restriction of \((W_1,(J_1,b),E_1)\) to  \(\mbox{Spec}(A_i)\). Then by
\ref{PropOrdDelta}, the transformation (\ref{rangle1}) induces a transformation
\begin{equation}\label{rangle2}
(W,(\delta(\langle f\rangle),b-1),E) \longleftarrow (W_1,(\delta(\langle f\rangle)_1,b-1),E_1).
\end{equation}

Note here that \(\langle f\rangle_1=\langle g_i\rangle\subset A_i\), and that \(x_i\cdot\delta \) is a
derivation on  \(A_i\). Finally check that
\begin{equation*}
        \frac{\delta(f)}{x_i^{b-1}}=
        \frac{x_i\delta(x_i^b\cdot g_i)}{x_i^{b}}=
        \frac{x_i\delta(x_i^b)}{x_i^b}g_i+
        x_i^b\frac{(x_i\delta)(g_i)}{x_i^b})=
        b\cdot\delta(x_{i})\cdot g_i+ (x_i \delta)(g_i).
\end{equation*}
Using this formula, and the fact that $\delta: \calo_W \to \calo_W $ is a derivation, we also conclude
that $$(\delta(\langle f\rangle))_1 \subset \langle g_{i},\delta(f)/x_i^{b-1}\rangle.$$ The formula also
shows that $\langle g_{i},\delta(f)/x_i^{b-1}\rangle=\langle g_{i},x_i\delta(g_i)\rangle$; and hence,
that globally:
\begin{equation}\label{rangle3}
(\delta(\langle f\rangle))_1\subset\Delta(\langle f\rangle_1),
\end{equation}
and now by \ref{PropOrdDelta} again we have that
\begin{equation*}
        \Sing(\langle f\rangle_1,b) \cap
        \Sing((\delta(\langle f\rangle))_1,b-1)=
        \Sing(\langle f\rangle_1,b).
\end{equation*}
Our argument also shows  that this equality is stable by  any {\em sequence} of transformations.  \qed

\begin{corollary}\label{nose} 
For  any  basic object  \((W,(J,b),E)\) we have that:
\begin{itemize}
\item[1)] $ (W,(J,b),E)\subset (W,(\Delta(J),b-1),E)$ (i.e. $
(W,(J,b),E)\cap (W,(\Delta(J),b-1),E)=(W,(J,b),E)).$

\item[2)] Whenever $$(W,(\Delta(J),b-1),E)\longleftarrow
(W_1,({\Delta(J)}_1,b-1),E_1)$$ is defined in terms of $$(W,(J,b),E)\longleftarrow (W_1,(J_1,b),E_1)$$
(in the sense of (1)), then:
\begin{equation}\label{rangle4}
{\Delta(J)}_1 \subset \Delta(J_1).
\end{equation}
\end{itemize}
\end{corollary}

\noindent{\em Proof:} By  \ref{PropOrdDelta} we have  that $\Sing(J,b) \subset \Sing(\Delta(J),b-1)$,
which proves (1). To prove (2) , note that  if $$(W,(J,b),E)\longleftarrow (W_1,(J_1,b),E_1)$$ is a
projection (as in \ref{projections}), then it is simple to check that \({\Delta(J)}_1 = \Delta(J_1)\).
Now assume that $$(W,(J,b),E)\longleftarrow (W_1,(J_1,b),E_1)$$ is a permissible transformation (as in
\ref{transformationop}).
  If \(J\) is a principal sheaf of ideals, then the proof  follows
directly from  the proof of Lemma \ref{invertible}. If \(J\) is not principal, there is an open
covering of \(W\), \(\{U_i\}_{i=1,\ldots,k}\), such that    \(J\mid_{U_i}=\langle
f_{i,1},...,f_{i,s_i}\rangle\). Then
\begin{equation*}
        (W,(J,b), E)\mid_{U_i}=\cap (W,(\langle
        f_{i_j}\rangle,b),E ),
\end{equation*}
and the statement reduces to the principal case.  \qed

\begin{corollary} \label{corgir}
Let  \((W,(J,b),E)\) be a basic object,  and assume that there is a closed smooth hypersurface
\(\widetilde{W}\subset W\) such that
\begin{equation}
       \mathcal{I}(\widetilde{W}) \subset \Delta^{b-1}(J).
\end{equation}
Then $$(W,(J,b),E)\subset (W,(\cali(\widetilde{W}),1),E).$$
\end{corollary}

\noindent{\em Proof:} By applying Corollary \ref{nose} several times, we have that
\begin{equation*}
        (W,(J,b),E) \subset (W,(\Delta^{b-1}(J),1),E).
\end{equation*}
Now the proof follows from the fact that $\cali(\widetilde{W})\subset \Delta^{b-1}(J)$ and the observation
in Example \ref{yqs}.  \qed

\

\begin{corollary}
\label{cambiada} Let ${\mathcal B}=(W,(J,b),E)$ be a basic object which is within the case \(\ord_B \cong
1\). Then locally at any point $\xi\in Sing(J,b)$, there is  a smooth   hypersurface \(\widetilde{W}\),
such that, identifying $W$ with a suitable neighborhood, \(\widetilde{W}\) is closed in $W$, and
\begin{equation}
\label{forminduc1} (W,(J,b),E) \subset (W, (\mathcal{I}(\widetilde{W}),1),E).
\end{equation}
\end{corollary}

\noindent{\em Proof: } Since ${\mathcal B}$ is within the case $\cong 1$, by Remark \ref{withincase1oh}
(d)  the order of $\Delta^{b-1}(J)$ is at most one at points of $W$. So, locally at each point $\xi \in
V(\Delta^{b-1}(J))=\Sing(J,b)$, there is an element $h_{\xi}\in \Delta^{b-1}(J)_{\xi}$ which defines a
smooth hypersurface \(\widetilde{W}\) at a suitable open neighborhood of \(\xi\). Now the result follows
from Corollary \ref{corgir}. \qed

\

\begin{itemize}
\item {\bf The restriction of basic objects to smooth hypersurfaces.}
\end{itemize}

\begin{lemma}\label{lemRestricZ}
         Let \( (W,(J,b),E) \) be a basic object, and let \(
         Z\subset W \) be a closed smooth subscheme. Assume that
$Z$ and $E$ are as in Definition \ref{DefCoeff}.
         If \( (J\mathcal{O}_{Z})_{\xi}\neq 0 \) for any \( \xi\in
         Z \), then there is an inclusion of basic objects:
         \begin{equation*}
        (W,(J,b),E)\cap (W,(\mathcal{I}(Z),1),E)
        \subset
        (Z,(J\mathcal{O}_{Z},b),E\cap Z).
         \end{equation*}
\end{lemma}

\begin{proof} Set $(W,(J,b),E)=(W_0,(J_0,b),E_0)$ and $Z_0=Z$.
The inclusion $$ \Sing(J_{0},b)\cap\Sing(\mathcal{I}(Z_{0}),1)\subset
             \Sing((J_{0}\mathcal{O}_{Z_{0}}),b)$$ is clear, and it is
clear also that such relation will hold after a
             transformation as in \ref{projections}. So we consider a
sequence of transformation of pairs,
         \begin{equation*}
             (W_0,E_0)\longleftarrow (W_{1},E_{1})
             \longleftarrow\cdots\longleftarrow (W_{k},E_{k}),
         \end{equation*}
         with centers \( C_{i}\subset W_{i} \), \( i=0,\ldots,k-1 \).
Assume
         that this sequence defines sequences of transformations of
         basic objects,

         \begin{equation}
             (W_0,(J_0,b),E_0) \longleftarrow
       (W_{1},(J_{1},b),E_{1})
       \longleftarrow \cdots \longleftarrow
       (W_{k},(J_{k},b),E_{k}),
             \label{eqSeqBas1}
         \end{equation}
         \begin{equation}
             (W_0,(\mathcal{I}(Z_0),1),E_0) \longleftarrow
       (W_{1},(\mathcal{I}(Z_{1}),1),E_{1})
       \longleftarrow \cdots \longleftarrow
       (W_{k},(\mathcal{I}(Z_{k}),1),E_{k}),
             \label{eqSeqBas2}
        \end{equation}
         \begin{multline}
             (Z_0,(J_0\mathcal{O}_{Z_0},b),E_0\cap Z_0)
       \longleftarrow
       (Z_{1},((J_{0}\mathcal{O}_{Z_{0}})_{1},b),E_{1}\cap Z_{1})
       \longleftarrow \cdots   \\
             \cdots \longleftarrow
       (Z_{k},((J_{0}\mathcal{O}_{Z_{0}})_{k},b),E_{k}\cap Z_{k}),
             \label{eqSeqBas3}
        \end{multline}
         and assume, by induction,
         that
         for $i=0,\ldots,k-1$,
         \begin{equation*}
             (C_{i}\subset)
             \Sing(J_{i},b)\cap\Sing(\mathcal{I}(Z_{i}),1)\subset
            \Sing((J_{0}\mathcal{O}_{Z_{0}})_{i},b).
        \end{equation*}
        We want to prove that
         \begin{equation}
       \Sing(J_{k},b)\cap\Sing(\mathcal{I}(Z_{k}),1)\subset
             \Sing((J_{0}\mathcal{O}_{Z_{0}})_{k},b).
             \label{eqInclusSingk}
         \end{equation}

\

         Note that the first transformation of basic objects of
         sequence (\ref{eqSeqBas3}) is
        \begin{equation*}
       (Z_{0},(J_{0}\mathcal{O}_{Z_{0}},b),E_{0}\cap Z_{0})
       \longleftarrow
       (Z_{1},(J_{1}\mathcal{O}_{Z_{1}},b),E_{1}\cap Z_{1}),
         \end{equation*}
         since \( (J_{0}\mathcal{O}_{Z_{0}})_{1}=J_{1}\mathcal{O}_{Z_{1}}
         \).  Now we have
         \begin{equation*}
       \Sing(J_{1},b)\cap\Sing(\mathcal{I}(Z_{1}),1)\subset
       \Sing(J_{1}\mathcal{O}_{Z_{1}},b)=
       \Sing((J_{0}\mathcal{O}_{Z_{0}})_{1},b),
         \end{equation*}
         where the first inclusion is clear, and hence
         inclusion (\ref{eqInclusSingk}), and the lemma follow by
induction on \( k
         \).

\end{proof}

\begin{remark} \label{RemCoeffComplet}
         Let \( Z\subset W \) be a smooth closed subscheme, let \( \xi\in
         Z \) be a closed point, and let
$$\{
         z_{1},\ldots,z_{r},x_{1},\ldots,x_{n} \}$$
be a regular system of
         parameters in \( \mathcal{O}_{W,\xi} \) such that \(
         \mathcal{I}(Z)_{\xi}=(z_{1},\ldots,z_{r}) \).  Consider the
isomorphisms
         \begin{equation*}
             \hat{\mathcal{O}}_{W,\xi}\cong
             k(\xi)[[z_{1},\ldots,z_{r},x_{1},\ldots,x_{n}]],
       \qquad
       \hat{\mathcal{O}}_{Z,\xi}\cong k(\xi)[[x_{1},\ldots,x_{n}]],
         \end{equation*}
         where the right sides are the rings of formal series.
Given \( f\in\mathcal{O}_{W,\xi} \), let \( \hat{f} \) denote
         the image  in \( \hat{\mathcal{O}}_{W,\xi} \).  Assume that
         \begin{equation*}
       \hat{f}=\sum_{i_{1},\ldots,i_{r}=0}^{\infty}
       a_{i_{1},\ldots,i_{r}}z_{1}^{i_{1}}\cdots z_{r}^{i_{r}},
         \end{equation*}
        where each \( a_{i_{1},\ldots,i_{r}}\in
         k(\xi)[[x_{1},\ldots,x_{n}]] \).  Note that
         \begin{equation} \label{EqParcCoeff}
       (i_{1}!\cdots i_{r}!)a_{i_{1},\ldots,i_{r}}=
       \varphi\left(
       \frac{\partial^{i_{1}+\cdots+i_{r}} f}{\partial z_{1}^{i_{1}}\cdots
       \partial z_{r}^{i_{r}}}
       \right),
         \end{equation}
         where \(
         \varphi:k(\xi)[[z_{1},\ldots,z_{r},x_{1},\ldots,x_{n}]]\rightarrow
         k(\xi)[[x_{1},\ldots,x_{n}]] \) is the quotient map induced by
         the inclusion \( Z\subset W \) at \( \xi \). Note also that,
for a fixed integer \( b \),
         \begin{equation} \label{EqCoeffOrd}
       \nu_{\xi}(f)\geq b
       \Longleftrightarrow
       \nu_{\xi}(a_{i_{1},\ldots,i_{r}})\geq b-(i_{1}+\cdots+i_{r}),
\end{equation}
for all $ i_{1},\ldots,i_{r}$ with
           $0\leq i_{1}+\cdots+i_{r}<b$ (here
    the left hand side is the order at \(\calo_{W, \xi}\), and
         the right hand side is the order at \(\calo_{Z, \xi}\)).
\end{remark}

\begin{proposition} \label{PropBasInduc}
        Let \( Z\subset W \) be a closed and smooth subscheme, and let
        \( (W,(J,b),E) \) be a basic object.  Assume that $\{E,Z\}$
        have normal crossings. Then with the same notation as in
        Definition~\ref{DefCoeff},
  if $\Coeff_{Z}(J)_{\xi}\neq 0$  for any $\xi\in Z$,
         there is an equivalence  of basic objects,
\begin{equation*}
(Z,(\Coeff_{Z}(J),b!),E\cap Z) \cong
       (W,(J,b),E)\cap
       (W,(\mathcal{I}(Z),1),E).
         \end{equation*}
\end{proposition}

\begin{proof} Set $(W_{0},(J_{0},b),E_{0})=(W,(J,b),E)$ and  let $Z_0=Z$.
        We will restrict attention to transformations as in
\ref{transformation} and leave
        the case of \ref{projections} to the reader. Consider a
sequence of transformation of pairs:
        \begin{equation*}
      (W_{0},E_{0})\longleftarrow (W_{1},E_{1})
      \longleftarrow\cdots\longleftarrow (W_{k},E_{k})
        \end{equation*}
        with centers \( C_{i}\subset W_{i} \), \( i=0,\ldots,k-1 \).
        Assume that this sequence defines the following three sequences
of transformations:
        \begin{equation*}
      (W_{0},(J_{0},b),E_{0}) \longleftarrow
      (W_{1},(J_{1},b),E_{1}) \longleftarrow \cdots \longleftarrow
      (W_{k},(J_{k},b),E_{k}),
        \end{equation*}
        \begin{equation*}
      (W_{0},(\mathcal{I}(Z_{0}),1),E_{0}) \longleftarrow
      (W_{1},(\mathcal{I}((Z_{0})_{1}),1),E_{1}) \longleftarrow \cdots
      \longleftarrow (W_{k},(\mathcal{I}((Z_{0})_{k}),1),E_{k}),
        \end{equation*}
        \begin{multline*}
      (Z_{0},(\Coeff_{Z_{0}}(J_{0}),b!),E_{0}\cap Z_{0})
      \longleftarrow
      (Z_{1},((\Coeff_{Z_{0}}(J_{0}))_{1},b!),E_{1}\cap Z_{1})
      \longleftarrow \cdots \\
      \cdots \longleftarrow
      (Z_{k},((\Coeff_{Z_{0}}(J_{0}))_{k},b!),E_{k}\cap Z_{k}).
        \end{multline*}
Assume now, by induction on \(k\), that for $i=0,\ldots,k-1$,
        \begin{equation*}
      (C_{i}\subset) \Sing(J_{i},b)\cap Z_{i}=
      \Sing((\Coeff_{Z_{0}}(J_{0}))_{i},b!),
        \end{equation*}
        and recall that
        $\mathcal{I}((Z_{0})_{k})=\mathcal{I}(Z_{k})$, where $Z_k$ is
the strict transform of $Z_0$ (see
       Remark \ref{so} (2)),
        and of course \(\Sing (\mathcal{I}(Z_k),1)=Z_k\).  We
will  prove that
        \begin{equation}
      \Sing(J_{k},b)\cap
      Z_{k}=\Sing((\Coeff_{Z_{0}}(J_{0}))_{k},b!).
      \label{eqInterSingk}
       \end{equation}

       Consider the basic objects \(
        (W_{0},(\Delta^{j}(J_{0}),b-j),E_{0}) \) and \(
        (Z_{0},(\Delta^{j}(J_{0})\mathcal{O}_{Z_{0}},b-j),E_{0}\cap
        Z_{0}) \), for \( j=0,\ldots,b-1\).  From
Definition~\ref{DefCoeff} it follows  that
        \begin{equation*}
     (Z_{0},(\Coeff_{Z_{0}}(J_{0}),b!),E_{0}\cap Z_{0})=
      \bigcap_{j=0}^{b-1}
      (Z_{0},(\Delta^{j}(J_{0})\mathcal{O}_{Z_{0}},b-j),E_{0}\cap
      Z_{0}),
        \end{equation*}
    (see also \ref{intersection});   and by Corollary \ref{nose},
        \begin{equation*}
      (W_{0},(J_{0},b),E_{0}))=
      \bigcap_{j=0}^{b-1}(W_{0},(\Delta^{j}(J_{0}),b-j),E_{0})).
        \end{equation*}
        Finally, Lemma~\ref{lemRestricZ} asserts that
        \begin{equation*}
      (W_{0},(J_{0},b),E_{0}) \cap
      (W_{0},(\mathcal{I}(Z_{0}),1),E_{0}) \subset
      (Z_{0},(\Coeff(J_{0}),b!),E_{0}\cap Z_{0}),
        \end{equation*}
and hence,
\begin{equation*}
      \Sing(J_{k},b)\cap
      Z_{k}\subset\Sing((\Coeff_{Z_{0}}(J_{0}))_{k},b!).
        \end{equation*}
To prove the reverse inclusion we first
       prove the following claim:

\

\noindent {\bf Claim (k):} {\em For any closed point \(
        \xi_{k}\in\Sing((\Coeff_{Z_{0}}(J_{0}))_{k},b!)\subset Z_{k} \),
        there is a regular system of parameters at \( O_{W_{k},\xi_{k}}
        \),
      $\{z_{k,1},\ldots,z_{k,r},x_{k,1},\ldots,x_{k,n}\}$,
        such that
        \begin{enumerate}
      \item[(a)] \( \mathcal{I}(Z_{k})_{\xi_{k}}= \left\langle
      z_{k,1},\ldots,z_{k,r}\right\rangle \).
      \item[(b)] Given  an  isomorphism of complete rings as in
Remark \ref{RemCoeffComplet},
      \begin{gather*}
          \hat{\mathcal{O}}_{W_{k},\xi_{r}}=R_{k}\cong
          k(\xi_{k})[[z_{k,1},\ldots,z_{k,r},x_{k,1},\ldots,x_{k,n}]],
          \\
         \hat{\mathcal{O}}_{Z_{k},\xi_{k}}=\bar{R}_{k}\cong
          k(\xi_{k})[[x_{k,1},\ldots,x_{k,n}]],
      \end{gather*}
      there is a set of
      generators \( \{\hat{f}_{k}^{(\lambda)}\} \) of \( J_{k}R_{k}
      \) such that
      \begin{equation*}
          \hat{f}_{k}^{(\lambda)}=
          \sum_{i_{1},\ldots,i_{r}=0}^{\infty}
          a_{k,i_{1},\ldots,i_{r}}^{(\lambda)}
          z_{k,1}^{i_{1}}\cdots z_{k,r}^{i_{r}},
     \end{equation*}
      as in Remark \ref{RemCoeffComplet}, and
      \begin{equation*}
          \left(a_{k,i_{1},\ldots,i_{r}}^{(\lambda)}
          \right)^{\frac{b!}{b-(i_{1}+\cdots+i_{r})}}\in
          \left(\Coeff_{Z_{0}}(J_{0})\right)_{k}\bar{R}_{k}.
      \end{equation*}
        \end{enumerate}
}

Note  that Claim (k) asserts that
        \begin{equation*}
      \Sing((\Coeff_{Z_{0}}(J_{0}))_{k},b!)\subset
      \Sing(J_{k},b)\cap Z_{k}.
        \end{equation*}
In fact, if $\xi_k \in \Sing((\Coeff_{Z_{0}}(J_{0}))_{k},b!)$, then the order of $
        \left(\Coeff_{Z_{0}}(J_{0})\right)_{k}\bar{R}_{k}$ is $b!$
        at $\bar{R}_{k}$, and hence  by \ref{EqCoeffOrd}
$\xi_k \in \Sing(J_k,b)$.

\

\noindent{\em Proof of  Claim (k):} Claim (0) follows from (\ref{EqParcCoeff}).  Assume that
        Claim (\(k-1\)) holds for some \( k\geq 1 \).  Let \( \xi_{k-1} \)
        be the image of \( \xi_{k} \) in  \( W_{k-1} \), and let
        \begin{equation*}
      \{z_{k-1,1},\ldots,z_{k-1,r},x_{k-1,1},\ldots,x_{k-1,n}\}
        \end{equation*}
        be the regular system of parameters, at \( \xi_{k-1} \), provided
        by Claim (\(k-1\)).  After a finite extension of the base field,
        and a linear change involving only the variables \(
        x_{k-1,1},\ldots,x_{k-1,n} \) in \( R_{k-1} \), we may define a
        regular system of parameters
        \begin{equation*}
      \{z_{k,1},\ldots,z_{k,r},x_{k,1},\ldots,x_{k,n}\}
        \end{equation*}
        with
        \begin{equation*}
      z_{k,i}=\frac{z_{k-1,i}}{x_{k-1,1}}, \quad i=1,\ldots,r,
\qquad x_{k,1}=x_{k-1,1},\quad
      x_{k,i}=\frac{x_{k-1,i}}{x_{k-1,1}}, \qquad i=1,\ldots,n,
        \end{equation*}
       and \( \mathcal{I}(Z_{k})_{\xi_{k}}= \left\langle
        z_{k,1},\ldots,z_{k,r}\right\rangle \).  Set
        \begin{equation*}
      \hat{f}_{k}^{(\lambda)}=
      \frac{\hat{f}_{k-1}^{(\lambda)}}{x_{k-1,1}^{b}}=
      \sum_{i_{1},\ldots,i_{r}=0}^{\infty}
      a_{k,i_{1},\ldots,i_{r}}^{(\lambda)} z_{k,1}^{i_{1}}\cdots
      z_{k,r}^{i_{r}}, \mbox{ } \mbox{ }   \mbox{ and }
\mbox{ } \mbox{ } a_{k,i_{1},\ldots,i_{r}}^{(\lambda)}=
      \frac{a_{k-1,i_{1},\ldots,i_{r}}^{(\lambda)}}{x_{k-1,1}^{b-
      (i_{1}+\cdots+i_{r})}}.
        \end{equation*}
      This settles the first part of (b). To prove that
       $
      \left( a_{k,i_{1},\ldots,i_{r}}^{(\lambda)}
      \right)^{\frac{b!}{b-(i_{1}+\cdots+i_{r})}}\in
      \left(\Coeff_{Z_{0}}(J_{0})\right)_{k}\bar{R}_{k},
      $
        use first (\ref{EqParcCoeff}) to show that $
a_{k,i_{1},\ldots,i_{r}}^{(\lambda)}
      \in \Delta^{i_{1}+\cdots+i_{r}}(J_k)\calo_{Z_k,\xi_k}$, and
then Definition \ref{coeficientes}.
\end{proof}

\

\begin{parrafo}
\label{otraolvidada} {\em Proof of Proposition \ref{ordenuno}:} It follows from Propositions
\ref{cambiada} and \ref{PropBasInduc}. \qed
\end{parrafo}




\section{The monomial case.}
\label{SeccMonom} In this section we will consider  basic objects with the property that the function \(
\word \) is equal to zero at every point. For this special case  we define a resolution quite easily by
means of an upper-semi-continuous function, without making any induction on the dimension of the basic
object. We will see that the resolution, in this case, is purely combinatorial.

\begin{definition} \label{DefMonomBas}
Let \( (W,(J,b),E) \) be a basic object with  \( E=\{H_{1},\ldots,H_{r}\} \). We say the the basic object
is \emph{monomial} if for any point \( \xi\in \Sing(J,b) \) we have that
\begin{equation*}
J_{\xi}=\mathcal{I}(H_{1})^{a_{1}(\xi)}_{\xi}
       \mathcal{I}(H_{2})^{a_{2}(\xi)}_{\xi}\cdots
        \mathcal{I}(H_{r})^{a_{r}(\xi)}_{\xi},
\end{equation*}
where \( a_{i}:H_{i}\cap\Sing(J,b)\longrightarrow\mathbb{Z} \) is a locally constant function.
\end{definition}

Note that, for a monomial basic object, the closed set \( \Sing(J,b) \) is the  union of some of the
irreducible components of intersections of the hypersurfaces \( H_{i} \).  In fact, consider the
intersection \( H_{i_{1}}\cap\cdots\cap H_{i_{p}} \) and let \( C \) be an irreducible component.  Then
the  functions \( a_{i_{1}},\ldots, a_{i_{p}} \) are constant on  \( C \) and \( C \) is included in \(
\Sing(J,b) \) if and only  \( a_{i_{1}}+\cdots+a_{i_{p}}\geq b \) along \( C \).

\begin{definition} \label{DefGamma}
Let \( (W,(J,b),E) \) be a monomial basic object.  With the notation of Definition \ref{DefMonomBas}
define the function:
\begin{gather*}
        h:\Sing(J,b)\longrightarrow
        \Gamma=\mathbb{Z}\times\mathbb{Q}\times\mathbb{Z}^\mathbb{N} \\
        h(\xi)=
      (-p(\xi),\omega(\xi),\ell(\xi)).
\end{gather*}
If \(\xi\in\Sing(J,b)\) define,
\begin{equation} \label{DefGam1}
        p(\xi)=
        \min\left\{
         q\mid\exists i_1,\ldots,i_q,\
         \begin{array}{l}
             a_{i_1}(\xi)+\cdots+a_{i_q}(\xi)\geq b  \\
             \xi\in H_{i_1}\cap\cdots\cap H_{i_q},
         \end{array}\right\}
\end{equation}

\begin{multline} \label{DefGam2}
        \omega(\xi)=
        \max\left\{
        \frac{a_{i_1}(\xi)+\cdots+a_{i_q}(\xi)}{b}\mid \right.\\
        \left. q=p(\xi), \mbox{ } a_{i_1}(\xi)+\cdots+a_{i_q}(\xi)\geq
b,\mbox{ }
        \xi\in H_{i_1}\cap\cdots\cap H_{i_q}\right\}
\end{multline}
and
\begin{multline} \label{DefGam3}
        \ell(\xi)=
        \max\left\{(i_1,\ldots,i_q,0,\ldots)\mid\right. \\
        \left. q=p(\xi), \mbox{ }
\frac{a_{i_1}(\xi)+\cdots+a_{i_q}(\xi)}{b}=\omega(\xi),
        \mbox{ }\xi\in H_{i_1}\cap\cdots\cap H_{i_q}\right\}.
\end{multline}
In the last formula we consider the lexicographical order in \( \mathbb{Z}^{\mathbb{N}} \), and the
convention that \( i_{1}<i_{2}<\cdots<i_{q} \).
\end{definition}

Fix a point \( \xi\in\Sing(J,b) \) and let \( C_{1},\ldots,C_{s} \) be the irreducible components of \(
\Sing(J,b) \) at \( \xi \). The first coordinate of \( h(\xi) \), \( -p(\xi) \),  will indicate the
minimal codimension of \( C_{1},\ldots,C_{s} \).  Denote by \( C'_{1},\ldots,C'_{s'} \) the components
with minimum codimension \( p(\xi) \).  The second coordinate of \( h(\xi) \) is \(
\omega(\xi)=\dfrac{b'}{b} \) where \( b' \) is the maximum of the order of \( J \) along the sets \(
C'_{1},\ldots,C'_{s'} \). Denote by \( C''_{1},\ldots,C''_{s''} \) the components with maximum order. Now
the last coordinate of \( h(\xi) \), \( \ell(\xi) \), corresponds to \( C''_{j} \) for some index \( j \).

So far, fixed a point \( \xi \), with \( p(\xi) \) we have selected the irreducible components of \(
\Sing(J,b) \) at \( \xi \) of highest dimension.  With \( \omega(\xi) \) we have select, among the
previous components, those where the order of \( J \) is maximum.  Finally with \( \ell(\xi) \) we select
a unique component containing \( \xi \).

\begin{parrafo}\label{monoparrafo1}
Now one can check that the function \( h \) is upper-semi-continuous and that the set \( \Max{h} \) is
regular and a union of connected components of the regular scheme \( H_{i_{1}}\cap\cdots\cap
H_{i_{p_{0}}} \) if \( \max{h}=(-p_{0},w_{0},(i_{1},\ldots,i_{p_{0}},0,\ldots)) \).  It is clear that \(
\Max{h} \) is a permissible center for  the basic object \( (W,(J,b),E) \). Let
\begin{equation*}
       (W,(J,b),E)\stackrel{\Pi}{\longleftarrow}(W_{1},(J_{1},b),E_{1})
\end{equation*}
be the transformation with center \( \Max{h} \), and let  \( E_{1}=\{H_{1},\ldots,H_{r},H_{r+1}\} \),
where, by abuse of notation, \( H_{i} \) is the strict transform of \( H_{i} \), for \( i=1,\ldots,r,\)
and \( H_{r+1} \) is the exceptional divisor of the transformation \( \Pi \).  The basic object \(
(W_{1},(J_{1},b),E_{1}) \) is also monomial, in fact for \( \xi\in\Sing(J_{1},b) \) we have
\begin{equation}
        J_{\xi}=
       \mathcal{I}(H_{1})^{a'_{1}(\xi)}_{\xi}\cdots
        \mathcal{I}(H_{r})^{a'_{r}(\xi)}_{\xi}
        \mathcal{I}(H_{r+1})^{a'_{r+1}(\xi)}_{\xi},
        \label{eq:ExpJ1}
\end{equation}
where the functions \( a'_{i} \) are given by:
\begin{equation}
        \begin{array}{lll}
            a'_{i}(\xi)=a_{i}(\Pi(\xi)) & \forall \xi\in H_{i} &
             \text{and}\ i=1,\ldots,r; \\
       a'_{r+1}(\xi)=a_{i_{1}}(\Pi(\xi))+\cdots+a_{i_{p_{0}}}(\Pi(\xi))-b &
       \forall \xi\in H_{r+1}. \\
         \end{array}
          \label{eq:Expaprim}
  \end{equation}

As in Definition \ref{DefGamma},  a function \( h_{1} \) has been associated to the basic object \(
(W_{1},(J_{1},b),E_{1}) \), and  one can check that the maximum value has dropped:
  \begin{equation*}
  \max{h}>\max h_{1}.
  \end{equation*}

In fact, for any point \( \xi\in\Sing(J_{1},b) \):
\begin{equation}
        \begin{array}{lll}
        h_{1}(\xi)=h(\Pi(\xi)) & \text{if} &
        \Pi(\xi)\not\in\Max h  \\
       h_{1}(\xi)<h(\Pi(\xi)) & \text{if} & \Pi(\xi)\in\Max h.
         \end{array}
          \label{eq:IgualMenor}
  \end{equation}
  It is not hard to check that
this function  \( h \) defines
  a resolution of the basic object. To illustrate
this fact and the statement in (\ref{eq:IgualMenor}),
  consider the following example:
  \end{parrafo}

\begin{example}\label{exaplemono}
  Let \( (W,(J,b),E) \) where   $W$ now denotes the real analytic space
${\mathbb R}^4$, $H_1,\ldots,H_4$ are  the coordinate hyperplanes,  and
  $$  E=\{H_{1},H_{2},H_{3},H_{4}\}, \qquad b=9, \qquad
          J=\mathcal{I}(H_{1})^{6}\mathcal{I}(H_{2})^{4}
       \mathcal{I}(H_{3})^{2}\mathcal{I}(H_{4})^{2}.$$

Note that the intersection \( H_{1}\cap H_{2}\cap H_{3}\cap H_{4}
  \) is a closed point. The singular locus of the basic object is:
  \begin{equation*}
          \Sing(J,9)=\left(H_{1}\cap H_{2}\right)\cup
          \left(H_{1}\cap H_{3}\cap H_{4}\right),
  \end{equation*}
  and the function \( h \) is given by
  \begin{equation*}
          h(\xi)=\left\{
         \begin{array}{lcl}
              (-2,\frac{10}{9},(1,2,0,0)) & \text{if} & \xi\in
H_{1}\cap H_{2}  \\
              (-3,\frac{10}{9},(1,3,4,0)) & \text{if} &
       \xi\in \left(H_{1}\cap H_{3}\cap H_{4}\right)\setminus
        \left(H_{1}\cap H_{2}\right),
          \end{array}
          \right.
\end{equation*}
  so that \( \max h=(-2,\frac{10}{9},(1,2,0,0)) \) and \( \Max
  h=H_{1}\cap H_{2} \). Note that the two irreducible components of
  \( \Sing(J,9) \) are of different dimension and that the function \( p
  \) takes different values along each component. Consider the
  transformation with center \( H_{1}\cap H_{2} \):
  \begin{equation*}
         (W,(J,9),E)\longleftarrow(W_{1},(J_{1},9),E_{1}),
  \end{equation*}
  where \( E_{1}=\{H_{1},H_{2},H_{3},H_{4},H_{5}\} \) and
  \begin{equation*}
          J_{1}=\mathcal{I}(H_{1})^{6}\mathcal{I}(H_{2})^{4}
          \mathcal{I}(H_{3})^{2}\mathcal{I}(H_{4})^{2} \mathcal{I}(H_{5}).
  \end{equation*}
  The singular locus of the new basic object is
$$
  \Sing(J_{1},9)=\left(H_{1}\cap H_{3}\cap H_{4}\right)\cup
         \left(H_{1}\cap H_{3}\cap H_{5}\right)\cup
          \left(H_{1}\cap H_{4}\cap H_{5}\right)\cup
          \left(H_{2}\cap H_{3}\cap H_{4}\cap H_{5}\right).
$$
  The function \( h_{1} \) corresponding to the basic object \(
(W_{1},(J_{1},9),E_{1}) \) is given as follows:
  \begin{itemize}
         \item  If \( \xi\in H_{1}\cap H_{3}\cap H_{4} \) then
          \begin{equation*}
             h_{1}(\xi)=(-3,\frac{10}{9},(1,3,4,0)).
         \end{equation*}
        \item  If \( \xi\in H_{1}\cap H_{4}\cap H_{5} \) and
          \( \xi\not\in H_{1}\cap H_{3}\cap H_{4} \) then
          \begin{equation*}
              h_{1}(\xi)=(-3,\frac{9}{9},(1,4,5,0)).
          \end{equation*}
         \item  If \( \xi\in H_{1}\cap H_{3}\cap H_{5} \) and
          \( \xi\not\in \left(H_{1}\cap H_{3}\cap H_{4}\right)\cup
         \left(H_{1}\cap H_{4}\cap H_{5}\right) \) then
          \begin{equation*}
              h_{1}(\xi)=(-3,\frac{9}{9},(1,3,5,0)).
          \end{equation*}
          \item  If \( \xi\in H_{2}\cap H_{3}\cap H_{4}\cap H_{5} \) then
         \begin{equation*}
              h_{1}(\xi)=(-4,\frac{9}{9},(2,3,4,5)).
          \end{equation*}
  \end{itemize}
  Here \( \Sing(J_{1},9) \) has four irreducible components and the
  function \( h_{1} \) takes different values along each component.
Then  \( \max h_{1}=(-3,\frac{10}{9},(1,3,4,0)) \)
  and  \( \Max h_{1}=H_{1}\cap H_{3}\cap H_{4} \). Note that \( \max
  h>\max h_{1} \). We construct the transformation of basic objects
  with center \( \Max h_{1} \):
  \begin{equation*}
          (W_{1},(J_{1},9),E_{1})\longleftarrow (W_{2},(J_{2},9),E_{2}).
  \end{equation*}
  Now \( E_{2}=\{H_{1},H_{2},H_{3},H_{4},H_{5},H_{6}\} \) and
  \begin{equation*}
          J_{2}=\mathcal{I}(H_{1})^{6}\mathcal{I}(H_{2})^{4}
          \mathcal{I}(H_{3})^{2}\mathcal{I}(H_{4})^{2}
          \mathcal{I}(H_{5})\mathcal{I}(H_{6}).
  \end{equation*}
  The singular locus of \( (W_{2},(J_{2},9),E_{2}) \) is
$$ \Sing(J_{2},9)=
          \left(H_{1}\cap H_{3}\cap H_{5}\right)\cup
          \left(H_{1}\cap H_{4}\cap H_{5}\right)\cup $$
$$\cup\left(H_{1}\cap H_{3}\cap H_{6}\right)\cup
         \left(H_{1}\cap H_{4}\cap H_{6}\right)\cup
          \left(H_{2}\cap H_{3}\cap H_{4}\cap H_{5}\right)\cup
          \left(H_{2}\cap H_{3}\cap H_{4}\cap H_{6}\right).
$$
  There are four irreducible components of dimension one and two
  components of dimension zero.  The function \( h_{2} \) takes
  different values along each component and one concludes that \(
  \max h_{2}=(-3,\frac{9}{9},(1,4,6,0)) \) and that \( \Max
  h_{2}=H_{1}\cap H_{4}\cap H_{6} \).  Consider the transformation
  with center \( \Max h_{2} \)
  \begin{equation*}
          (W_{2},(J_{2},9),E_{2})\longleftarrow (W_{3},(J_{3},9),E_{3}),
  \end{equation*}
  where \( E_{3}=\{H_{1},H_{2},H_{3},H_{4},H_{5},H_{6},H_{7}\} \)
  and
  \begin{equation*}
          J_{3}=\mathcal{I}(H_{1})^{6}\mathcal{I}(H_{2})^{4}
          \mathcal{I}(H_{3})^{2}\mathcal{I}(H_{4})^{2}
         \mathcal{I}(H_{5})\mathcal{I}(H_{6}).
  \end{equation*}
  Now the singular locus is
$$ \Sing(J_{3},9)=
          \left(H_{1}\cap H_{3}\cap H_{5}\right)\cup
          \left(H_{1}\cap H_{4}\cap H_{5}\right)\cup $$
$$\cup\left(H_{1}\cap H_{3}\cap H_{6}\right)\cup
          \left(H_{2}\cap H_{3}\cap H_{4}\cap H_{5}\right)\cup
          \left(H_{2}\cap H_{3}\cap H_{4}\cap H_{6}\right),
$$
  and we have that \( \max h_{3}=(-3,\frac{9}{9},(1,4,5,0)) \)
  and \( \Max h_{3}=H_{1}\cap H_{4}\cap H_{5} \). After five
  transformations of basic objects we achieve a resolution.
  \end{example}


\section{On Hironaka's trick.}
  \label{Hironaka}
  The purpose of this section is to state and prove
  Proposition \ref{HiroTrick}. A consequence of this
  result is that the function $\ord$ introduced in   \ref{ordn}
    is well defined
  for general basic objects, and moreover,
    equivariant (see Remark \ref{Hiro1}).

  \begin{proposition} \label{HiroTrick}
  \cite{AHV1} Let \( (\mathcal{F},(W,E)) \) be a \( d \)-dimensional
  general basic object. Assume that there is   an open covering  \(
  \{U_{\alpha}\}_{\alpha\in \Lambda} \)
     of $W$ as in Definition
\ref{defbo} consisting of  the  single open set   \( W \). Assume
also that we have two \( d
  \)-dimensional basic objects \( (W',(B',d'),E') \) and \(
  (W'',(B'',d''),E'') \) defining \( (\mathcal{F},(W,E)) \), (i.e.
  both basic objects  describe the closed sets associated to \(
  (\mathcal{F},(W,E)) \)).  If \(
  x\in\Sing(B',d')=\Sing(B'',d'') \) then
  \begin{equation*}
          \frac{\Fnu_{B'}(x)}{d'}=\frac{\Fnu_{B''}(x)}{d''}.
  \end{equation*}
  \end{proposition}
  \begin{proof}
  Set \( \omega'=\Fnu_{B'}(x) \) and \( \omega''=\Fnu_{B''}(x) \).
  We shall prove the proposition  by constructing a number of \( x
  \)-extendable sequences of transformations of general basic
  objects with suitable properties. Let
  \begin{equation}
          (\mathcal{F},(W,E))\stackrel{\Pi_{0}}{\longleftarrow}
          (\mathcal{F}_{0},(W_{0},E_{0}))\stackrel{\Pi_{1}}{\longleftarrow}
          (\mathcal{F}_{1},(W_{1},E_{1}))\stackrel{\Pi_{2}}{\longleftarrow}
          \cdots\stackrel{\Pi_{k}}{\longleftarrow}
          (\mathcal{F}_{k},(W_{k},E_{k}))
         \label{SeqGBOHiro}
  \end{equation}
  be a sequence of transformations of general basic objects. Note
  that   sequence (\ref{SeqGBOHiro}) defines sequences of
  transformations of basic objects
  \begin{multline}
          (W',(B',d'),E')\stackrel{\Pi'_{0}}{\longleftarrow}
          (W'_{0},(B'_{0},d'),E'_{0})
          \stackrel{\Pi'_{1}}{\longleftarrow} (W'_{1},(B'_{1},d'),E'_{1})
          \stackrel{\Pi'_{2}}{\longleftarrow} \cdots \\
          \cdots
         \stackrel{\Pi'_{k}}{\longleftarrow} (W'_{k},(B'_{k},d'),E'_{k}),
          \label{SeqBasHiro1}
  \end{multline}
  and
  \begin{multline}
          (W'',(B'',d''),E'')\stackrel{\Pi''_{0}}{\longleftarrow}
          (W''_{0},(B''_{0},d''),E''_{0})
          \stackrel{\Pi''_{1}}{\longleftarrow} (W''_{1},(B''_{1},d''),E''_{1})
          \stackrel{\Pi''_{2}}{\longleftarrow} \cdots \\
         \cdots
          \stackrel{\Pi''_{k}}{\longleftarrow} (W''_{k},(B''_{k},d''),E''_{k}).
          \label{SeqBasHiro2}
  \end{multline}
  The first transformation \( \Pi_{0} \) of (\ref{SeqGBOHiro}) is a
  projection (as allowed in \ref{defbo}), so the first transformations of
  (\ref{SeqBasHiro1}) and (\ref{SeqBasHiro2}) are  projections too.
  All the other transformation will be permissible transformations (as in
  \ref{transformation}). For each index \( k \),   sequence
  (\ref{SeqGBOHiro}) is defined as follows:
  \begin{enumerate}
          \item Identify \( L_{0}=\Pi_{0}^{-1}(x) \) with \(
          \mathbb{A}^{1}_k \) and set \( x_{0} = 0 \in L_0 \).
          Note that \( L_{0} \subset F_{0} \), the singular locus of \(
          (\mathcal{F}_{0},(W_{0},E_{0})) \).
         \item Given an index \( s \geq 0 \), a line \( L_s \subset F_s \)
          and a point \( x_s \in L_s \), define the transformation \(
          \Pi_{s+1} \) with center \( x_s \).  Now set:
          \begin{description}
       \item[i] \(L_{s+1}\subset F_{s+1}\) as the strict transform
       of \( L_s \).
        \item[ii] \( H_{s+1}(\in E_{s+1}) \) as the exceptional
        locus of \( \Pi_{s+1} \);
       \item[iii] \( x_{s+1}=H_{s+1}\cap L_{s+1} \).
          \end{description}
  \end{enumerate}
  Note that (1) together with (2) provides us with a  rule for  constructing
   a sequence as (\ref{SeqGBOHiro}) with length \( s \) for any \( s
  \geq 1 \). By construction the sequence (\ref{SeqGBOHiro}) is \( x
  \)-extendable (see \ref{ExtendTrans}).
  Note that \( L_{s}\subset F_{s} \) for any \( s \), so in
  particular \( x_{s}\in F_{s} \), and by assumption:
  \begin{equation*}
         x_{s}\in\Sing(B'_{s},d')=\Sing(B''_{s},d'')
          \qquad \forall s\geq 0.
  \end{equation*}
  There are expressions as in (\ref{expresionj}),
  \begin{equation} \label{ExprHiro}
        (B'_{s})_{ x_{s}}=
         \mathcal{I}(H'_{s})_{ x_{s}}^{a'_{s}}(\overline{B'}_{s})_{ x_{s}}
         \qquad
         (B''_{s})_{ x_{s}}=
         \mathcal{I}(H''_{s})_{ x_{s}}^{a''_{s}}(\overline{B''}_{s})_{ x_{s}}.
  \end{equation}
  Note that here  \( H'_{s}=H_{s}\cap W'_{s} \) and \(
  H''_{s}=H_{s}\cap W''_{s} \). On may check, by induction on \( s
  \), that
  \begin{equation*}
          a'_s = s(\omega' - d') \qquad a''_{s}=s(\omega''-d'').
  \end{equation*}
  Since the first term of all sequences is a projection, then for \(
  s \geq 1 \), \( \dim(W'_{s})=\dim(W''_{s})=d+1 \). It follows that
\begin{equation*}
          \begin{split}
             \dim(F_s\cap H_s)=d
       & \Leftrightarrow a'_{s}=s(\omega'-d')\geq d'  \\
             & \Leftrightarrow a''_{s}=s(\omega''-d'')\geq d''.
          \end{split}
  \end{equation*}
  Note that \( \dim{H'_{s}}=\dim{H''_{s}}=d \), so if \(
  \dim(F_s\cap H_s)=d \) then \( F_{s}\cap H_{s}=H'_{s}=H''_{s} \).
  Furthermore if \( \dim(F_s\cap H_s)=d \), then \( F_s\cap H_s \)
  is a permissible center for the general basic object. At this
  point we start a new stage of the process by choosing
    \( F_s \cap H_s \) as a center of a
  transformation \( \Pi_{s+1} \). It turns out that in the
  expressions of (\ref{ExprHiro}),
  \begin{equation*}
          a'_{s}=s(\omega'-d')-d', \qquad a''_{s}=s(\omega''-d'')-d''.
  \end{equation*}
  Fix the index \( s \) and set, if possible, the center of
  transformations \( \Pi_{s+j} \) to be \( F_{s+j}\cap H_{s+j} \),
for \( j\geq 0 \). Note that
\begin{equation*}
          \begin{split}
              \dim(F_{s+j}\cap H_{s+j})=d
       & \Leftrightarrow a'_{s+j}=s(\omega'-d')-jd'\geq d'  \\
             & \Leftrightarrow a''_{s+j}=s(\omega''-d'')-jd''\geq d''.
          \end{split}
\end{equation*}
And we conclude that
  \begin{equation*}
          \begin{split}
              \dim(F_{s+j}\cap H_{s+j})=d\ \text{(in which case is a permissible
              center)}
        & \Leftrightarrow j\leq \ell'_{s}  \\
              & \Leftrightarrow j\leq \ell''_{s},
         \end{split}
  \end{equation*}
where
\begin{equation*}
          \ell'_{s}=\left\lfloor\frac{s(\omega'-d')}{d'}\right\rfloor
        \qquad
          \ell''_{s}=\left\lfloor\frac{s(\omega''-d'')}{d''}\right\rfloor
  \end{equation*}
  and \( \lfloor\cdot\rfloor \) denotes the integer part.

Finally note that
  \begin{gather*}
          \frac{\Fnu_{B'}(x)}{d'}=\frac{w'}{d'}=
        \lim_{s\to\infty}\frac{1}{s}\ell'_{s}+1,  \\
          \frac{\Fnu_{B''}(x)}{d''}=\frac{w''}{d''}=
          \lim_{s\to\infty}\frac{1}{s}\ell''_{s}+1,
  \end{gather*}
  and we have expressed the rational numbers \(
  \dfrac{\Fnu_{B'}(x)}{d'} \) and \( \dfrac{\Fnu_{B''}(x)}{d''} \)
   in terms of permissible sequences of the general basic object  \(
  (\mathcal{F},(W,E)) \), and all sequences are \( x
  \)-extendable.  Hence \(
  \dfrac{\Fnu_{B'}(x)}{d'}=\dfrac{\Fnu_{B''}(x)}{d''} \).
  \end{proof}

\begin{remark}\label{Hiro1}
  It  follows from Proposition~\ref{HiroTrick} that the function \(
  \ord \) is well defined for general basic objects, and note also
  that the proof of Proposition \ref{HiroTrick} implies the
  equivariance of this  function:  Given an isomorphism of \( d
  \)-dimensional general basic objects
    \( \Theta:(\mathcal{F},(W,E))\to
  (\mathcal{F}',(W',E')) \) (see Definition \ref{isoGBO}) and a
  point \( x\in\Sing(\mathcal{F}) \), set \(
  x'=\Theta(x)\in\Sing(\mathcal{F}') \).  By Definition \ref{isoGBO}
  (ii), the isomorphism \( \Theta \) defines a bijection between \(
  x \)-extendable sequences and \( x' \)-extendable sequences, say
  \( \alpha_{\Theta}:C_{x}(\mathcal{F})\to C_{x'}(\mathcal{F}') \).
  Note that the dimension \( n \) is fixed and the rational number
  \( \ord_{\mathcal{F}}(x) \) (resp. \( \ord_{\mathcal{F}'}(x') \))
    is expressed in terms of sequences of \( C_{x}(\mathcal{F}) \)
  (resp. \( C_{x'}(\mathcal{F'}) \)). We conclude that \(
  \ord_{\mathcal{F}}(x)=\ord_{\mathcal{F}'}(x') \).
  \end{remark}


\begin{thebibliography}{W}



\bibitem{Ab} S.S. Abhyankar, {\em Good points of a
  Hypersurface}, Adv. in Math. {68} (1988) 87-256.

\bibitem{AbrJong97} D. Abramovich and A.J. de Jong,
  `Smoothness, semistability and toroidal geometry', {\em  Journal of
  Algebraic Geometry} { 6} (1997) 789-801.





\bibitem{AbrWang97} D. Abramovich and J. Wang,  `Equivariant
  resolution of singularities in characteristic 0', {\em Mathematical
  Research Letters}  4 (1997) 427-433.





\bibitem{AHV1} J.M. Aroca, H. Hironaka and J.L. Vicente,
  `The theory of maximal contact', {\em Mem. Mat. Ins.
  Jorge Juan (Madrid)}  29 (1975).



\bibitem{Artin} M.~Artin.
  \newblock  `Algebraic approximation of structures over complete
  local rings', {\em Pub. Math. I.H.E.S.} 36    (1969) 23-58.

\bibitem{BM91} E.~Bierstone and P.~Milman.
\newblock A simple constructive proof of canonical resolution of
singularities,
\newblock Effective Methods in Algebraic Geometry,
Prog. in Math. 94, pp. 11--30, Birkhauser 1991.


\bibitem{BM97} E. Bierstone and P. Milman,  `Canonical
  desingularization in characteristic zero by blowing-up the maxima
  strata of a local invariant', {\em Inv. Math.}  128 (1997) 207-302.


\bibitem{BM02} E. Bierstone and P. Milman, `Desingularization algorithms I. Role of
exceptional divisors'. Preprint: ArXiv math.AG/0207098.

\bibitem{BodEncSch} G.~Bodn\'{a}r, S.~Encinas and J.~Schicho,
`Localizing Villamayor's Algorithm for Resolution of Singularities'. Preprint.



\bibitem{GabSch98} G.~Bodn\'{a}r and J.~Schicho
   `A Computer Program for the Resolution of Singularities',
  {\em Resolution of singularities. A research
  book in tribute of Oscar Zariski}
  (eds H. Hauser, J. Lipman, F. Oort, A. Quir\'os), Progr. Math. 181
  (Birkh\"auser, Basel, 2000) pp.
    231-238.





\bibitem{GabSch99}
   G.~Bodn\'{a}r and J.~Schicho,
  `Automated resolution of singularities for hypersurfaces',
{\em J. Symbolic Comput.}  (4) 30 (2000) 401-428.



\bibitem{BogPan96} F. Bogomolov and T. Pantev,  `Weak
  Hironaka Theorem', {\em Mathematical Research Letters}  3 (1996)
  299-307.





\bibitem{BV1} A. Bravo and O. Villamayor,
  `Strengthening a Theorem of Embedded Desingularization,'
  {\em Math. Res. Letters} { 8}  (2001) 1-11.



\bibitem{BV2} A. Bravo and O. Villamayor,
  `A Strengthening of resolution of singularities in characteristic zero',
  {\em Proc. London Math. Soc.} 86:2 (2003) 327-357.


\bibitem{EH} S. Encinas and H. Hauser, `Strong Resolution of
   Singularities'. {\em Comentarii Mathematici Helvetici},
 77,  4 (2002) 821-845.



\bibitem{EncinasNobileVillamayor} S. Encinas, A. Nobile and O. Villamayor, `On
  algorithmic Equiresolution and
  stratification of Hilbert schemes',
 {\em Proc. London Math. Soc.} 86:3 (2003) 607-648.



\bibitem{EncVil98} S. Encinas and O. Villamayor, `Good
  points and constructive resolution of singularities', {\em Acta Math.}
    181:1  (1998) 109-158.



\bibitem{EncVil97:Tirol}
  S.~Encinas and O.~Villamayor,
  \newblock `A Course on Constructive Desingularization and
  Equivariance',
  \newblock  {\em Resolution of Singularities. A research textbook
  in tribute to Oscar Zariski} (eds H. Hauser, J.
  Lipman, F. Oort, A. Quir\'os),  Progr. in Math.  181
  (Birkh\"auser, Basel, 2000) pp.  147-227.



\bibitem{EncVil99} S. Encinas and O. Villamayor,  `A new
  theorem of desingularization over fields of characteristic zero'.
  Preprint. arXiv:math.AG/0101208.



\bibitem{Gir} J. Giraud, `Sur la theorie du contact maximal',
      {\em Math. Zeit.}, 137  (1972), 285-310.






\bibitem{GrothendieckRaynaud1971}
  A.~Grothendieck and M.~Raynaud.
  \newblock `S\'eminaire de G\'eometrie Alg\'ebrique de Bois Marie 1960/61.
        Rev\^etements Etales et Groupe Fondamental', {\em
  Lecture Notes
        in Mathematics,}  224
  \newblock (Springer-Verlag, Berlin, 1971).

\bibitem{Hauser} H. Hauser, `The Hironaka Theorem on Resolution of
Singularities'. To appear in Bull. Amer. Math. Soc..

\bibitem{Hir64} H.~Hironaka, `Resolution of singularities
  of an algebraic variety over a field of characteristic zero I-II',
  {\em Ann. Math.},  79 (1964) 109--326.



\bibitem{Hironaka77} H. Hironaka, `Idealistic exponent of a
  singularity', {\em Algebraic Geometry},  The John Hopkins centennial
  lectures, Baltimore, John Hopkins University Press (1977), 52-125.


\bibitem{laz} R. Lazarsfeld, Positivity in Algebraic Geometry. In preparation.


\bibitem{Lipman}  J.~Lipman, `Introduction to Resolution of Singularities',
Proc. Symp. in Pure Math. 29 (1975) pp. 187-230.

\bibitem{LipmanTirol}  J.~Lipman,
  `Equisingularity and Simultaneous Resolution of Singularities',
    {\it Resolution of Singularities. A research textbook in tribute to
  Oscar Zariski}, (eds H. Hauser, J. Lipman, F.
  Oort and A. Quir\'os),
  Progr. in Math. 181, (Birkh\"auser,
  Basel, 2000) pp. 485-505.



\bibitem{Matsuki}  K. Matsuki,
  `Notes on the inductive algorithm of resolution of singularities by
S. Encinas and O.
  Villamayor'. Preprint. arXiv:math.math.AG/0103120



\bibitem{Matsumura1980} H.~Matsumura.
  `Commutative algebra',  {\em Mathematics
     Lecture Note Series, 56},
  Benjamin/Cummings Publishing Company, Inc., 2nd ed. edition, 1980.

\bibitem{Moh92} T. T. Moh.
\newblock  Canonical uniformization of hypersurface singularities of
characteristic zero,
\newblock Communication of Algebra 20 (1992)
pp. 3207--3251.


  \bibitem{Mumford} D.~Mumford,
  `Geometric Invariant Theory',
  {\em Ergebnisse der Mathematik,} 34, Springer-Verlag, Berlin (1965).

\bibitem{Oda1987} T.~Oda,
  \newblock `Infinitely very near singular points',
  \newblock  {\em Complex analytic singularities},  Adv.
  Studies in Pure Math. 8  (North-Holland, 1987) pp. 363--404.

\bibitem{Stromme} S.~A.~Stromme,
  `Elementary introduction to representable functors and
  Hilbert Schemes',
  {\em Parameter Spaces,} Banach Center Publications,  36.
  Institute of Mathematics. Polish Academy of Science
  Warszawa 1996.





\bibitem{Villa89} O.~Villamayor, `Constructiveness of
  Hironaka's resolution', {\em
  Ann. Scient. Ec. Norm. Sup. \(4^{\rm{e}}\)} serie  22 (1989) 1-32.

\bibitem{Villa92} O.~Villamayor, Patching local
uniformizations, {\em Ann. Scient. Ec. Norm. Sup.},  25 (1992), 629-677.

\bibitem{Villa02} O. Villamayor, `Introduction to the
  algorithm of resolution',  Progress in Mathematics 134, (Birkh\"auser Verlag
  Basel, Switzerland).
  \end{thebibliography}
  \end{document}